\documentclass[3p,times]{elsarticle}

\usepackage{ecrc}

\volume{00}

\firstpage{1}

\journalname{Computer Methods in Applied Mechanics and Engineering}

\runauth{S. Kang and T. Bui-Thanh}

\jid{cmame}

\usepackage{amssymb}
\usepackage{amsmath}
\usepackage{amsthm}

\usepackage{float}
\usepackage[pdfborder={0 0 0}, colorlinks=true, linkcolor=blue, citecolor=red]{hyperref}
\hypersetup{pdfborder=0 0 0}
\usepackage{algorithm}
\usepackage{algorithmic}
\usepackage{mathrsfs}
\usepackage{subfigure}
\usepackage{multirow}

\usepackage{leftidx}
\usepackage[latin1]{inputenc}
\usepackage{ulem}

\usepackage{IEEEtrantools}

\usepackage{lineno,hyperref}
\modulolinenumbers[5]

\usepackage[figuresright]{rotating}

\def\e{{\epsilon}}

\def\norm#1{\|#1\|} 

\newcommand{\beq} {\begin{equation}}
\newcommand{\eeq} {\end{equation}}
\newcommand{\bdm} {\begin{displaymath}}
\newcommand{\edm} {\end{displaymath}}
\newcommand{\bit}{\begin{itemize}}
\newcommand{\eit}{\end{itemize}}
\newcommand{\bde}{\begin{description}}
\newcommand{\ede}{\end{description}}
\newcommand{\bce}{\begin{center}}
\newcommand{\ece}{\end{center}}
\newcommand{\ben} {\begin{enumerate}}
\newcommand{\een} {\end{enumerate}}
\newcommand{\bea} {\begin{eqnarray}}
\newcommand{\eea} {\end{eqnarray}}
\newcommand{\barr} {\begin{array}}
\newcommand{\earr} {\end{array}}
\newcommand{\bean} {\begin{eqnarray*}}
\newcommand{\eean} {\end{eqnarray*}}
\newcommand{\edoc} {\end{document}}

\newcommand{\dd}[2] {\ensuremath{\frac{\partial {#1}}{\partial {#2}}}}
\newcommand{\DD}[2] {\ensuremath{\frac{d {#1}}{d {#2}}}}
\newcommand{\Grad} {\ensuremath{\nabla}}

\newcommand{\Div} {\ensuremath{\nabla\cdot}}

\newcommand{\eval}[2][\right]{\relax
  \ifx#1\right\relax \left.\fi#2#1\rvert}

\providecommand{\norm}[1]{\left\hspace{-0.5ex}\lVert#1\right\hspace{-0.5ex}\rVert}

\newcommand{\Nel} {\ensuremath{{N_\text{el}}}}

\newcommand{\half} {\ensuremath{\frac{1}{2}}}

\newcommand{\mc}[1]{\mathcal{#1}}

\newcommand{\mb}[1]{\mathbf{#1}}

\newcommand{\nor}[1]{\left\vert\kern-0.35ex\left\vert #1 \right\vert\kern-0.35ex\right\vert}
\newcommand{\norL}[1]{\left\vert\kern-0.5ex\left\vert #1 \right\vert\kern-0.5ex\right\vert}
\newcommand{\nort}[1]{{\left\vert\kern-0.35ex\left\vert\kern-0.35ex\left\vert #1 
    \right\vert\kern-0.35ex\right\vert\kern-0.35ex\right\vert}}
\newcommand{\snor}[1]{\left| #1 \right|}
\newcommand{\LRp}[1]{\left( #1 \right)}
\newcommand{\LRs}[1]{\left[ #1 \right]}
\newcommand{\LRa}[1]{\left< #1 \right>}
\newcommand{\LRc}[1]{\left\{ #1 \right\}}
\newcommand{\pp}[2]{\frac{\partial #1}{\partial #2}}

\newcommand{\td}[2]{\frac{d #1}{d #2}}

\newcommand{\bs}[1]{\boldsymbol{#1}}

\newcommand{\jump}[1] {\ensuremath{\LRs{\!\!\LRs{{#1}}\!\!}}}
\newcommand{\jumpL}[1] {\ensuremath{\LRs{\hspace{-1ex}\LRs{#1}\hspace{-1ex}}}}

\newcommand{\average}[1] {\ensuremath{\LRc{\hspace{-0.9ex}\LRc{#1}\hspace{-0.9ex}}}}
\newcommand{\averageM}[1] {\ensuremath{\LRc{\hspace{-1.ex}\LRc{#1}\hspace{-1.ex}}}}
\newcommand{\averageL}[1] {\ensuremath{\LRc{\hspace{-2ex}\LRc{#1}\hspace{-2ex}}}}

\newcounter{tablecell}

\newtheorem{coro}{Corollary}
\newtheorem{theo}{Theorem}
\newtheorem{lem}{Lemma}

\newtheorem{rema}{Remark}

\newcommand{\figlab}[1]{\label{fig:#1}}
\newcommand{\eqnlab}[1]{\label{eq:#1}}
\newcommand{\theolab}[1]{\label{theo:#1}}
\newcommand{\corolab}[1]{\label{coro:#1}}

\newcommand{\lemlab}[1]{\label{lem:#1}}

\newcommand{\tablab}[1]{\label{tab:#1}}

\newcommand{\figref}[1]{\ref{fig:#1}}
\newcommand{\theoref}[1]{\ref{theo:#1}}

\newcommand{\cororef}[1]{\ref{coro:#1}}

\newcommand{\lemref}[1]{\ref{lem:#1}}
\newcommand{\eqnref}[1]{\eqref{eq:#1}}

\newcommand{\seclab}[1]{\label{sect:#1}}
\newcommand{\secref}[1]{\ref{sect:#1}}
\newcommand{\tabref}[1]{\ref{tab:#1}}

\renewcommand{\v}{v}
\renewcommand{\u}{u}
\newcommand{\un}{\u_h}
\renewcommand{\e}{e}
\newcommand{\w}{w}
\newcommand{\z}{z}

\newcommand{\q}{q}

\newcommand{\qh}{\hat{\q}}

\newcommand{\qb}{{\mb{\q}}}

\newcommand{\ut}{\tilde{\u}}

\newcommand{\vb}{\mb{\v}}
\newcommand{\ub}{\mb{\u}}

\newcommand{\PL}{\mathbb{L}}
\newcommand{\PN}{\mathbb{N}}

\newcommand{\R}{{\mathbb{R}}}

\newcommand{\f}{f}

\newcommand{\fb}{\mb{\f}}

\newcommand{\pOmega}{\partial \Omega}

\newcommand{\K}{K}
\newcommand{\Kp}{K^+}
\newcommand{\Km}{K^-}
\newcommand{\pK}{{\partial \K}}
\newcommand{\Kj}{{\K_j}}
\newcommand{\Gh}{{\mc{E}_h}}
\newcommand{\Gho}{{\mc{E}_h^o}}
\newcommand{\Ghb}{{\mc{E}_h^{\partial}}}
\newcommand{\Poly}{{\mc{P}}}

\newcommand{\p}{{p}}

\newcommand{\Vb}{{\mb{V}}}
\newcommand{\V}{{V}}
\newcommand{\Lam}{{\Lambda}}
\newcommand{\Lamb}{\bs{\Lambda}}

\newcommand{\lambdab}{\bs{\lambda}}

\newcommand{\Vh}{{\V_h}}
\newcommand{\Vbh}{{\Vb_h}}
\newcommand{\Lamh}{{\Lam_h}}
\newcommand{\Lambh}{{\Lamb_h}}

\newcommand{\VhK}{{\V_h\LRp{\K}}}

\newcommand{\VbhK}{{\Vb_h\LRp{\K}}}
\newcommand{\Lamhe}{{\Lam_h\LRp{\e}}}
\newcommand{\Lambhe}{{\Lamb_h\LRp{\e}}}

\newcommand{\Lte}{{L^2\LRp{\Gh}}}
\newcommand{\n}{\mb{n}}
\newcommand{\nm}{\n^-}
\newcommand{\np}{\n^+}
\newcommand{\uh}{{\hat{\u}}}

\renewcommand{\c}{{c}}
\newcommand{\Oh}{{\Omega_h}}

\newcommand{\cIT}{\c_{IT}}
\newcommand{\cq}{\c^q}
\newcommand{\cu}{\c^u}

\renewcommand{\sb}{\mb{s}}

\newcommand{\veps}{{\varepsilon}}

\newcommand{\vepsh}{\veps^h}
\newcommand{\xin}{\xi^n}
\newcommand{\rhon}{\rho^n}
\newcommand{\vepsI}{\veps^I}

\newcommand{\E}{E}

\newcommand{\Scal}{\mathcal{S}}
\newcommand{\pOmegah}{{\pOmega_h}}
\newcommand{\Omegah}{{\Omega_h}}

\newcommand{\vtest}{{\bf v}}
\newcommand{\vhtest}{{\hat{\bf v}}}

\newcommand{\xb}{{\bf x}}

\newcommand{\unm}{{u_\nu}}
\newcommand{\nb}{{\bf n}}

\newcommand{\Fcal}{\mathcal{F}}

\newcommand{\Fcals}{\mathcal{F}^*}

\newcommand{\Res}{{\mathcal Res}}

\newcommand{\dt}{{\triangle t}}

\newcommand{\pres}{{{p}}}

\newcommand{\Ical}{\mathcal{I}}

\newcommand{\gammam}{{\tilde{\gamma}}}

\usepackage{color}
\usepackage{soul,xargs}
\usepackage[pdftex,dvipsnames]{xcolor}

\usepackage[colorinlistoftodos,prependcaption,textsize=tiny]{todonotes}

\newcommandx{\question}[2][1=]{\todo[linecolor=red,backgroundcolor=red!25,bordercolor=red,#1]{#2}}
\newcommandx{\change}[2][1=]{\todo[linecolor=blue,backgroundcolor=blue!25,bordercolor=blue,#1]{#2}}
\newcommandx{\add}[2][1=]{\todo[linecolor=OliveGreen,backgroundcolor=OliveGreen!25,bordercolor=OliveGreen,#1]{#2}}
\newcommandx{\improve}[2][1=]{\todo[linecolor=Plum,backgroundcolor=Plum!25,bordercolor=Plum,#1]{#2}}
\newcommandx{\thiswillnotshow}[2][1=]{\todo[disable,#1]{#2}}
\newcommandx{\remove}[2][1=]{\todo[linecolor=yelllow,backgroundcolor=yellow!10,bordercolor=red,#1]{#2}}

\begin{document}

\begin{frontmatter}



\dochead{}

\title{A scalable exponential-DG approach for nonlinear conservation laws: with application to Burger and Euler equations}


\author[AddrShinhoo]{Shinhoo Kang\corref{mycorrespondingauthor}}
\cortext[mycorrespondingauthor]{Corresponding author  }
\ead{shinhoo.kang@anl.gov}
\author[AddrTan]{Tan Bui-Thanh}

\address[AddrShinhoo]{Mathematics and Computer Science Division Argonne National Laboratory Lemont, IL 60439, USA.}

\address[AddrTan]{Department of Aerospace Engineering and Engineering Mechanics, and Oden Institute for Computational Engineering
  \& Sciences, The University of Texas at Austin, Austin, TX 78712, USA.}
\ead{tanbui@oden.utexas.edu}

\begin{abstract} We propose an Exponential DG approach for numerically solving partial differential equations (PDEs). 
The idea is to decompose the governing PDE operators 
into linear (fast dynamics extracted by linearization) and nonlinear (the remaining after removing the former) parts, on which we apply the discontinuous Galerkin (DG) spatial discretization.
The resulting semi-discrete system is then integrated using exponential time-integrators: exact for the former and approximate for the latter. By construction, our approach 
  i) is stable with a large Courant number ($Cr>1$); 
  ii) supports high-order solutions both in time and space;
  iii) is computationally favorable compared to IMEX DG methods 
  with no preconditioner;
  iv) 
  requires comparable computational time compared to explicit RKDG methods, while having time stepsizes orders magnitude larger than maximal stable time stepsizes for  explicit RKDG methods;
  v) is scalable in a modern massively parallel computing architecture 
  by exploiting Krylov-subspace matrix-free  exponential time integrators and compact communication stencil of DG methods. 
  Various numerical results for both Burgers and Euler equations are presented to showcase these expected properties.
For Burgers equation, we present a detailed stability and convergence analyses for the exponential Euler DG scheme.
  
\end{abstract}



\end{frontmatter}



\section{Introduction}

The discontinuous Galerkin (DG) method has gain popularity for decades as a spatial discretization.
The DG method\textemdash originally developed \cite{ReedHill73,LeSaintRaviart74,johnson1986analysis} for the neutron transport equation\textemdash has
been studied extensively for various types of partial
differential equations (PDEs)
including Poisson type equation \cite{wheeler1978elliptic,arnold1982interior,cockburn2000development,arnold2002unified},
poroelasticity \cite{liu2004discontinuous}, 
shallow water equations \cite{nair2005discontinuous,GiraldoWarburton08,GiraldoRestelli10,wintermeyer2017entropy},
Euler and Navier-Stokes equations \cite{bassi2005discontinuous,gassner2016split},
Maxwell equations \cite{fezoui2005convergence,bui2012analysis}, 
solid dynamics \cite{noels2008explicit},
magma dynamics \cite{tirupathi2015modeling}, to name a few.
One of the reason is that DG methods are well-suited for parallel-computing
 due to the local nature of the methods. 
DG methods combine advantages of finite volume and finite element methods
in the sense that a global solution is approximated by a finite set of local functions,
and each local element communicates with its adjacent element through numerical flux on element boundary.  
Since the numerical flux is calculated using the state variables on the face,
DG methods have compact stencil, hence reduces inter-communication cost. 
Another reason can be the positive properties of the scheme, i.e.,  
flexibilty for handling complex geometry, hp-adaptivity, 
high-order accuracy, upwid stabilization, etc
\cite{wheeler1978elliptic,arnold1982interior,cockburn2000development,arnold2002unified,demkowicz2006computing}.


To fully discretize a time-dependent partial differential equation (PDE),
temporal discretization is also necessary. Explicit time integrators
such as Runge-Kutta methods are popular due to their simplicity and 
ease in computer implementation. However, scale-separated or geometrically-induced stiffness limits the time-step size
severely for high-order DG methods (see, e.g., \cite{kanevsky2007application,GiraldoWarburton08}). 
For long-time integration
this can lead to an excessive number of time steps,
and hence substantially taxing computing and storage resources. On the
other hand, fully-implicit methods could be expensive, especially for
nonlinear PDEs for which Newton-like methods are typically
required. Semi-implicit time-integrators have been designed to relax
the time-step size restriction caused by the stiffness in order to reduce the computational burden
arising from the linear solve
\cite{ascher1997implicit,Kennedy2003additive,pareschi2005implicit}.
 In the context of low-speed fluid flows, including Euler, Navier-Stokes, and
shallow water equations, implicit-explicit (IMEX) DG methods have been proposed and demonstrated 
to be more advantageous than either explicit or fully-implicit DG
methods \cite{FeistauerDolejsiKucera07,RestelliGiraldo09,kang2019imex}. 
 The common feature of these methods is that they
relax the stiffness condition by employing implicit time-stepping schemes for handling the linear stiff part of the
PDE. Therefore the performance highly depends on a linear solver, 
which means 
 an appropriate preconditioner needs to be constructed for achieving decent performance. 
However, developing such a preconditioner is not a trivial task and it is problem-specific. 
 
Alternatively, exponential time integrators have been received great attention due to the positive characteristics such as stability and accuracy.
  The methods have been applied 
  to various types of PDEs
  including linear advection-diffusion equations
  \cite{caliari2004interpolating}, 
  Schr{\"o}dinger equation \cite{celledoni2008symmetric},
  Maxwell equations \cite{botchev2006gautschi},
  magnetohydrodynamics (MHD) equations \cite{tokman2002three},
  Euler equations \cite{li2018exponential},
  incompressible Navier-Stokes equations \cite{kooij2018exponential}, 
  compressible Navier-Stokes equations \cite{schulze2009exponential,li2019exponential}, 
  shallow water equations \cite{clancy2013use},
 among others. 

Exponential time integrators 
 is similar to IMEX methods in the sense of splitting a governing equation into stiff and non-stiff parts. 
However, 
  exponential time integrators exactly integrate the linear stiff part by multiplying an integrating factor instead of using a quadrature in time. 
Compared to IMEX methods, exponential time integrators 
replaces a linear solve at each time step with a  computationally demanding matrix exponential.

Many researchers have conducted various studies to mitigate the challenge,
 one way is to use Krylov subspace, 
  where a large matrix is projected onto a small Krylov subspace
  so that computing the matrix exponential becomes less expensive. 
  To improve the Krylov subspace projection-based algorithm,
  rational Krylov method \cite{moret2007rd,ragni2014rational}, 
  restart Krylov method \cite{tal2007restart,afanasjew2008implementation},
  block Krylov method \cite{botchev2013block,kooij2017block},
  adaptive Krylov method \cite{niesen2012algorithm} have been developed.
  Lately, \cite{gaudreault2018kiops} and \cite{luan2019further} 
  enhance the computational efficiency of the adaptive Krylov method 
  by replacing the Arnoldi procedure \cite{arnoldi1951principle} 
  with the incomplete orthogonalization procedure 
  \cite{koskela2015approximating,vo2017approximating}. 
  The work in \cite{tokman2006efficient} shows that
  the exponential propagation iterative (EPI) schemes can outperform the standard implicit Newton-Krylov integrators 
  with no preconditioning.
  The work in \cite{clancy2013use} observes the second-order EPI2 provides comparable results to the explicit fourth-order Runge-Kutta (RK4). 
  For elastodynamic problems, 
  the second-order Gautschi-type exponential integrator
  outperform the backward Euler integrators \cite{michels2014exponential,michels2017stiffly}. 
 
In this study, we propose an Exponential DG framework for partial differential equations. 
 To that end, we separate governing equations into linear and nonlinear parts, 
  to which we apply the DG spatial discretization. 
  The former is integrated analytically, whereas the latter is approximated. 
  Since the method does not require any linear solve, 
  it has a potential 
to be scalable in a modern massively parallel computing architecture.
  The proposed Exponential DG method: 
  i) is stable with a large Courant number ($Cr>1$); 
  ii) exploits high-order solutions both in time and space;
  iii) is more efficient than IMEX DG methods with no preconditioner;
  iv) is comparable to explicit RKDG methods on uniform mesh and more beneficial on non-uniform grid for Euler equations;
  v) provides promising weak and strong scalable parallel solutions. 

In the following, we discuss Exponential DG framework  
in Section \secref{ExponentialDG}.
In Section \secref{ModelProblems}, 
we apply Exponential DG framework to Burgers equation and Euler equations,
where we show the construction of linear operator based on a flux Jacobian. 
Then, we presents a detailed analysis on the stability and convergence of the exponential DG scheme for Burgers equation. 
The performance of the proposed method will be discussed in Section \secref{NumericalResults}
with several numerical examples for both Burgers and Euler equations.
We finally conclude the paper in Section \secref{Conclusion}.

\section{Exponential DG framework}
\seclab{ExponentialDG}

In this section, we present the key idea behind Exponential DG framework. We first split a given PDE into a linear and a nonlinear parts, to which DG discretization and exponential time integrators are applied. We begin with notations and conventions used in the paper. 

\subsection{Finite element definitions and notations}
Let $\Omega$ 
be an open and bounded subset of $\R^d$, where $d = \LRc{1,2,3}$ is
the spatial dimention.
We denote by
$\Omega_h := \cup_{i=1}^\Nel \K_i$ the mesh containing a finite
collection of non-overlapping elements, $\K_i$, that partition
$\Omega$.  Here, $h$ is defined as $h := \max_{j\in
  \LRc{1,\hdots,\Nel}}diam\LRp{\Kj}$. Let $\pOmega_h := \LRc{\pK:\K
  \in \Omega_h}$ be the collection of the faces of all elements. 
Let us define $\Gh$
as the skeleton
of the mesh which consists of the set of all uniquely defined faces.
For
two neighboring elements $\Kp$ and $\Km$ that share an interior
interface $\e = \Kp \cap \Km$, we denote by $q^\pm$ the trace of their
solutions on $\e$. 
We define $\nm$ as the unit outward normal vector on
the boundary $\pK^-$ of element $\Km$, and $\np = -\nm$ the unit outward
normal of a neighboring element $\Kp$.  On the interior interfaces $\e
\in \Gho$, we define the mean/average operator $\average{\bf v}$, where $\bf v$ is
either a scalar or a vector quantify, as
$\average{{\bf v}}:=\LRp{{\bf v}^- + {\bf v}^+}/2$, 
and the jump operator $\jump{\bf v} := 2\average{\bf v \cdot \nb}$.
 On the boundary faces $\e \in \Ghb$, unless otherwise stated, we define the mean and jump operators as 
$\average{{\bf v}}:={\bf v}, \quad \jump{{\bf v}} :={\bf v}\cdot\nb$.

Let $\Poly^{k}\LRp{D}$ denote the space of polynomials of degree at
most $k$ on a domain $D$. Next, we introduce discontinuous
piecewise polynomial spaces for scalars and vectors as
\begin{align*}
\Vh\LRp{\Omega_h} &:= \LRc{v \in L^2\LRp{\Omega_h}:
  \eval{v}_{\K} \in \Poly^k\LRp{\K}, \forall \K \in \Omega_h}, \\
\Lamh\LRp{\Gh} &:= \LRc{\lambda \in \Lte:
  \eval{\lambda}_{\e} \in \Poly^k\LRp{\e}, \forall \e \in \Gh},\\
\Vbh\LRp{\Omega_h} &:= \LRc{{\bf v} \in \LRs{L^2\LRp{\Omega_h}}^m:
  \eval{{\bf v}}_{\K} \in \LRs{\Poly^k\LRp{\K}}^m, \forall \K \in \Omega_h},\\
\Lambh\LRp{\Gh} &:= \LRc{\lambdab \in \LRs{\Lte}^m:
  \eval{\lambdab}_{\e} \in \LRs{\Poly^k\LRp{\e}}^m, \forall \e \in \Gh}.
\end{align*}
and similar spaces $\VhK$, $\Lamhe$, $\VbhK$, and $\Lambhe$ by
replacing $\Omega_h$ with $\K$ and $\Gh$ with $\e$. Here, $m$ is the
number of components of the vector under consideration.

We define $\LRp{\cdot,\cdot}_\K$ as the $L^2$-inner product on an
element $\K \in \R^d$, and $\LRa{\cdot,\cdot}_{\pK}$ as the
$L^2$-inner product on the element boundary $\pK \in
\R^{d-1}$.  We also define the broken inner products as
$\LRp{\cdot,\cdot}_\Omega := \LRp{\cdot,\cdot}_{\Omega_h} :=
\sum_{\K\in \Omega_h}\LRp{\cdot,\cdot}_\K$ and
$\LRa{\cdot,\cdot}_{\pOmega} := \LRa{\cdot,\cdot}_{\pOmega_h} :=
\sum_{\pK\in \pOmega_h}\LRa{\cdot,\cdot}_\pK$, and on the mesh
skeleton as $\LRa{\cdot,\cdot}_\Gh := \sum_{\e\in
  \Gh}\LRa{\cdot,\cdot}_\e$.
  We also define associated norms as 
$\norm{\cdot}_{\Omegah}:= \LRp{ \sum_{K\in \Omegah} \norm{\cdot}_{K}^2 }^\half$ and
$\norm{\cdot}_{\pOmegah}:= \LRp{ \sum_{K\in \Omegah} \norm{\cdot}_{\pK}^2 }^\half$ where   
$\norm{\cdot}_{\K}=\LRp{\cdot,\cdot}_{\K}^\half$ and 
$\norm{\cdot}_{\pK} = \LRa{\cdot,\cdot}_{\pK}^\half$. 


\subsection{Constructing linear and nonlinear DG operators for conservation laws}

We consider conservation laws governed by a generic system of partial differential equation (PDE): 
\begin{align}
  \eqnlab{pde}
   \dd{\qb}{t} + \Div \Fcal = \sb, \quad \text{ in } \Omega,
\end{align}
where $\qb$ is the conservative variable, $\Fcal=\Fcal(\qb)$ is 
the flux tensor, and $\sb$ is the source vector. 
We seek a stiff linear flux $\Fcal_L$ 
that, we assume, captures the rapidly changing dynamics in the system.
Inspired by the works in \cite{tokman2006efficient, giraldo2010high}, we use a flux Jacobian to define the linear flux, i.e., 
\begin{align}
  \eqnlab{linearflux}
  \Fcal_L := \dd{\Fcal}{\qb} \Bigg\vert_{\tilde{\qb}} \qb,
\end{align}
where $\tilde{\qb}$ is a reference state.
By adding and subtracting the linear flux $\Fcal_L$ in \eqnref{pde},
 we split the divergence term into a linear (stiff) part $\Div{\Fcal_L}$ and a nonlinear (non-stiff)
part $\Div{\LRp{\Fcal - \Fcal_L}}$.
Similarly, we decompose the source term $\sb$ into a linear term $\sb_L$,
\begin{align}
  \eqnlab{linearsource}
  \sb_L := \dd{\sb}{\qb} \Bigg\vert_{\tilde{\qb}} \qb,
\end{align}
 and a nonlinear term $\sb - \sb_L$.
Thus \eqnref{pde} becomes 
\begin{align}
  \eqnlab{split-goveq-continuous}
   \dd{\qb}{t} + \Div \Fcal_L + \Div \Fcal_{\mc{N}} = \sb_L + \sb_{\mc{N}}
\end{align}
at the continuous level, 
where $\Fcal_{\mc{N}} = \Fcal - \Fcal_L$ and $\sb_{\mc{N}} = \sb - \sb_L$. 
The decomposition at continuous level avoids complicated derivatives
of the stabilization parameter coming from a numerical flux in DG
methods when applying exponential time integrators.
To the rest of the paper, except for the analysis in Section \secref{AppendixAnalysis}, we use the same notations for the exact and the DG solutions for simplicity. A semi-discrete form of \eqnref{split-goveq-continuous}  using DG
discretization for spatial derivatives
 reads: 
find $\qb \in \Vbh(\Omegah)$ such that
    \begin{align}
    \eqnlab{EXPO-SemiDiscretizedForm_ode}
      \LRp{\pp{\qb}{t},\vtest}_\Omegah 
       &= \LRa{L\qb,\vtest} + \LRa{\mc{N}(\qb),\vtest},
    \end{align}
 for all ${\bf v}\in \Vbh(\Omegah)$, 
  where 
  \begin{align*}
    \LRa{\mc{N}(\qb),\vtest} &:= - \LRp{\Div \Fcal_{\mc{N}}\LRp{\qb}, \vtest }_\Omegah
                  + \LRp{\sb_{\mc{N}}(\qb),\vtest}_\Omegah
                  - \LRa{ \LRp{ \Fcals_{\mc{N}}\LRp{\qb^\pm}
                              - \Fcal_{\mc{N}}\LRp{\qb} } \cdot \n,\vtest}_\pOmegah,\\
    \LRa{L\qb,\vtest}  &:= - \LRp{\Div \Fcal_L\LRp{\qb}, \vtest }_\Omegah
                  + \LRp{\sb_L(\qb),\vtest}_\Omegah
                  - \LRa{ \LRp{ \Fcals_L\LRp{\qb^\pm}
                              - \Fcal_L\LRp{\qb} } \cdot \n,\vtest}_\pOmegah.
  \end{align*} 
Here, 
$\Fcal_L^*$ and $\Fcal_{\mc{N}}^*:=\Fcal^* - \Fcal_L^*$ 
are a linear and a nonlinear DG numerical flux, respectively, 
such that
\begin{align*}
  \LRa{\jump{ \Fcals_{\mc{N}}\LRp{\qb^\pm} \cdot \n},\vhtest}_\Gh = 0
\text{ and } 
  \LRa{\jump{ \Fcals_{L}\LRp{\qb^\pm} \cdot \n},\vhtest}_\Gh = 0,
\end{align*}
for all $\vhtest \in \Lamh\LRp{\Gh}$. At this point both spatial linear and the nonlinear operators are discretized with DG, and  
we discuss exponential time integrators for temporal derivative next.

\subsection{Exponential time integrators}

For the clarity of the exposition, let us rewrite \eqnref{EXPO-SemiDiscretizedForm_ode} as
\begin{align}
    \eqnlab{splitDG}
    \DD{\qb}{t} = L\qb + \mc{N}(\qb), \quad t\in\LRp{0,T},
  \end{align}  
with an initial condition $\qb_0 = \qb(0)$. An abuse of notations has been made for brevity: 
first, the proper form for both $L$ and $\mc{N}$  would be a composition with a projection operator onto $\Vh\LRp{\Oh}$ as in Section \secref{AppendixAnalysis}; second, we don't distinguish $\qb$ with its nodal (or modal) vector; and third,  $L$ and $\mc{N}$ are used interchangeably with their matrix representations from $\Vh\LRp{\Oh}$ to $\Vh\LRp{\Oh}$ (see also Section \secref{AppendixAnalysis}).
Now multiplying \eqnref{splitDG} with integrating factor $e^{-\dt L}$ yields
\begin{align}
\eqnlab{nonlinear_ivp_sol}
  \qb(t^{n+1}) = e^{\dt L} \qb(t^n) 
    + \int_0^\dt e^{(\dt - \sigma) L} \mc{N}(\qb(t^n+\sigma) ) d\sigma
\end{align}
via a simple application of the method of variation of constants. 
At this point, \eqnref{nonlinear_ivp_sol} is exact.
The first term $e^{\dt L}\qb(t^n)$ is the homogeneous solution, 
whereas the second term is the particular solution that involves a convolution integral with the matrix exponential.
Various exponential integrators have been proposed to approximate 
\eqnref{nonlinear_ivp_sol} in different ways. 
In particular, 
 a $p$th-order time polynomial approximation to the nonlinear map $\mc{N}$ can be written as
 \begin{align*}
  \mc{N}(\qb(t^n+\sigma))
   = \sum_{j=0}^{p-1} \frac{\LRp{t^n+\sigma}^j}{j!}\vb_{j+1} + \mc{O}\LRp{\dt^p}, 
\end{align*}
with appropriate choice \cite{skaflestad2009scaling} for  $\vb_j$.
 This allows us to express an approximation of $\qb(t^{n+1})$ in \eqnref{nonlinear_ivp_sol}, 
 denoted as $\qb^{n+1}$,
as a linear combination of $\varphi$-functions \cite{niesen2012algorithm}, i.e., 
\begin{align}
  \eqnlab{linear_combination_of_phifunctions}
  \qb^{n+1} = \sum_{i=0}^p (\dt)^i \varphi_i (\dt L) {\bf b}_i,
\end{align}
where we have defined ${\bf b}_0:=\qb(t^n)$ and 
${\bf b}_i:=\sum_{j=0}^{p-i}\frac{(t^n)^j}{j!} \vb_{i+j}$.
Here, $\varphi_i$-functions for a scalar $\tau$ are defined by
\begin{align*}
  \varphi_i(\tau) := \int_0^1 e^{(1-z)\tau} \frac{z^{i-1} }{(i-1)!} dz
\end{align*}
with 
$\varphi_0(\tau):= e^{\tau}$. It is easy to see the recurrence relation
 $\LRs{\varphi_i(\tau) - \varphi_i(0)} \tau^{-1} = \varphi_{i+1} (\tau)$ and $\varphi_i(0) = \frac{1}{i!}$ hold true. The definition of $\varphi_i$-functions for matrices is straightforward, e.g., based on Jordan canonical form \cite{horn1994topics}.



For efficient computation of \eqnref{linear_combination_of_phifunctions},
Krylov subspace methods \cite{friesner1989method,
gallopoulos1992efficient,
hochbruck1998exponential} with the exponential of the augmented matrix can be used \cite{sidje1998expokit,al2009new,niesen2012algorithm,gaudreault2018kiops}, 
in which an augmented matrix is constructed and projected onto a small Krylov subspace so that matrix exponential is amenable to compute. 

In this paper, we use \texttt{KIOPS} \cite{gaudreault2018kiops} algorithm 
\footnote{
The performance of the adaptive Krylov subspace solver depends on several parameters such as the size of Krylov space. We empirically determine the parameters in this studies.
} 
for serial computations.
For parallel computation, we have implemented an exponential time integrator based on the \texttt{KIOPS} algorithm for our C\texttt{++} DG finite element library
(a spin-off from \texttt{mangll}
\cite{wilcox2010high}).
Thanks to DG discretization and explicit nature of the exponential integrators, 
the proposed method is highly 
parallel as the communication cost can be  effectively overlapped by computation. 

\section{Model problems}
\seclab{ModelProblems}

The key in  the operator splitting in \eqnref{EXPO-SemiDiscretizedForm_ode}
is the linearized flux $\Fcal_L$ \eqnref{linearflux}. 
In this section we
choose
 Burgers  and Euler equations as prototypes for the generic conservation law \eqnref{pde} and construct $\Fcal_L$ (hence $L$ and $\mc{N}$) for these equations.

\subsection{Burgers equation}
Burgers equation is a quasi-linear parabolic PDE
that comprises of nonlinear convection and linear diffusion:
\begin{align}
\eqnlab{pde-burgers-eq}
 \dd{u}{t} + \half\dd{u^2}{x} 
 = \dd{}{x} \LRp{\kappa \dd{u}{x} } \text{ in } \Omega,
\end{align}
where $u$ is a scalar quantity and $\kappa >0 $ is the constant viscosity. 
The linearization $\Fcal_L$
\eqnref{linearflux}  of $\Fcal :=
u^2/2 - \kappa \dd{u}{x}$ evaluated at $\tilde{u}$  is given by
\begin{align*}
  \Fcal_L:= \tilde{u} u - \kappa \dd{u}{x},
\end{align*}
and thus
\[
  \Fcal_{\mc{N}} := \frac{u^2}{2} - \tilde{u}u.
\]
We can now write \eqnref{pde-burgers-eq} in the form \eqnref{split-goveq-continuous} as
\begin{align}
\eqnlab{pde-burgers-eq-split}
 \dd{u}{t} 
 + 
  \dd{}{x}\underbrace{ \LRp{ \tilde{u} u  - \kappa \dd{u}{x}}
  }_{\Fcal_L}
   + 
   \dd{}{x}\underbrace{ \LRp{ \frac{u^2}{2} - \tilde{u}u } 
   }_{\Fcal_{\mc{N}}} = 0
  \text{ in } \Omega
\end{align}
where $\tilde{u}$ is a reference state (for example $\tilde{u} = u^n$: the numerical solution at $t^n$). 
The DG weak formulation of \eqnref{pde-burgers-eq-split} reads: seek $q,u \in \Vh\LRp{\Omegah}$ such that 
\begin{subequations}
\eqnlab{gov-burgers1d}
\begin{align}
 \eqnlab{gov-burgers1d-gradu}
 \LRp{q,p}_\Omegah &= 
   \LRp{\dd{u}{x},p}_\Omegah + \LRa{\nb(u^{**} - u),p}_\pOmegah, \\
   \eqnlab{gov-burgers1d-u}
 \LRp{\dd{u}{t},v}_\Omegah 
   &:= \LRa{Lu,v} + \LRa{\mc{N}(u),v},
\end{align}
\end{subequations}
where
\begin{subequations}
\eqnlab{Burger}
\begin{align}
\eqnlab{BurgerL}
 \LRa{Lu,v} &:= -\LRp{\kappa q - \tilde{u} u  ,\pp{v}{x}}_\Omegah
       + \LRa{ \nb 
                      \LRp{ \kappa q^{**} - \LRp{\tilde{u}u}^* }  ,v }_\pOmegah, \\
                    \eqnlab{BurgerN}
 \LRa{\mc{N}(u),v} &:= -\LRp{\tilde{u} u - \half u^2 ,\pp{v}{x}}_\Omegah
       + \LRa{ \nb
                      \LRp{ \LRp{\tilde{u}u}^* -  \LRp{\frac{u^2}{2}}^* } ,v }_\pOmegah,
\end{align}
\end{subequations}
for all $p,v \in \Vh\LRp{\Omega_h}$ 
with $\nb=\pm 1$ for one-dimensional problems.
Here, we use the central flux for $u^{**}$ and $q^{**}$ for diffusion operator, i.e.,
$q^{**} = \average{\q}$ and $u^{**} = \average{\u}$, 
the entropy flux for inviscid Burgers part, i.e., 
$\LRp{\frac{u^2}{2}}^* := \frac{1}{3}\LRp{\averageM{\half u^2} + \average{u}^2} + \frac{\sigma}{h}\jump{\u}$,
 where $\sigma \ge 0$ is a constant (the Lax-Friedrichs flux  
$\LRp{\frac{u^2}{2}}^* := \average{\half u^2} + \half \max(|u^\pm|)\jump{u}$ is also considered
to compare with the entropy flux), 
and the Lax-Friedrichs flux for linear Jacobian part, i.e., 
%
$ \LRp{\tilde{u}u}^* := \average{\tilde{u} u} + \half \max(|\tilde{u}^\pm|)\jump{
  u}$. 
With the central flux $u^{**}$, $q$ can be computed locally element-by-element
from \eqnref{gov-burgers1d-gradu}, and the only actual (global) unknown is $u$ in \eqnref{gov-burgers1d}.

\subsection{Euler equations}

We consider the compressible Euler equations written in the following form
\begin{subequations}
\eqnlab{euler-gov}
\begin{align}
  \dd{\rho}{t}    + \Div\LRp{\rho\ub} &= 0,\\
  \dd{\rho\ub}{t} + \Div\LRp{\rho\ub\otimes\ub + \pres\Ical} &= 0, \\ 
  \dd{\rho E}{t} + \Div\LRp{\rho \ub H} &= 0, 
\end{align}
\end{subequations}
where $\rho$ is the density,
$\ub$  the velocity,
$\pres$  the pressure,
$\rho E = \rho e + \half \rho \norm{\ub}^2$  the total energy,
$e=\frac{\pres}{\rho(\gamma -1)}$  the internal energy,
$\pres$ the pressure,
$H = E + \frac{\pres}{\rho} = \frac{a^2}{\gamma -1} + \frac{1}{2}\norm{\ub}^2$ the total specific enthalpy,
$a = \sqrt{\gamma \pres/\rho}$ 
 the sound speed, 
$\gamma$ the ratio of the specific heat, and
$\Ical$ the $d\times d$ 
identity matrix.
In a compact form, \eqnref{euler-gov} can be written as
\begin{align}
\eqnlab{euler3dd-gov}
    \dd{\qb}{t} + \Div{\Fcal(\qb)} = 0,
\end{align}
with $\qb = (\rho, \rho \ub, \rho E)^T$, $\sb=(s_\rho, s_{\rho \ub}, s_{\rho E})^T$, and 
$\Fcal(\qb) = (\rho u, \rho \ub \otimes \ub + \pres \Ical, \rho \ub H)^T$.
Let us define 
$\unm:=\nb \cdot \ub$,
$\phi:=\LRp{\frac{\gamma-1}{2}} \norm{\ub}^2$, $\gammam:=\gamma-1$,
the flux Jacobian $A:=\dd{\Fcal}{\qb}$ where
\begin{align*}
  A= 
  \begin{pmatrix}
    0 & \nb^T & 0\\
    \phi \nb - \ub \unm & \ub \otimes \nb - \gammam \nb \otimes \ub + \unm \Ical & \gammam \nb \\
    \LRp{\phi - H}\unm  & H \nb^T -\gammam \ub^T \unm & \gamma \unm
  \end{pmatrix}.
\end{align*}
The linearized flux \eqnref{linearflux} in this case is defined
$\Fcal_L:= A(\tilde{\qb}) \qb =: \tilde{A} \qb$,
 and \eqnref{split-goveq-continuous} now reads
\begin{align}
  \eqnlab{euler-split-3dd}
  \dd{\qb}{t} +
  \Div \underbrace{ \LRp{\tilde{A}\qb} }_{\Fcal_L} 
  + \Div \underbrace{ \LRp{  \Fcal(\qb)  -  \tilde{A}\qb }
  }_{\Fcal_{\mc{N}}} = {\bf 0}.
\end{align}
By multiplying a test function $\vtest$ to \eqnref{euler-split-3dd},
integrating by parts 
for each element, 
and summing all the elements we arrive at the semi-discretization with DG:  
seek $\qb \in \Vbh(\Omegah)$ such that
\begin{align*}
  \LRp{\dd{\qb}{t},\vtest}_\Omegah
  = \LRa{L\qb,\vtest} + \LRa{\mc{N}(\qb),\vtest}, \quad \forall \vtest \in \Vbh(\Omegah),
\end{align*}
where
\begin{align*}
  \LRa{L\qb,\vtest} &:= 
   \LRp{\tilde{A}\qb, \Grad\vtest}_\Omegah 
  - \LRa{\LRp{\tilde{A}\qb}^*\cdot\nb, \vtest}_\pOmegah, \nonumber \\
  \LRa{\mc{N}(\qb),\vtest} &:=
  -\LRp{\tilde{A}\qb - \Fcal(\qb) , \Grad\vtest}_\Omegah 
   + \LRa{ 
     \LRp{ \LRp{\tilde{A}\qb}^* - \Fcal^*(\qb) }\cdot\nb, \vtest}_\pOmegah.
\end{align*}
For this work, we use the Roe flux \cite{Roe81} for both linear and nonlinear fluxes, i.e.,
\begin{align*}
  \Fcal^* (\qb^\pm) & = \average{\Fcal(\qb)} + \half |A(\qb^*_{Roe})| \jump{\nb\cdot \qb},\\
  \LRp{\tilde{A}\qb^\pm}^* & = \average{\tilde{A}\qb} 
  + \half |A(\tilde{\qb}^*_{Roe})| \jump{\nb\cdot \qb},
\end{align*}
where $\qb^*_{Roe}$ and $\tilde{\qb}^*_{Roe}$ are Roe average states\footnote{
The Roe average state for $\qb^*_{Roe}$, for example, is defined as 
$\rho^*_{Roe} = \sqrt{\rho^- \rho^+}$,
${\ub}^*_{Roe} = \frac{\sqrt{\rho^-}\ub^- 
+ \sqrt{\rho^+}\ub^+ }{\sqrt{\rho^-} + \sqrt{\rho^+}}$,
$H^*_{Roe} = \frac{\sqrt{\rho^-}H^- 
+ \sqrt{\rho^+}H^+ }{\sqrt{\rho^-} + \sqrt{\rho^+}}$,
and $a^*_{Roe} = \sqrt{(\gamma-1) (H^*_{Roe} - \half \norm{\ub^*_{Roe}}^2 )} $. 
}, and $|A|:= R |\Lambda| R^{-1}$ with $\LRp{\Lambda,R}$ as the eigen-pairs of $A$ (see, e.g  
\cite{masatsuka2013like}, for more details).

\subsubsection{Artificial viscosity}

Solving nonlinear Euler equations is a notoriously challenging task. 
A smooth solution can turn into a discontinuous one due to nonlinearity. 
For numerical stability, a sufficient numerical diffusion needs to be equipped with the Euler system. 
The questions are how to measure the regularity of a solution, 
and how to choose a reasonable amount of the artificial viscosity.
The authors in \cite{guermond2011entropy, zingan2013implementation} introduced 
entropy viscosity to stabilize numerical solution in Runge Kutta time stepping.
Since a large entropy is produced in the vicinity of strong shocks,
the size of the entropy residual can be used to measure the solution regularity 
and the viscosity coefficient.

We employ the entropy viscosity method in the context of exponential DG methods for handling sharp gradient solutions.
The procedure is as following: 
given an entropy pair $(\Scal, \ub\Scal)$ for Euler systems, we first define the entropy residual, 
$ \Res(\qb):= \dd{\Scal(\qb)}{t} + \Div\LRp{\ub\Scal(\qb)}, $ 
an effective viscosity
$ \nu_E := c_E h_{max}^2 \max\LRp{ \snor{Res},\jump{ \ub \Scal }h^{-1} }$,
an upper bound to the viscosity 
$ \nu_{\max} := c_{\max} h_{\max} \max \snor{ \fb^\prime(\qb) } $ 
and the entropy viscosity  
$ \nu_{EV} := Smooth(\min(\nu_{max}, \nu_E))$. 
Here, 
$ \mc{S}: = \frac{\rho}{\gamma-1} \log \LRp{\frac{p}{\rho^\gamma} }$
is the physical entropy functional for Euler equations; 
$c_E$ and $c_{max}$ are tunable parameters;
$h$ is an element size;
\footnote{
  We define $h:=\frac{2 r}{k}$ with 
  $r:=\frac{1}{4|\K|}\Pi_{i=1}^3s_i$ as the radius of the circumscribed circle on the $\K$th triangular element.
}  
$\fb^\prime(\qb) = \dd{\fb}{\qb}$ is the flux Jacobian; 
and $Smooth$ is a smooth function. 
\footnote{
  For smoothing, we first compute vertex averaged entropy vicosity 
  and then linearly reconstruct the entropy vicosity on each element. 
} 
Then, we add the artifical diffusion term to Euler systems in \eqnref{euler3dd-gov}, 
which leads to 
\begin{align}
\eqnlab{euler3dd-gov-ev}
    \dd{\qb}{t} + \Div{\Fcal(\qb)} = \Div \Fcal^{EV} 
\end{align}
with $\Fcal^{EV}:= \LRp{\nu_{EV} \Grad \qb}$. 
Note that this choice of the viscous flux makes the diffusion term linear  
if we compute the entropy viscosity $\nu_{EV}$ using the solutions at current and previous timestep.
\footnote{
  In this study, we approximate the entropy residual by 
   $\Res \approx \frac{\Scal(\qb^n) - \Scal(\qb^{n-1}) }{\dt} + \half \LRp{\Div \LRp{\ub^n \Scal(\qb^n)} + \Div \LRp{ \ub^{n-1} \Scal (\qb^{n-1}) } }$.
}
Thus, the linearized flux in \eqnref{euler-split-3dd} becomes $\Fcal_L := \tilde{A} \qb - \nu_{EV} \Grad \qb$. 
Treating the linear diffusion for Euler systems is similar to that for Burgers equation, 
and hence omitted here.

\renewcommand{\un}{{\unm}}

\section{An analysis of the exponential DG method for Burgers equation}
\seclab{AppendixAnalysis}
In this Section we shall provide a rigorous analysis of the exponential DG approach in Section \secref{ExponentialDG}
for the Burger equation \eqnref{gov-burgers1d}. For the simplicity of the exposition, we assume that the integrals can be computed exactly though we use LGL quadrature for computing the integrals. Note that aliasing errors from LGL quadrature and interpolation are typically negligible for well-resolved solutions or can be made vanished by using a split form of the flux (see, e.g., \cite{bui2012analysis,gassner2018br1}).
For brevity, we assume zero boundary conditions\footnote{At the Dirichlet boundary, we take $\u^+ = \u_D = 0$ for the numerical fluxes $\LRp{\tilde{\u}\u}^*$ and $\LRp{\frac{\u^2}{2}}^*$, $\u^{**} = \u_D = 0$, and $\q^{**} = \q^-$.} on both sides of the domain $\Omega$ or periodic boundary conditions. Without any ambiguity, we also neglect the dependency of the (semi-discrete and exact) solutions on time $t$, e.g. $\u = \u\LRp{t}$, except for cases where this dependence is important. We assume that $\kappa > 0$ is a constant.
The stability is trivial (see, e.g. \cite{gassner2018br1}, and the references therein), thanks to the entropy numerical flux.  (Note that unlike the standard entropy numerical flux, ours has an additional jump term.) Indeed, by 
 taking $\p = \kappa \q$ in \eqnref{gov-burgers1d-gradu} and $\v = \u$ in \eqnref{gov-burgers1d-u}, and then adding the resulting equations together we obtain
\begin{equation}
\eqnlab{preStability}
\half\td{}{t}\nor{\u}_\Oh^2  + k\nor{q}_\Oh^2 = \LRp{\frac{1}{6}\pp{\u^3}{x},1}_\Oh -  \LRa{\LRp{\frac{\u^2}{2}}^*,\jump{\u}}_\Gh = - \frac{\sigma}{h}\nor{\jump{\u}}_\Gh^2.    
\end{equation}
Let us define the DG-norm\footnote{Recall our convention that on the boundary faces $\Ghb$ we have $\jump{\u} = \u$. For periodic boundary condition, $\nor{\u}_{DG}$ in Lemma \lemref{BurgerStability} and Theorem \theoref{spatialError} treats the boundary and interior interfaces the same.} for $H^1\LRp{\Oh}$ as follows
\begin{equation}
\eqnlab{DGnorm}
\nor{\u}_{DG}^2 = \nor{\pp{\u}{x}}_\Oh^2 + \frac{1}{h}\nor{\jump{\u}}^2_\Gh. 
\end{equation}
\begin{lem}[Semi-discrete stability and uniqueness]
\lemlab{BurgerStability}
The semi-discretization with DG in \eqnref{gov-burgers1d}, with $\u \in \Vh\LRp{\Omega_h}$, is stable in the following sense
\[
\nor{\u}_\Oh^2 + C\int_0^t\nor{\u}_{DG}^2\,d\tau \le \nor{\u\LRp{0}}_\Oh^2,
\]
where $C$ is some positive constant and $\u\LRp{0}$ is the DG initial condition. Hence, the DG solution is unique.
\end{lem}
\begin{proof}
From an inverse  and a multiplicative trace inequalities we have
\begin{equation}
\eqnlab{ITinequality}
h\nor{\averageL{\pp{\u}{x}}}_\Gh^2 \le \cIT\nor{\pp{\u}{x}}_\Oh^2, \quad \forall \pp{\u}{x} \in \Vh\LRp{\Omega_h},
\end{equation}
where $\cIT$ in a constant independent of the meshsize $h$.
Taking $\p = \pp{\u}{x}$ in \eqnref{gov-burgers1d-gradu} we obtain
\begin{multline*}
   \frac{c_1}{2}\nor{\q}_\Oh^2 + \frac{1}{2c_1} \nor{\pp{\u}{x}}_\Oh^2\ge \LRp{\q,\pp{\u}{x}}_\Oh = \nor{\pp{\u}{x}}_\Oh^2 - \LRa{\jump{\u},\averageL{\pp{\u}{x}}}_\Gh  \ge \\
   \nor{\pp{\u}{x}}_\Oh^2 - \frac{c_2}{2h}\nor{\jump{\u}}_\Gh^2 - \frac{h}{2c_2}\nor{\averageL{\pp{\u}{x}}}^2_\Gh \ge \LRp{1 - \frac{\cIT}{2c_2}}\nor{\pp{\u}{x}}_\Oh^2 - \frac{c_2}{2h}\nor{\jump{\u}}_\Gh^2,
\end{multline*}
where we have used \eqnref{ITinequality} in the third inequality, and $c_1,c_2$ are arbitrary positive constant. It follows that
\[
\frac{c_1}{2}\nor{\q}_\Oh^2 \ge \LRp{1 - \frac{\cIT}{2c_2} - \frac{1}{2c_1}}\nor{\pp{\u}{x}}_\Oh^2 - \frac{c_2}{2h}\nor{\jump{\u}}_\Gh^2,
\]
which, together with \eqnref{preStability}, yields
\[
\td{}{t}\nor{\u}_\Oh^2 + \frac{2\kappa}{c_1}\LRp{2 - \frac{\cIT}{c_2} - \frac{1}{c_1}}\nor{\pp{\u}{x}}_\Oh^2 + \frac{2}{h}\LRp{\sigma - \frac{\kappa c_2}{c_1}}\nor{\jump{\u}}_\Gh^2 \le 0.
\]
Now, by choosing $c_1,c_2$ large enough such that
\[
\sigma - \frac{\kappa c_2}{c_1} > 0, \text{ and } 2 - \frac{\cIT}{c_2} - \frac{1}{c_1} > 0,
\]
and defining
\[
C := \min\LRc{\sigma - \frac{\kappa c_2}{c_1},2 - \frac{\cIT}{c_2} - \frac{1}{c_1}},
\]
we arrive at
\[
\td{}{t}\nor{\u}_\Oh^2 + C\nor{\u}^2_{DG} \le 0,
\]
which concludes the proof.
\end{proof}
\begin{rema}
We can choose $c_1$ sufficiently large relative to $\kappa c_2$ so that $\sigma$ can be chosen to be (very) small. 
\end{rema}

Note that the $L^2$-stability in Lemma \lemref{BurgerStability} does not imply $L^\infty$-stability in general. The fact that $\Vh\LRp{\Omega_h}$ is piecewise continuous implies\footnote{In fact, by an inverse inequality and shape regularity, we can obtain the estimate $\nor{\u}_\infty \le C h^{-1}\nor{\u}_\Oh$.} $\Vh\LRp{\Omega_h} \subset L^\infty\LRp{\Oh}$, i.e., $\nor{\u}_\infty < \infty$. For the convergence analysis, we assume that $\nor{u}_\infty$ is bounded uniformly for the time horizon $\LRp{0,T}$, i.e., there exists $M < \infty$ so that $\nor{u}_\infty \le M$ at any $t \in \LRp{0,T}$. Let us denote by $\uh$ and $\qh = \pp{\uh}{x}$ the exact solution and its gradient, and we assume that $\uh, \pp{\uh}{t} \in H^s\LRp{\Omega}$ with $s > 3/2$ for $t \in \LRp{0,T}$. By the Sobolev embedding theorem, $\uh$ and $\qh$ are continuous and without loss of generality we assume $\nor{\uh}_\infty \le M$. It is easy to see that $\uh$ and $\qh$ satisfy the DG weak form \eqnref{gov-burgers1d}. Let $\Pi$ be the $L^2$-projection onto $\Vh\LRp{\Omega_h}$ and let us define
\begin{align*}
    \veps_\u &:= \uh - \u = \underbrace{\uh - \Pi\uh}_{=: \vepsI_\u} + \underbrace{\Pi\uh - \u}_{=:\vepsh_\u} = \vepsI_\u + \vepsh_\u, \\
    \veps_\q &:= \qh - \q = \underbrace{\qh - \Pi\qh}_{=: \vepsI_\q} + \underbrace{\Pi\qh - \q}_{=:\vepsh_\q} = \vepsI_\q + \vepsh_\q.
\end{align*}
Since both the exact and the DG solutions satisfy the DG weak form \eqnref{gov-burgers1d}, we subtract their corresponding equations, take $\v = \vepsh_\u$ and $\p = \kappa \vepsh_\q$, and add the resulting equations altogether to obtain
\begin{multline}
    \half \td{}{t}\nor{\vepsh_\u}_\Oh^2 + \kappa\underbrace{\nor{\vepsh_\q}_\Oh^2}_{\E_1} = \underbrace{-\kappa\LRp{\vepsI_\q,\vepsh_\q}_\Oh + \kappa\LRp{\pp{\vepsI_\u}{x},\vepsh_\q}_\Oh - \kappa\LRa{\jump{\vepsI_\u},\averageM{\vepsh_\q}}_\Gh}_{\E_{2}} \\
    \underbrace{-\kappa\LRp{\vepsI_\q,\pp{\vepsh_\u}{x}}_\Oh + \kappa\LRa{\averageM{\vepsI_\q},\jump{\vepsh_\u}}_\Gh}_{\E_{3}^a} \underbrace{-\LRp{\pp{\vepsI_\u}{t},\vepsh_\u}_\Oh}_{\E_{4}} \\
    + \underbrace{\LRp{\frac{\uh^2 - \u^2}{2},\pp{\vepsh_\u}{x}}_\Oh  - \LRa{\LRp{\frac{\uh^2}{2}}^* - \LRp{\frac{\u^2}{2}}^*,\jump{\vepsh_\u}}_\Gh}_{\E_5}.  
    \eqnlab{errorEquation}
\end{multline}
\begin{lem}[Estimate for $\E_1$]
\lemlab{E1}
Assume that the mesh $\Oh$ is regular. There holds:
\begin{align*}
    \E_1 \ge \frac{1}{c_1}\LRp{2 - \frac{1}{c_1} - \frac{\cIT}{c_2}}\nor{\pp{\vepsh_\u}{x}}_\Oh^2 - \frac{c_2}{h c_1}\nor{\jump{\vepsh_\u}}^2_\Gh 
     + \underbrace{\frac{2}{c_1}\LRp{\pp{\vepsI_\u}{x}-\vepsI_\q,\pp{\vepsh_\u}{x}}_\Oh-\frac{2}{c_1}\LRa{\jump{\vepsI_\u},\averageL{\pp{\vepsh_\u}{x}}}_\Gh}_{\E_{3}^b}
\end{align*}
\end{lem}
\begin{proof}
Since both the exact and DG solutions satisfy \eqnref{gov-burgers1d-gradu}, taking $\p = \pp{\vepsh_\u}{x}$ in \eqnref{gov-burgers1d-gradu} yields
\begin{multline*}
    \frac{c_1}{2}\nor{\vepsh_\q}_\Oh^2 + \frac{1}{2c_1} \nor{\pp{\vepsh_\u}{x}}_\Oh^2\ge \LRp{\vepsh_\q,\pp{\vepsh_\u}{x}}_\Oh = \nor{\pp{\vepsh_\u}{x}}_\Oh^2 - \LRa{\jump{\vepsh_\u},\averageL{\pp{\vepsh_\u}{x}}}_\Gh  \\- \LRp{\vepsI_\q,\pp{\vepsh_\u}{x}}_\Oh + \LRp{\pp{\vepsI_\u}{x},\pp{\vepsh_\u}{x}}_\Oh - \LRa{\jump{\vepsI_\u},\averageL{\pp{\vepsh_\u}{x}}}_\Gh. 
\end{multline*}
Now using Cauchy-Schwarz 
we have
\[
\LRa{\jump{\vepsh_\u},\averageL{\pp{\vepsh_\u}{x}}}_\Gh \le \frac{\cIT}{2c_2}\nor{\pp{\vepsh_\u}{x}}_\Oh^2 + \frac{c_2}{2h}\nor{\jump{\vepsh_\u}}_\Gh^2,
\]
where $\cIT$, independent of $h$, is constant resulting from the inverse  and multiplicative trace inequalities \eqnref{ITinequality}. Now combining the two inequalities ends the proof.
\end{proof}
\begin{lem}[Estimate for $\E_{2}$]
\lemlab{E2}
Assume that the mesh $\Oh$ is regular. There holds:
\begin{align*}
    \frac{\E_2}{\kappa} \le  \frac{\cq_1}{2}\nor{\vepsI_\q}_\Oh^2 + \frac{\cq_2}{2}\nor{\pp{\vepsI_\u}{x}}_\Oh^2  +  \frac{\cu_1}{2h}\nor{\jump{\vepsI_\u}}_\Gh^2 + \half\LRp{\frac{1}{\cq_1}+ \frac{1}{\cq_2} + \frac{\cIT}{\cu_1}}\nor{\vepsh_\q}^2_\Oh,
\end{align*}
where $\cIT$, independent of $h$, is constant resulting from an inverse  and a multiplicative trace inequalities, while $\cq_1,\cq_2$ and $\cu_1$ are any positive constants.
\end{lem}
\begin{proof}
The proof is straightforward using an Cauchy-Schwarz, an inverse, and a multiplicative trace inequalities.
\end{proof}
\begin{lem}[Estimate for $\E_3^a - \E_3^b/2$]
\lemlab{E3}
Assume that the mesh $\Oh$ is regular. There holds:
\begin{multline*}
    \E_3^a-\E_3^b/2 \le  \LRp{\frac{\cu_2\LRp{\kappa - 2}}{2} + \frac{\kappa\cIT\cu_3}{2}}\nor{\vepsI_\q}_\Oh^2 + \frac{\cu_4}{2c_1^2}\nor{\pp{\vepsI_\u}{x} - \vepsI_\q}_\Oh^2  +   
    \frac{\cu_5}{2hc_1^2}\nor{\jump{\vepsI_\u}}_\Gh^2
    \\+
    \LRp{\frac{\kappa - 2}{2\cu_2} + \frac{1}{2\cu_4} + \frac{\cIT}{2\cu_5}}\nor{\pp{\vepsh_\u}{x}}_\Oh^2
    +\frac{\kappa}{2h\cu_3}\nor{\jump{\vepsh_\u}}^2_\Gh,
\end{multline*}
where $\cIT$, independent of $h$, is constant resulting from an inverse  and a multiplicative trace inequalities, while $\cu_2,\cu_3,\cu_4$ and $\cu_5$ are any positive constants.
\end{lem}
\begin{proof}
The proof is straightforward using an Cauchy-Schwarz, and an inverse and a multiplicative trace inequality similar to \eqnref{ITinequality}.
\end{proof}
\begin{lem}[Estimate for $\E_4$]
\lemlab{E4}
Assume that the mesh $\Oh$ is regular. There holds:
\[
   \E_4 \le \frac{\cu_6}{2}\nor{\pp{\vepsI_\u}{t}}_\Oh^2 + \frac{1}{2\cu_6}\nor{\vepsh_\u}^2_\Oh.
\]
where $\cu_6$ is any positive constant.
\end{lem}

\begin{lem}[Estimate for $\E_5$]
\lemlab{E5}
Assume that the mesh $\Oh$ is regular, $\nor{\uh}_\infty \le M$, and $\nor{\u}_\infty \le M$ over the time horizon of interest. There holds:
\begin{align*}
   \E_5 \le \frac{M\cIT\cu_8}{2}\LRp{\nor{\vepsI_\u}_\Oh^2 + h^2 \nor{\pp{\vepsI_\u}{x}}_\Oh^2} + 
    \frac{\sigma\cu_7}{2h}\nor{\jump{\vepsI_\u}}_\Gh^2 
     + \frac{M}{\cu_8}\nor{\vepsh_\u}_{DG}^2 + \frac{M\cIT\cu_8}{2}\nor{\vepsh_\u}_\Oh^2 +\frac{\sigma}{h}\LRp{\frac{1}{2\cu_7} - 1}\nor{\jump{\vepsh_\u}}_\Gh^2,
\end{align*}
where $\cu_7$ and $\cu_8$ are any positive constant.
\end{lem}
\begin{proof}
We begin with the Lipschitz continuity of the flux function and the numerical flux
\begin{align*}
    \half\LRp{\uh^2 - \u^2} &\le M\snor{\veps_\u}, \text{ and } 
    \LRp{\frac{\uh^2}{2}}^* - \LRp{\frac{\u^2}{2}}^* 
    \le M\LRp{\snor{\veps_\u^+} + \snor{\veps_\u^-}} + \frac{\sigma}{h}\jump{\veps_\u},
\end{align*}
and thus by Cauchy-Schwarz inequality we have
\begin{multline*}
    \E_5 \le \underbrace{M\LRp{\snor{\vepsI_\u},\snor{\pp{\vepsh_\u}{x}}}_\Oh  + M\LRa{\snor{\vepsI_{\u^+}} + \snor{\vepsI_{\u^-}},\snor{\jump{\vepsh_\u}}}_\Gh}_{B\LRp{\vepsI_\u,\vepsh_\u}} \\
    +\underbrace{M\LRp{\snor{\vepsh_\u},\snor{\pp{\vepsh_\u}{x}}}_\Oh  + M\LRa{\snor{\vepsh_{\u^+}} + \snor{\vepsh_{\u^-}},\snor{\jump{\vepsh_\u}}}_\Gh}_{B\LRp{\vepsh_\u,\vepsh_\u}} \\
    + \frac{\sigma\cu_7}{2h}\nor{\jump{\vepsI_\u}}_\Gh^2 + 
    \frac{\sigma}{h}\LRp{\frac{1}{2\cu_7} - 1}\nor{\jump{\vepsh_\u}}_\Gh^2,
\end{multline*}
where we have defined $\vepsh_{\u^\pm} := \LRp{\vepsh_\u}^\pm$. 

To estimate $B\LRp{\vepsI_\u,\vepsh_\u}$ we apply Cauchy-Schwarz inequality and a multiplicative trace inequality to obtain 
\begin{multline*}
     B\LRp{\vepsI_\u,\vepsh_\u} \le M\LRp{\nor{\vepsI_\u}_\Oh^2 + h\nor{\vepsI_\u}_\Gh^2}^\half\nor{\vepsh_\u}_{DG} \\ \le \frac{M\cu_8}{2}\LRp{\nor{\vepsI_\u}_\Oh^2 + h\nor{\vepsI_\u}_\Gh^2} + \frac{M}{2\cu_8}\nor{\vepsh_\u}_{DG}^2 \\ \le
     \frac{M\cIT\cu_8}{2}\LRp{\nor{\vepsI_\u}_\Oh^2 + h^2 \nor{\pp{\vepsI_\u}{x}}_\Oh^2} + \frac{M}{2\cu_8}\nor{\vepsh_\u}_{DG}^2.
\end{multline*}
Similarly, together with an inverse inequality, we have
\begin{align*}
    B\LRp{\vepsh_\u,\vepsh_\u} \le \frac{M\cu_8}{2}\LRp{\nor{\vepsh_\u}_\Oh^2 + h\nor{\vepsh_\u}_\Gh^2} + \frac{M}{2\cu_8}\nor{\vepsh_\u}_{DG}^2 
     \le
     \frac{M\cIT\cu_8}{2}\nor{\vepsh_\u}_\Oh^2 + \frac{M}{2\cu_8}\nor{\vepsh_\u}_{DG}^2. 
\end{align*}
Combining the above estimates concludes the proof.
\end{proof}

Now combining Lemmas \lemref{E1}--\lemref{E5} and \eqnref{errorEquation} we arrive at
\begin{multline}
   \td{}{t}\nor{\vepsh_\u}_\Oh^2 + C_1\nor{\vepsh_\q}^2_\Oh + C_2\nor{\vepsh_\u}_{DG}^2 \le M\cIT\cu_8\LRp{\nor{\vepsI_\u}_\Oh^2 + h^2 \nor{\pp{\vepsI_\u}{x}}_\Oh^2} \\
   + \LRp{\kappa\cq_1 + \cu_2\LRp{\kappa - 2} + \kappa\cIT\cu_3+2\cu_4}\nor{\vepsI_\q}_\Oh^2   
     + \LRp{\kappa\cq_2+ 2\cu_4}\nor{\pp{\vepsI_\u}{x}}_\Oh^2
      \\+  \cu_6\nor{\pp{\vepsI_\u}{t}}_\Oh^2+ 
    \LRp{\kappa\cu_1 + \cu_5 + \sigma\cu_7}\frac{1}{h}\nor{\jump{\vepsI_\u}}_\Gh^2 
    + C_3\nor{\vepsh_\u}^2_\Oh,
    \eqnlab{errorFinal}
\end{multline}
where we have defined
\begin{align*}
    C_1 &:= {\frac{2\kappa}{c_1}\LRp{2 - \frac{1}{c_1} - \frac{\cIT}{c_2}} - \frac{1}{\cq_1}- \frac{1}{\cq_2} - \frac{\cIT}{\cu_1}}, \\
    C_2^1 &:= \LRp{\frac{\kappa}{2}-\frac{\kappa - 2}{\cu_2} - \frac{2}{\cu_4} - \frac{\cIT}{\cu_5}-\frac{2M}{\cu_8}}, \\
    C_2^2 &:= \sigma\LRp{1-\frac{1}{\cu_7} }-\frac{\kappa}{\cu_3}-\frac{2M}{\cu_8}-  \frac{2 \kappa c_2}{c_1}, \\
    C_2 &:= \min\LRc{C_2^1,C_2^2}, \\
    C_3 &:= \LRp{\frac{1}{\cu_6}+ M\cIT\cu_8}.
\end{align*}
As can be seen, for a given $\kappa$, we can choose $c_1, c_2$, $\sigma$, $\cq_1, \cq_2$, and $\cu_i, i=1,\hdots,8$ 
such that 
all the constants $C_1, C_2$ and $C_3$ are positive.
\begin{theo}[Semi-discrete error estimate]
\theolab{spatialError}
 Assume $\uh, \pp{\uh}{t} \in H^s\LRp{\Omega}$ with $s > 3/2$ for $t \in \LRs{0,T}$, and $\sigma > 0$. There exist positive constants $C$ independent of the meshsize $h$ and $t$ such that 
\begin{align}
\td{}{t}\nor{\vepsh_\u}_\Oh^2 + \nor{\vepsh_\q}^2_\Oh + \nor{\vepsh_\u}_{DG}^2 \le
C\nor{\vepsh_\u}^2_\Oh 
+ C h^{2\min\LRc{s,k}}\LRp{\nor{\uh}^2_{H^s\LRp{\Omega}} + \nor{\pp{\uh}{t}}^2_{H^s\LRp{\Omega}}},
\eqnlab{error1}
\end{align}
and thus, in addition, if $\uh, \pp{\uh}{t} \in L^2\LRp{\LRp{0,T};H^s\LRp{\Omega}}$ and $\u\LRp{0} = \Pi\uh\LRp{0}$, there exists a constant $C$ independent of the meshsize $h$ and $t$ such that
\begin{align*}
\nor{\veps_\u}^2_\Oh + \int_0^t\LRp{\nor{\veps_\q}^2_\Oh + \nor{\veps_\u}_{DG}^2}\,ds  \le 
C h^{2\min\LRc{s,k}}\LRp{\exp\LRp{Ct}-1}.
\end{align*}
\end{theo}
\begin{proof}
The first assertion is the direct consequence of the error estimation \eqnref{errorFinal}, the $L^2$-projection error \cite{Babuska71,BabuskaSuri87,BabuskaSuri94,ErnGuermond04,Ciarlet02}, and the following definition of $C$:
\[
C := \max\LRc{\frac{C_3}{\min\LRc{1,C_1,C_2}},\frac{C_4}{\min\LRc{1,C_1,C_2}}}.
\]
The second assertion is straightforward by 1) integrating the first assertion and then applying a Gronwall's lemma to obtain,
\begin{align*}
\nor{\vepsh_\u}^2_\Oh + \int_0^t\LRp{\nor{\vepsh_\q}^2_\Oh + \nor{\vepsh_\u}_{DG}^2}\,ds 
\le 
Ch^{2\min\LRc{s,k}}\LRp{\exp\LRp{Ct}-1}\int_0^t\LRp{\nor{\uh}^2_{H^s\LRp{\Omega}} + \nor{\pp{\uh}{t}}^2_{H^s\LRp{\Omega}}}\,ds.
\end{align*}
and 2) using the the $L^2$-projection error \cite{Babuska71,BabuskaSuri87,BabuskaSuri94,ErnGuermond04,Ciarlet02} and triangle inequalities for $\veps_\u$ and $\veps_\q$, e.g.,
\[
\nor{\veps_\u}_\Oh \le \nor{\vepsI_\u}_\Oh + \nor{\vepsh_\u}_\Oh.
\]
\end{proof}
\begin{rema}
Theorem \theoref{spatialError} shows that though the convergence of the solution in the DG norm \eqnref{DGnorm} is optimal,
 the convergence in $L^2$-error is suboptimal, 
 i.e. $\nor{\veps_\u}_\Oh = \mc{O}\LRp{h^k}$ if $k \le s$,
  and this seems to be sharp as we observe this rate in the diffusion-dominated numerical results especially for odd $k$. 
  When $Ct \ll 1$, then $\exp\LRp{Ct}-1 \approx Ct$, and thus the error increases at most linearly in time. 
\end{rema}

We next analyze the temporal discretization error using the exponential integrator. We begin with a few important lemmas.
\begin{lem}
\lemlab{qBound}
There exists a constant $C$, independent of the meshsize $h$, such that:
\[
    h\nor{\average{\q}}_\Gh^2 + \nor{\q}_\Oh^2 \le C\nor{u}_{DG}^2
\]
\end{lem}
\begin{proof}
Taking $\p = \q$ in \eqnref{gov-burgers1d-gradu}, then using Cauchy-Schwarz and multiplicative trace and inverse inequality similar to \eqnref{ITinequality} give
\begin{align*}
 \frac{h}{4\cIT}\nor{\average{\q}}_\Gh^2 + \frac{3}{4}\nor{\q}_\Oh^2 \le \nor{\q}_\Oh^2 = \LRp{\pp{\u}{x},\q}_\Oh - \LRa{\jump{\u},\average{\q}}_\Gh 
  \le \half\nor{\pp{\u}{x}}_\Oh^2 + \half\nor{\q}_\Oh^2 + \frac{2\cIT}{h}\nor{\jump{\u}}_\Gh^2 + \frac{h}{8\cIT}\nor{\average{\q}}_\Gh^2,    
\end{align*}
then taking $C = \max\LRc{\half,2\cIT}/\min\LRc{\frac{1}{8\cIT},\frac{1}{4}}$ concludes the proof.
\end{proof}

Using the definition of the $L^2$-projection, e.g., $\Pi\mc{N}\LRp{\u} \in \Vh\LRp{\Oh}$ such that $\LRa{\mc{N}\LRp{\u},\v} = \LRp{\Pi\mc{N}\LRp{\u},\v}_\Oh$ for all $\v \in \Vh\LRp{\Oh}$, we can write \eqnref{gov-burgers1d-u} as
\[
\pp{u}{t} = \Pi L\u + \Pi\mc{N}\LRp{\u}.
\]
For the clarity of the exposition let us define $\PL := \Pi L$ and $\PN := \Pi\mc{N}$.
For the rest of the analysis, we do not distinguish the operator
$\PL\u$ (and hence $\PN\LRp{\u}$ ) and its matrix representation from
$\Vh\LRp{\Oh}$ to $\Vh\LRp{\Oh}$ since in finite dimension all norms
are equivalent and $\Vh\LRp{\Oh}$ is homeomorphic to $\R^M$, where $M$
is the dimension of $\Vh\LRp{\Oh}$. This allows us to work conveniently and directly on $\PL$, $\PN$, $\Vh\LRp{\Oh}$, and the DG-norm.

\begin{lem}[Uniform boundedness of $L$ and $\PL$ on $\Vh\LRp{\Oh}$]
\lemlab{Lboundedness}
Suppose that $\nor{\tilde{\u}}_\infty \le M$, the linear operator $\PL$ defined in \eqnref{gov-burgers1d-u}, as a linear operator from $\Vh\LRp{\Oh}$ to $\Vh\LRp{\Oh}$, is bounded for any $t\in \LRp{0,T}$ in the following sense: there exists a constant $C_L$, independent of the meshsize $h, \u, \v$, 
such that
\[
\snor{\LRp{\PL\u,\v}_\Oh} = \snor{\LRa{L\u,\v}} \le C_L\nor{\u}_{DG}\nor{\v}_{DG},\quad \forall \u,\v \in \Vh\LRp{\Oh}.
\]
\end{lem}
\begin{proof}
We have
\begin{align*}
    \LRa{\n\LRp{\tilde{u}u}^*,\v}_\pOmegah & \le M\LRa{\snor{\u},\snor{\jump{\v}}}_\pOmegah, \text{ and } \LRa{\n\kappa\q^{**},\v}_\pOmegah \le \kappa\LRa{\snor{\average{\q}},\snor{\jump{\v}}}_\Gh, \\
    \LRp{\kappa\q-\tilde{\u}\u,\pp{\v}{x}}_\Oh &\le \kappa\LRp{\snor{\q},\snor{\pp{\v}{x}}}_\Oh + \kappa M\LRp{\snor{\u},\snor{\pp{\v}{x}}}_\Oh.
\end{align*}
Now from definition of $L$ in \eqnref{BurgerL} using Cauchy-Schwarz inequality  we obtain
\begin{align*}
    \LRa{L\u,\v}^2 \le \LRp{\frac{\kappa M}{2}\nor{\u}_\Oh^2 + \frac{\beta h}{2}\nor{\u}_\pOmegah^2 + \frac{\kappa}{2}\nor{\q}_\Oh^2 + \frac{\kappa h}{2}\nor{\average{\q}}_\Gh^2} \times 
    \LRp{\frac{\kappa M}{2}\nor{\pp{v}{x}}_\Oh^2 + \frac{\beta}{h}\nor{\jump{\v}}_\Gh^2 + \frac{\kappa}{2}\nor{\pp{\v}{x}}_\Oh^2 + \frac{\kappa}{2h}\nor{\jump{\v}}_\Gh^2},
\end{align*}
Now using the result of Lemma \lemref{qBound}, a Poincar\'e-Friedrichs inequality for $H^1\LRp{\Oh}$, and a multiplicative trace inequality for $\nor{\u}_\Gh$ we arrive at
\begin{align*}
    \LRa{L\u,\v}^2 \le \LRs{C_{PF}\LRp{\frac{\kappa M}{2} + \frac{\beta\cIT}{2}}+\frac{\kappa C}{2}}\nor{\u}_{DG}^2 \times 
    \max\LRc{\frac{\kappa M}{2} + \frac{\kappa}{2},\beta + \frac{\kappa}{2}}\nor{\v}_{DG}^2,
\end{align*}
where $C_{PF}$ is the constant in the Poincar\'e-Friedrichs inequality. The result follows by taking 
\[
C_L^2 := \LRs{C_{PF}\LRp{\frac{\kappa M}{2} + \frac{\beta\cIT}{2}}+\frac{\kappa C}{2}} \times \max\LRc{\frac{\kappa M}{2} + \frac{\kappa}{2},\beta + \frac{\kappa}{2}}.
\]
\end{proof}
\begin{coro}
\corolab{semigroup}
Suppose the assumptions for Lemma \lemref{Lboundedness} hold. We have that
\[
e^{t\PL} := \sum_{m=0}^\infty \frac{\LRp{t \PL}^m}{m!}, \quad t \ge 0,
\]
is a uniform continuous operator semigroup in $\Vh\LRp{\Oh}$. In particular, $\PL$ is the infinitesimal generator of the semigroup $e^{t\PL}$ with
\[
\nort{e^{t\PL}}\le e^{C_L t} \text{ and } \nort{\varphi_i\LRp{\tau \PL}} \le C_{\varphi_i}\LRp{\tau} := e^{\tau C_L} \int_0^1\frac{z^{i-1}}{\LRp{i-1}!}\,dz,
\]
where $\nort{\cdot}$ is the operator norm from $\Vh\LRp{\Oh}$ to  $\Vh\LRp{\Oh}$.
\end{coro}
\begin{proof}
The proof is straightforward using the boundedness of $\PL$ in Lemma \lemref{Lboundedness} \cite{engel2001one,hochbruck2010exponential}.
\end{proof}
\begin{lem}[Lipschitz continuity of $\mc{N}$]
\lemlab{Nlipschitz}
Let $\u$ and $\w$ be in $\Vh\LRp{\Oh}$, $\nor{\u}_\infty \le M$ and $\nor{\w}_\infty \le M$ for $t \in \LRp{0,T}$, there exists a constant $C_{\mc{N}}$ independent of the meshsize such that
\[
\snor{\LRa{\mc{N}\LRp{\u} - \mc{N}\LRp{\w},\v}} \le C_{\mc{N}}\nor{\u -\w}_{DG}\nor{\v}_{DG}, \quad \forall \v \in \Vh\LRp{\Oh}.
\]
\end{lem}
\begin{proof}
From the definition of $\mc{N}$ in \eqnref{BurgerN}, it is easy to see that
\begin{align*}
    \LRp{\frac{u^2}{2}}^* - \LRp{\frac{w^2}{2}}^* &\le \frac{M}{3}\LRp{\snor{\u^- - \w^-} + \snor{\u^+ - \w^+}} + \frac{\sigma}{h}\jump{\u-\w}, \\
    \LRp{\tilde{u}\u}^* - \LRp{\tilde{u}\w}^* &\le 2M \LRp{\snor{\u^- - \w^-} + \snor{\u^+ - \w^+}}, \\
    \LRp{\tilde{u} \u -\half u^2} - \LRp{\tilde{u}\w - \half \w^2} &\le 2M\snor{\u-\w}. \\
\end{align*}
Thus,
\begin{multline*}
\snor{\LRa{\mc{N}\LRp{\u} - \mc{N}\LRp{\w},\v}} \le 2M\nor{\u -\w}_\Oh\nor{\pp{\v}{x}}_\Oh + \frac{7M}{3}\nor{\u-\w}_\pOmegah\nor{\jump{\v}}_\Gh   \\
+\frac{\sigma}{h}\nor{\jump{\u-\w}}_\Gh\nor{\jump{\v}}_\Gh\le \LRs{M^2\LRp{4 + \frac{98}{9}\cIT}\nor{\u-\w}_\Oh^2 + \frac{2\sigma^2}{h}\nor{\jump{\u-\w}}^2_\Gh}^\half\nor{\v}_{DG} \\
\le \underbrace{\LRs{C_{PF}M^2\LRp{4 + \frac{98}{9}\cIT} + 2\sigma^2}^\half}_{C_{\mc{N}}}\nor{\u -\w}_{DG}\nor{\v}_{DG},
\end{multline*}
where we have used Cauchy-Schwarz, inverse trace, and Poincar\'e-Friedrichs inequalities.
\end{proof}

A direct consequence of Lemmas \lemref{Lboundedness} and \lemref{Nlipschitz} is that $\pp{u}{t} = \PL\u + \PN\LRp{\u} \in \Vh\LRp{\Oh} \subset L^\infty\LRp{\Oh}$. 
Thus from \eqnref{BurgerN}, $\pp{\PN}{t}$ given by
\begin{align*}
\LRa{\pp{\PN\LRp{\u}}{t},\v} := \LRp{\LRp{\u - \ut}\pp{\u}{t},\pp{\v}{t}}_\Oh + \LRa{\averageL{\ut\pp{\u}{t}} + \half\max(|\tilde{u}^\pm|)\jumpL{\pp{\u}{t}},\jump{\v}}_\Gh 
-\LRa{\frac{1}{3}\averageL{\u\pp{\u}{t}} + \frac{2}{3}\average{\u}\averageL{\pp{\u}{t}} + \frac{\sigma}{h}\jumpL{\pp{\u}{t}},\jump{\v}}_\Gh
\end{align*}
is well-defined. Indeed the next results show that $\pp{\PN}{t}$ is a Lipschitz continuous map from $\Vh\LRp{\Oh}$ to $\Vh\LRp{\Oh}$.
\begin{lem}[Lipschitz continuity of $\pp{\PN}{t}$]
\lemlab{Nplipschitz}
Suppose $\u,\pp{\u}{t}, \ut,\w$ reside in $\Vh\LRp{\Oh}$, $\nor{\u}_\infty \le M$, $\nor{\w}_\infty \le M$, $\nor{\ut}_\infty \le M$, and $\norL{\pp{\u}{t}}_\infty \le M$ for $t \in \LRp{0,T}$. There exists a constant $C'_{\mc{N}}$ independent of the mesh size $h$ such that
\[
\snor{\LRa{\pp{\PN\LRp{\u}}{t} - \pp{\PN\LRp{\w}}{t},\v}} \le C'_{\mc{N}}\nor{\u-\w}_{DG}\nor{\v}_{DG}
\]
\end{lem}
\begin{proof}
Let us define $\z$ as
\[
\z := \pp{\u}{t}-\pp{\w}{t} = \PL\LRp{\u-\w} + \PN\LRp{\u} - \PN\LRp{\w},
\]
and Lemmas \lemref{Lboundedness} and \lemref{Nlipschitz} imply $\nor{\z}_{DG} \le \LRp{C_L + C_{\mc{N}}}\nor{\u-\w}_{DG}$. We have
\begin{multline*}
    \LRa{T,\v} = \LRp{\LRp{\u-\w}\pp{\u}{t} + \LRp{\w-\ut}\z,\pp{\v}{x}}_\Oh 
    + \LRa{\average{\ut\z} + \half\max(|\tilde{u}^\pm|)\jump{\z},\jump{\v}}_\Gh
     \\
    -\LRa{\frac{1}{3}\averageL{\LRp{\u-\w}\pp{\u}{t} + \w\z} + \frac{2}{3}\average{\u-\w}\averageL{\pp{\u}{t}} + \frac{2}{3}\average{\w}\average{\z} +\frac{\sigma}{h}\jump{\z},\jump{\v}}_\Gh,
\end{multline*}
where $T := \pp{\PN\LRp{\u}}{t} - \pp{\PN\LRp{\w}}{t}$. Now following similar arguments as in the proof of Lemma \lemref{Nlipschitz} we obtain
\[
\LRa{T,\v} \le C\nor{\z}_{DG}\nor{\v}_{DG} \le C\LRp{C_L+ C_{\mc{N}}}\nor{\u-\w}_{DG}\nor{\v}_{DG},
\]
where $C = C\LRp{\cIT, M,\sigma,C_{PF}}$.
\end{proof}

Let us denote by $\u^n$ the approximation solution of \eqnref{gov-burgers1d} using a time discretization. We are now in the position to analyze error of the fully discrete system using the exponential-Euler in time and DG in space. The fully discrete system using the exponential Euler integrator reads
\[
\u^{n+1} := e^{\dt \PL}\u^n + \dt \varphi_1\LRp{\dt \PL} \PN\LRp{\u^n},
\]
while the semi-discrete solution $\u\LRp{t^n}$ satisfies
\[
\u\LRp{t^{n+1}} = e^{\dt \PL}\u\LRp{t^n} + \dt \varphi_1\LRp{\dt \PL} \PN\LRp{\u\LRp{t^n}} + \delta^{n+1},
\]
with
\[
\delta^{n+1} := \int_0^\dt e^{\LRp{\dt -\tau}L}\int_0^\tau\pp{\PN}{t}\LRp{\u\LRp{t^n+\theta}}\,d\theta\,d\tau.
\]
Let us define 
\footnote{    
    $\rho^0=0$.
} 
\[
\xin := \uh\LRp{t^n} - \u^n =  \underbrace{\uh\LRp{t^n} - \u\LRp{t^n}}_{\veps_\u\LRp{t^n}} +  \underbrace{\u\LRp{t^n} - \u^n}_{=:\rhon
} = \veps_\u\LRp{t^n} + \rhon.
\]
We thus have
\begin{align*}
    \rho^{n+1} &=  e^{\dt \PL}\rhon + \dt \varphi_1\LRp{\dt \PL} \underbrace{\LRs{\PN\LRp{\u\LRp{t^n}} - \PN\LRp{\u^n}}}_{=:f\LRp{t^n}} + \delta^{n+1}. 
\end{align*}
After some algebraic manipulations
we obtain
    \begin{align*}
        \rho^{n} &=  \dt\sum_{j=0}^{n-1}e^{\LRp{n-j-1}\dt \PL}\varphi_1\LRp{\dt \PL}f\LRp{t_j} + \sum_{j=0}^{n-1}e^{j\dt \PL}\delta^{n-j}. 
    \end{align*}
    
    \begin{theo}[Convergence of the exponential Euler-DG]
      \theolab{temporalError}
      Assume the conditions of Corollary \cororef{semigroup}, Lemma \lemref{Nplipschitz}, Lemma \lemref{Nlipschitz}, Lemma \lemref{BurgerStability}, and Theorem \theoref{spatialError} hold.
      There exists a constant $C$ depending only $t^n, C_L,
      C_{\mc{N}}, \nor{\u\LRp{0}}_\Oh, C_{PF}$ such that the following
      estimate for the total
      discrete error $\xi^n$ at $t^n = n\dt$ holds true:
      \[
        \nor{\xi^n}_{\Oh} \le 
        C\LRp{\dt + h^{\min\LRc{s,k}}}.
      \]
\end{theo}
\begin{proof}
Using Corollary \cororef{semigroup}, Lemma \lemref{Nplipschitz}, and Lemma \lemref{BurgerStability} we have
\begin{align*}
    \nor{\sum_{j=0}^{n-1}e^{j\dt \PL}\delta^{n-j}}_{DG} \le \dt e^{t^n C_L}C'_{\mc{N}} \sum_{j=1}^n\int_{t^{j-1}}^{t^{j}}\nor{\u\LRp{t}}_{DG}\,dt 
    \le \dt \sqrt{t^n}e^{t^n C_L}C'_{\mc{N}}
    \sqrt{\int_{0}^{t^n}\nor{\u\LRp{t}}_{DG}^2\,dt} \le \dt 
    \underbrace{C\nor{\u\LRp{0}}_\Oh}_{C_\delta}
\end{align*}
where we have absorbed all quantities depending only on $t^n, C_L, C_{\mc{N}}$ into $C$. 
Now combining the above estimate with Lipschitz continuity of $\PN$
(Lemma \lemref{Nlipschitz}) and Corollary \cororef{semigroup} we have
that there exists a constant $C_1$ depending on
$C_{\varphi_1}\LRp{\dt}$ and $C_\delta$ such that
\begin{align*}
\nor{\rho^{n}}_{DG} &\le  C_1\dt \sum_{j=0}^{n-1}\nor{\rho^j}_{DG} + C_1\dt, 
\end{align*}
which, together with a discrete Gronwall's lemma, yield
\[
\nor{\rho^n}_\Oh \le C_{PF}\nor{\rho^{n}}_{DG} \le  C \dt, 
\]
where $C = C\LRp{C_{PF}, C_1}$.
We thus, via Theorem \theoref{spatialError}, conclude
\[
\nor{\xi^n}_{\Oh} \le \nor{\rho^{n}}_{\Oh} + \nor{\veps_\u\LRp{t^{n}}}_{\Oh} \le C\LRp{\dt + h^{\min\LRc{s,k}}}.
\]
\end{proof}


\section{Numerical Results}
\seclab{NumericalResults}

  
In this section, we conduct several numerical experiments 
to evaluate the performance of the proposed exponential DG framework for both Burgers and Euler equations.
In particular, we examine the numerical stability for a wide range of Courant numbers\footnote{
	We define Courant numbers $Cr:=Cr_a + Cr_d$, where 
	$Cr_a:= \frac{c \dt}{dx}$ for convection part 
	and $Cr_d:= \frac{\kappa \dt}{(dx)^2}$ for diffusion part. 
	Here, $dx$ is the minimum distance between two LGL nodes; $c$ is the maximum speed in the system, e.g., $c=|u|$ for Burgers equation 
	and $c=|\ub \cdot \nb| + a$ for Euler equations. 
} larger than unity, i.e. $Cr>1$,
the high-order convergence in both space and time, the efficiency, and weak and strong parallel scalability.
We measure $L^2$ error for convergence studies by $\norm{u -\hat{u}}_\Omegah$ where $\hat{u}$ is either an exact solution or a reference solution.

\subsection{Viscous Burgers equation}
\subsubsection{An exact time-independent smooth solution}
We consider a time-independent manufactured solution with $\kappa=0.03$ for the Burgers equation
$$ u(x,t) = \sin(x^2) x (x-1) $$ 
by adding the corresponding source term to \eqnref{pde-burgers-eq}. 
We perform 
a spatial convergence study with both Lax-Friedrich (LF) flux and entropy flux (EF) (we take $\sigma=0$ and $\sigma=3\times 10^{-4}$). 
In order to prevent temporal discretization error from polluting the spatial one, we employ high-order accurate EXPRB32 scheme
\footnote{
EXPRB32 \cite{hochbruck2010exponential} is the third-order two-stage exponential method given by
\begin{align*}
  q^{(2)} &= q^n + \dt \varphi_1(L \dt) R^n,\\
  q^{n+1} &= q^n + \dt \varphi_1(L \dt) R^n
                 + 2 \dt \varphi_3(L \dt) D_{n,2},
\end{align*}
where $R(q)=Lq + \mc{N}(q)$, $R^n:=R(q^n)$ and $D_{n,2}:= \mc{N}(q^{(2)}) - \mc{N}(q^n)$. 
Here, we choose 
$L=\dd{R}{q}\big\vert_{\q^n}$.
}
with $\dt=5\times 10^{-5}$ and $N_e=\LRc{20,40,80,160}$ elements. 
The error is measured at $t=0.01$ and the results are summarized in Table \tabref{expo-burgers-sconv-mms}.
We observe that entropy flux with no additional stabilization (i.e. $\sigma = 0$) and Lax-Friedrichs flux provide similar results, that is, the convergence order is optimal for even solution orders but sub-optimal for odd ones. This similar behavior when using the central flux for diffusion term has been recorded in the literature (see, e.g. \cite{cockburn1998local}, and the references therein) and it is also consistent with our analysis in Section \secref{AppendixAnalysis} in which we have shown that the spatial convergence is sub-optimal. Entropy flux with small additional stabilization seems to asymptotically deliver convergence rates between $k+1/2$ and $k+1$ for all solution orders considered in this case, which is better than what we could prove.



\begin{table}[t]
\caption{A time-independent manufactured solution for the viscous Burgers equation:
a spatial convergence study using $N_e=\LRc{20,40,80,160}$ elements is conducted with Lax-Friedrich (LF) flux,
 entropy flux (EF). EXPRB32 scheme with $\dt=5\times10^{-5}$ is used as the time integrator.
}
\tablab{expo-burgers-sconv-mms}
\begin{center}
\begin{tabular}{*{1}{c}|*{1}{c}|*{2}{c}|*{2}{c}|*{2}{c}}
\hline
\multirow{2}{*}{ } 
  & \multirow{2}{*}{$h$}  
  & \multicolumn{2}{c|}{$LF$}   
  & \multicolumn{2}{c|}{$EF$ ($\sigma=3\times 10^{-4}$)}   
  & \multicolumn{2}{c|}{$EF$ ($\sigma=0$)}   
  \tabularnewline
  &    &  error & order &  error & order &  error & order       \tabularnewline
\hline\hline
\multirow { 4}{*}{ $k=1$ }&  1/20  & 4.093E-04 &      $-$ & 4.096E-04 &      $-$ & 4.097E-04 &      $-$ \tabularnewline
                          &  1/40  & 1.223E-04 &    1.743 & 1.225E-04 &    1.741 & 1.232E-04 &    1.734 \tabularnewline
                          &  1/80  & 4.494E-05 &    1.445 & 4.378E-05 &    1.484 & 4.620E-05 &    1.415 \tabularnewline
                          &  1/160 & 1.937E-05 &    1.214 & 1.562E-05 &    1.487 & 2.074E-05 &    1.155 \tabularnewline
\multicolumn{8}{c}{} \tabularnewline
\multirow { 4}{*}{ $k=2$ }&  1/20  & 2.630E-06 &      $-$ & 2.632E-06 &      $-$ & 2.634E-06 &      $-$ \tabularnewline
                          &  1/40  & 3.210E-07 &    3.034 & 3.213E-07 &    3.034 & 3.222E-07 &    3.031 \tabularnewline
                          &  1/80  & 3.966E-08 &    3.017 & 3.952E-08 &    3.023 & 3.996E-08 &    3.011 \tabularnewline
                          &  1/160 & 4.916E-09 &    3.012 & 4.833E-09 &    3.032 & 4.984E-09 &    3.003 \tabularnewline
\multicolumn{8}{c}{} \tabularnewline
\multirow { 4}{*}{  $k=3$}&  1/20  & 1.431E-07 &      $-$ & 1.459E-07 &      $-$ & 1.471E-07 &      $-$ \tabularnewline
                          &  1/40  & 1.709E-08 &    3.066 & 1.742E-08 &    3.066 & 1.819E-08 &    3.015 \tabularnewline
                          &  1/80  & 2.003E-09 &    3.093 & 1.876E-09 &    3.215 & 2.268E-09 &    3.003 \tabularnewline
                          &  1/160 & 2.251E-10 &    3.154 & 1.441E-10 &    3.703 & 2.834E-10 &    3.001 \tabularnewline
\multicolumn{8}{c}{} \tabularnewline
\multirow { 4}{*}{ $k=4$} &  1/20  & 5.946E-10 &      $-$ & 5.991E-10 &      $-$ & 6.013E-10 &      $-$ \tabularnewline
                          &  1/40  & 1.827E-11 &    5.024 & 1.837E-11 &    5.027 & 1.860E-11 &    5.015 \tabularnewline
                          &  1/80  & 5.626E-13 &    5.021 & 5.590E-13 &    5.039 & 5.796E-13 &    5.004 \tabularnewline
                          &  1/160 & 1.730E-14 &    5.023 & 1.690E-14 &    5.048 & 1.810E-14 &    5.001 \tabularnewline
\hline\hline
\end{tabular}
\end{center}
\end{table}

\subsubsection{A smooth solution}
\seclab{burgers-smoothsol}
We next consider a case with smooth solution generated by the following initial condition, 
\begin{align*}
    u(x,t=0) = \sin^3(2\pi x) (1-x)^{\frac{3}{2}}, 
\end{align*}
 zero Dirichlet boundary conditions, and 
$\kappa=0.03$.
The smooth initial profile is spread out due to the viscosity as time goes by.

We conduct 
a spatial convergence study with both LF and EF fluxes with 
$N_e=\LRc{20,40,80,160}$. 
We again use EXPRB32 for the time integrator with $\dt=5\times 10^{-5}$.
Since there is no exact solution we use RK4 solution with $\dt=5\times 10^{-7}$,  $k=10$, and
$N_e=160$ as the "ground truth" solution.
The error at $t=0.01$ is used to compute the convergence rate
 and the results are summarized in Table \tabref{expo-burgers-sconv-smooth}.
Similar to the case of the manufactured solution, 
we observe that the convergence rates for odd solution orders are sub-optimal, 
but optimal for even solution orders. Again, a little additional stabilization via $\sigma$ not only facilitates our convergence analysis but also seems to improve the convergence rates.

\begin{table}[t]
\caption{A smooth solution for the viscous Burgers equation:
a spatial convergence study using  $N_e=\LRc{20,40,80,160}$ elements is conducted with Lax-Friedrich (LF) flux and entropy flux (EF), and
EXPRB32 scheme with $\dt=5\times 10^{-5}$ as the time integrator.}
\tablab{expo-burgers-sconv-smooth}
\begin{center}
\begin{tabular}{*{1}{c}|*{1}{c}|*{2}{c}|*{2}{c}|*{2}{c}}
\hline
\multirow{2}{*}{} 
  & \multirow{2}{*}{$h$}  
  & \multicolumn{2}{c|}{$LF$}   
  & \multicolumn{2}{c|}{$EF (\sigma=3\times 10^{-4})$}   
  & \multicolumn{2}{c|}{$EF (\sigma=0)$}   
  \tabularnewline
  &    &  error & order &  error & order &  error & order       \tabularnewline
\hline\hline
\multirow { 4}{*}{   $k=1$ }&  1/20 & 7.061E-03 &      $-$ & 7.157E-03 &      $-$ & 7.161E-03 &      $-$ \tabularnewline
                         &  1/40 & 2.080E-03 &    1.764 & 2.285E-03 &    1.647 & 2.308E-03 &    1.633 \tabularnewline
                         &  1/80 & 6.731E-04 &    1.627 & 8.552E-04 &    1.418 & 9.186E-04 &    1.329 \tabularnewline
                         &  1/160 & 2.432E-04 &    1.469 & 3.070E-04 &    1.478 & 4.241E-04 &    1.115 \tabularnewline
\multicolumn{8}{c}{} \tabularnewline
\multirow { 4}{*}{   $k=2$ }&  1/20 & 2.511E-04 &      $-$ & 2.540E-04 &      $-$ & 2.541E-04 &      $-$ \tabularnewline
                         &  1/40 & 2.741E-05 &    3.196 & 2.744E-05 &    3.210 & 2.744E-05 &    3.211 \tabularnewline
                         &  1/80 & 3.317E-06 &    3.047 & 3.317E-06 &    3.048 & 3.317E-06 &    3.048 \tabularnewline
                         &  1/160 & 4.112E-07 &    3.012 & 4.113E-07 &    3.012 & 4.112E-07 &    3.012 \tabularnewline
\multicolumn{8}{c}{} \tabularnewline
\multirow { 4}{*}{   $k=3$ }&  1/20 & 2.574E-05 &      $-$ & 2.738E-05 &      $-$ & 2.748E-05 &      $-$ \tabularnewline
                         &  1/40 & 2.883E-06 &    3.158 & 3.243E-06 &    3.078 & 3.311E-06 &    3.053 \tabularnewline
                         &  1/80 & 3.414E-07 &    3.078 & 3.748E-07 &    3.113 & 4.087E-07 &    3.018 \tabularnewline
                         &  1/160 & 4.108E-08 &    3.055 & 3.722E-08 &    3.332 & 5.091E-08 &    3.005 \tabularnewline
\multicolumn{8}{c}{} \tabularnewline
\multirow { 4}{*}{   $k=4$ }&  1/20 & 9.225E-07 &      $-$ & 9.273E-07 &      $-$ & 9.273E-07 &      $-$ \tabularnewline
                         &  1/40 & 2.529E-08 &    5.189 & 2.785E-08 &    5.057 & 2.858E-08 &    5.020 \tabularnewline
                         &  1/80 & 7.135E-10 &    5.147 & 7.563E-10 &    5.202 & 8.445E-10 &    5.081 \tabularnewline
                         &  1/160 & 2.188E-11 &    5.027 & 2.190E-11 &    5.110 & 2.604E-11 &    5.019 \tabularnewline
\hline\hline
\end{tabular}
\end{center}
\end{table}

Before showing temporal convergence rates let us demonstrate the numerical stability
 of exponential integrators by using very large Courant numbers. For this purpose, it is sufficient to choose EPI2 scheme\footnote{
EPI2 \cite{tokman2006efficient} is the second-order exponential method given by
\begin{align*}
    q^{n+1} 
    &= q^n + \dt \varphi_1(L\triangle t) R^n
\end{align*}
where $L=\dd{R}{q}\big\vert_{\q^n}$. 
}. We compare EPI2 and RK2 solutions at $t=1$, both with EF flux, in Figure \figref{pde-burgers-smooth}. 
For RK2 solution, 
we take $\dt=10^{-4}$ ($Cr_d=0.16$) as approximately the maximal stable timestep size since $\dt=2\times10^{-4}$ leads to unstability. Unlike RK2,
with 1000 times larger timestep size, i.e. $Cr_d=161.0$, EPI2 still
produces a comparable result compared to RK2.
Even with 5000 times larger timestep size ($Cr_d=804.9$), 
EPI2 solution is still stable though less accurate (see Figure \figref{pde-burgers-ex-smooth-difference} in which we compare the accuracy of EPI2 and RK2 using RK4 with $\dt=5\times 10^{-6}$ as a reference solution).
As for the wallclock time, RK2 takes $8.3s$, whereas EPI2 does $1.6s$ with $Cr_d=161.0$ and $1.3s$ with $Cr_d=804.9$. EPI2 is five to six times faster than RK2 in this example.  

\begin{figure}[h!t!b!]
    \centering
    \subfigure[solutions ]{
      \includegraphics[trim=0.0cm 0.0cm 0.0cm 0cm,clip=true,width=0.4\textwidth]{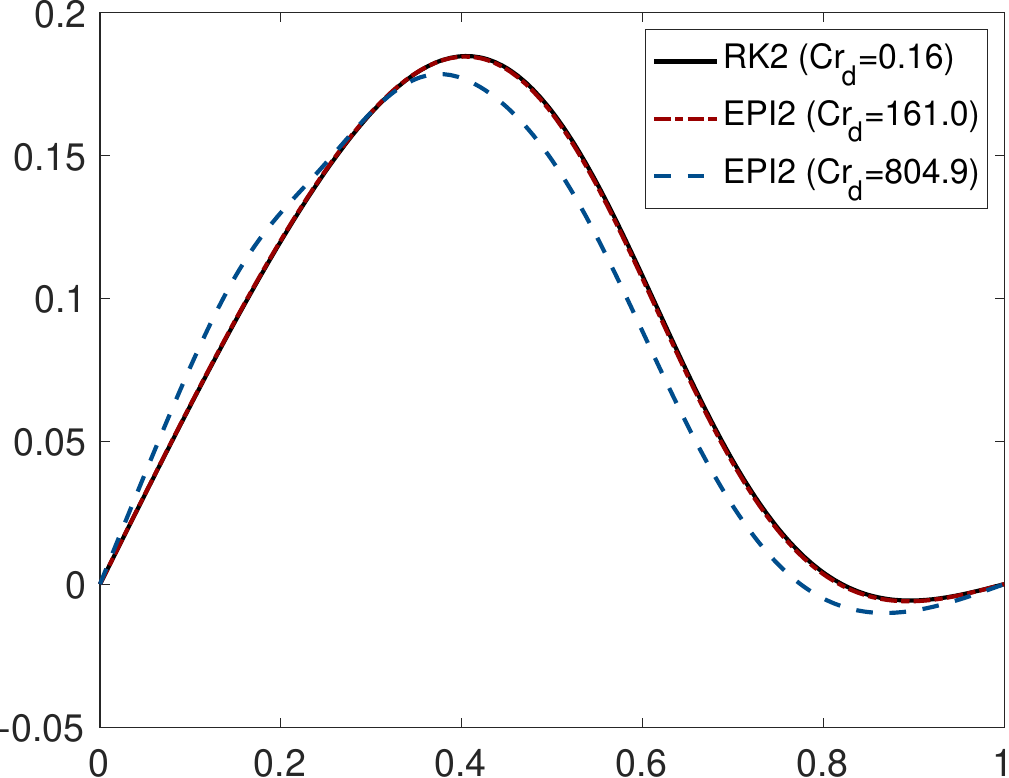}
      \figlab{pde-burgers-ex-smooth-all}
    }
    \subfigure[differences ]{
      \includegraphics[trim=0.0cm 0.0cm 0.0cm 0cm,clip=true,width=0.4\textwidth]{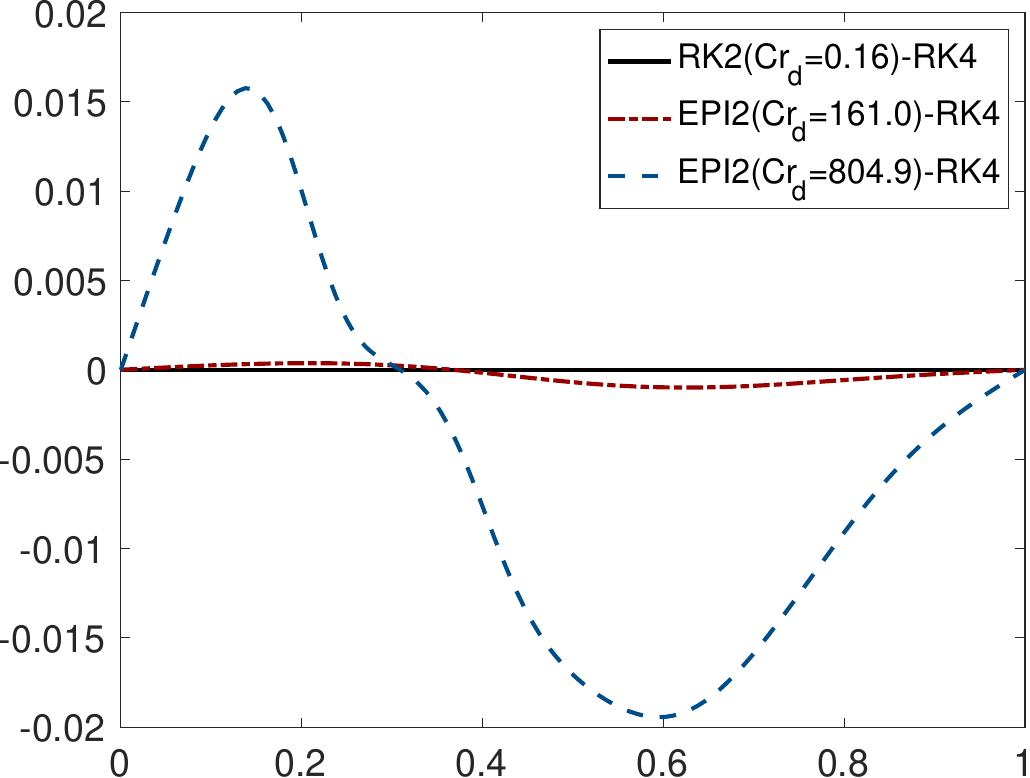}
      \figlab{pde-burgers-ex-smooth-difference}
    }
    \caption{Accuracy comparisons between EPI2 and RK2 using Burgers equation with smooth solutions at $t=1$ with $k=4$ and $N_e=40$: Figure \figref{pde-burgers-ex-smooth-all} shows RK2 with $Cr_d=0.16$ (black), EPI2 with $Cr_d=161.0$ (red-dashed), and EPI2 with $Cr_d=804.9$ (blue-dashed).  Figure  \figref{pde-burgers-ex-smooth-difference} plots their differences with RK4 solution with $dt=5\times 10^{-6}$.} 
    \figlab{pde-burgers-smooth}
  \end{figure}

We now compute the temporal convergence rates of two exponential integrators, EPI2 and EXPRB32, using our DG spatial discretization  with $k=4$ and $N_e=40$. 
To that end, we take the RK4 solution  with $\dt=5\times 10^{-6}$ as a ground truth.
The error is computed at $t=1$. 
As can be seen in Table \tabref{expo-burgers-tconv-smooth}, 
the numerical results with both LF and EF fluxes  
show second- and third-order convergence rates
for EPI2 and EXPRB32, respectively. 
We observe that the difference in the solutions of LF and EF are negligibly small (on the order of $\mc{O}(10^{-9})$). 
This, we believe, is due to the diffusion-dominated regime, for which different numerical fluxes for the nonlinear convection term do not make (much) difference on the solution. 

\begin{table}[h!t!b!]
\caption{Temporal convergence rates for EPI2 and EXPRB32 for a smooth solution of the viscous Burgers equation by computing the error at $t=1$ using RK4 solution with $\dt=5\times 10^{-6}$ as the reference solution.
Spatial discretization is carried out using our DG approach with $k=4$ and $N_e=40$. 
We observe second- and third-order convergence rates for EPI2
and EXPRB32, respectively.}
\tablab{expo-burgers-tconv-smooth}
\begin{center}
\begin{tabular}{c|cc||cc|cc}
\hline
\multirow{2}{*}{flux} 
& \multirow{2}{*}{$\dt$} 
& \multirow{2}{*}{$Cr_d$} 
& \multicolumn{2}{c|}{ EPI2 } 
& \multicolumn{2}{c}{ EXPRB32 } \tabularnewline
&  &  &  error & order & error & order \tabularnewline
\hline\hline
\multirow{5}{*}{LF} 
  & 0.50 & 804.9 &   1.171E-02 &      $-$ &       5.272E-03 &      $-$ \tabularnewline
  & 0.25 & 402.5 &   3.303E-03 &    1.827 &       1.077E-03 &    2.292 \tabularnewline
  & 0.10 & 161.0 &   5.411E-04 &    1.974 &       9.575E-05 &    2.641 \tabularnewline
  & 0.05 &  80.5 &   1.312E-04 &    2.044 &       1.300E-05 &    2.881 \tabularnewline
  & 0.01 &  16.1 &   4.943E-06 &    2.037 &       1.042E-07 &    2.999 \tabularnewline
  \\
\multirow{5}{*}{EF\small{($\sigma=3\times10^{-4}$)}}   
 & 0.50 & 804.9 &   1.171E-02 &      $-$ &       5.272E-03 &      $-$ \tabularnewline
 & 0.25 & 402.5 &   3.303E-03 &    1.827 &       1.077E-03 &    2.292 \tabularnewline
 & 0.10 & 161.0 &   5.411E-04 &    1.974 &       9.575E-05 &    2.641 \tabularnewline
 & 0.05 &  80.5 &   1.312E-04 &    2.044 &       1.300E-05 &    2.881 \tabularnewline
 & 0.01 &  16.1 &   4.943E-06 &    2.037 &       1.042E-07 &    2.999 \tabularnewline  
 \\
\multirow{5}{*}{EF\small{($\sigma=0$)}}    
 & 0.50 & 804.9 &   1.171E-02 &      $-$ &       5.272E-03 &      $-$ \tabularnewline
 & 0.25 & 402.5 &   3.303E-03 &    1.827 &       1.077E-03 &    2.292 \tabularnewline
 & 0.10 & 161.0 &   5.411E-04 &    1.974 &       9.575E-05 &    2.641 \tabularnewline
 & 0.05 &  80.5 &   1.312E-04 &    2.044 &       1.300E-05 &    2.881 \tabularnewline
 & 0.01 &  16.1 &   4.943E-06 &    2.037 &       1.043E-07 &    2.988 \tabularnewline
\hline\hline
\end{tabular}
\end{center}
\end{table}

\subsubsection{A solution with steep gradient}

We next consider a solution with steep gradient,
namely, a stationary shock that evolves in time from the following
initial condition 
\begin{align*}
  u(x,t=0) = \sin(2\pi x), \quad \text{ for } x \in [0,1],
\end{align*}
and  homogeneous boundary conditions.
We 
 perform the simulation for $t \in [0, 1]$ with $k = 4$ and $N_e = 40$.
 As time goes on, a sharp interface is progressively formed at $x=0.5$. 
 In this convection-dominated example, 
 EF flux with a uniform bound of $\sigma$ leads to an unstable solution, 
 whereas LF flux still produces a stable solution.
 This may be related to 
 the growth of aliasing errors arising from the sharp gradient, 
 or insufficient artificial diffusion between the elements 
 due to the lack of upwinding of the uniform value of $\sigma$. 
Inspired by LF flux, 
we set $\sigma=\kappa/100 + h\max(\snor{u^\pm})$ for EF flux in this example.
Figure \figref{pde-burgers-stationary-shock-kappa} are the EPI2 solution snapshots at $t=1$ with $dt=0.01$ and $\kappa\in \LRc{0.02,0.005,0.002}$.
By decreasing $\kappa$, the shock solutions with EF flux become stiffer, 
and the difference of EPI2 solutions between LF flux and EF flux increases up to 
$\mc{O}(10^{-4})$.


\begin{figure}[h!t!b!]
  \centering
    \subfigure[solutions]{
      \includegraphics[trim=0.0cm 0.0cm 0.0cm 0cm,clip=true,width=0.4\textwidth]{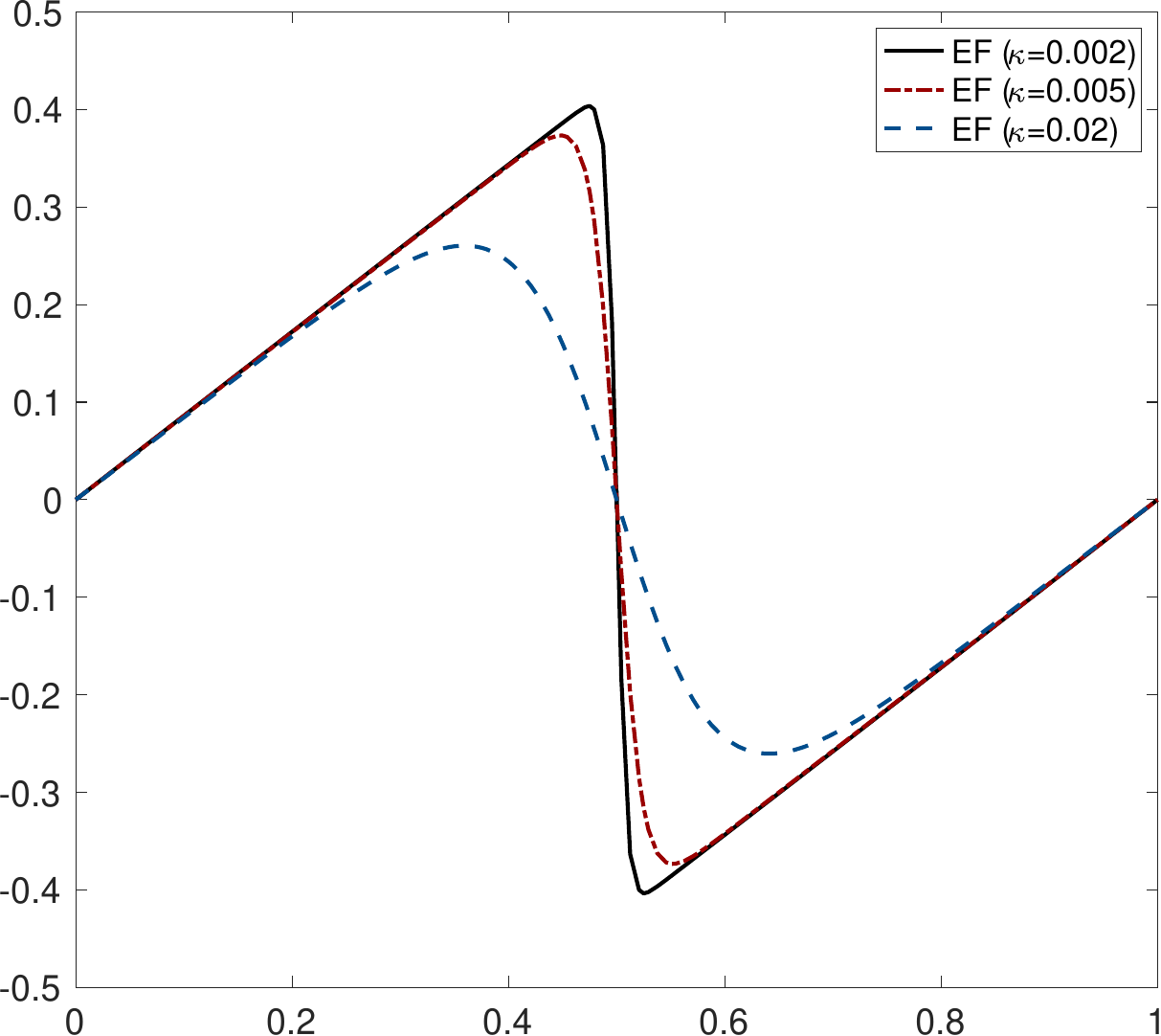}
      \figlab{pde-burgers-ex-stationary-shock-all-kappa}
    }
    \subfigure[differences]{
      \includegraphics[trim=0.0cm 0.0cm 0.0cm 0cm,clip=true,width=0.4\textwidth]{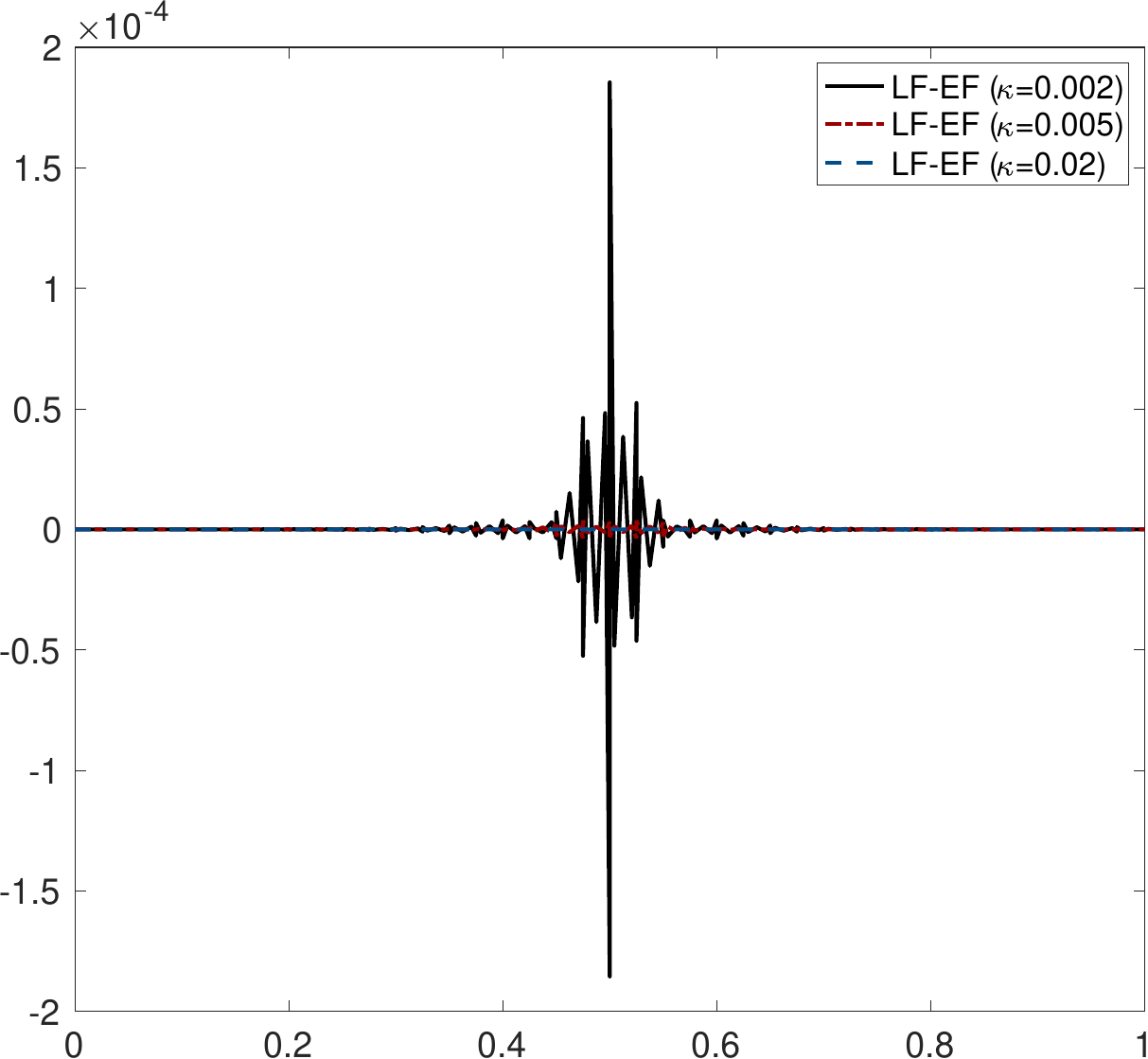}
      \figlab{pde-burgers-ex-stationary-shock-difference-kappa}
    }
    \caption{Burgers equation with a steep gradient solution with EPI2 time integrator, $k=4$, and $N_e=40$ at $t=1$: 
    (a) the snapshots with $\kappa=\LRc{0.02,0.005,0.0002}$ and 
    (b) the solution differences between $LF$ and $EF$ fluxes.}
  \figlab{pde-burgers-stationary-shock-kappa}
\end{figure} 

    

With $\kappa = 0.002$, we now perform spatial and temporal convergence studies.
For spatial convergence study, we use nested meshes with $N_e=\LRc{40,80,160}$
\footnote{
  All the meshes are chosen to align with the sharp interface at $x=0.5$.
} and the EXPRB32 integrator with $\dt=5\times 10^{-5}$.
RK4 solution (with LF flux, $\dt=5\times 10^{-7}$,  $k=10$, and
$N_e=160$) is used as the "ground truth"
solution for measuring the $L^2$ error at
$t=1$. 
We observe $k+\half$ rate of convergence for both LF and EF fluxes in Table \tabref{expo-burgers-sconv-shock}.

\begin{table}[h!t!b!]
\caption{A shock solution to the viscous Burgers equation with 
a spatial convergence study using  $N_e=\LRc{40,80,160}$ elements is conducted with Lax-Friedrich (LF) flux and entropy flux (EF), and
EXPRB32 scheme with $\dt=5\times 10^{-5}$ as the time integrator.
We take the RK4 solution (with $\dt=5\times 10^{-7}$,  $k=10$, and
$N_e=160$) as a reference solution to measure the $L^2$ error at
$t=1$. 
}
\tablab{expo-burgers-sconv-shock}
\begin{center}
\begin{tabular}{*{1}{c}|*{1}{c}|*{2}{c}|*{2}{c}}
\hline
\multirow{2}{*}{ } 
  & \multirow{2}{*}{$h$}  
  & \multicolumn{2}{c|}{$LF$}   
  & \multicolumn{2}{c|}{$EF$}   
  \tabularnewline
  &    &  error & order &  error & order    \tabularnewline
\hline\hline
\multirow { 3}{*}{ $k=1$ }&  1/40 & 1.295E-02 &      $-$ & 1.353E-02 &      $-$ \tabularnewline
                         &  1/80  & 5.558E-03 &    1.221 & 5.577E-03 &    1.279 \tabularnewline
                         &  1/160 & 2.001E-03 &    1.474 & 1.964E-03 &    1.506 \tabularnewline
\multicolumn{6}{c}{} \tabularnewline
\multirow { 3}{*}{ $k=2$ }&  1/40 & 4.223E-03 &      $-$ & 4.229E-03 &      $-$ \tabularnewline
                         &  1/80  & 6.722E-04 &    2.651 & 6.712E-04 &    2.656 \tabularnewline
                         &  1/160 & 9.852E-05 &    2.770 & 9.616E-05 &    2.803 \tabularnewline
\multicolumn{6}{c}{} \tabularnewline
\multirow { 3}{*}{ $k=3$ }& 1/40 & 9.976E-04 &      $-$ & 9.969E-04 &      $-$ \tabularnewline
                         &  1/80 & 1.435E-04 &    2.797 & 1.449E-04 &    2.782 \tabularnewline
                         &  1/160 & 1.224E-05 &    3.551 & 1.219E-05 &    3.571 \tabularnewline
\multicolumn{6}{c}{} \tabularnewline
\multirow { 3}{*}{ $k=4$ }& 1/40 & 3.606E-04 &      $-$ & 3.577E-04 &      $-$ \tabularnewline
                         &  1/80 & 2.728E-05 &    3.724 & 2.706E-05 &    3.725 \tabularnewline
                         &  1/160 & 8.084E-07 &    5.077 & 7.432E-07 &    5.186 \tabularnewline
\hline\hline
\end{tabular}
\end{center}
\end{table}

For temporal convergence study, 
  we consider the DG discretization with $k=4$, and $N_e=160$.
We take the RK4 solution (with LF flux, $\dt=5\times 10^{-7}$,  $k=10$, and
$N_e=160$) as a reference solution to  measure the $L^2$ error at $t=1$.
In Table \tabref{expo-burgers-tconv-shock}, 
we observe second- and third-order convergence rates for EPI2
and EXPRB32, respectively.

\begin{table}[h!t!b!]
\caption{Temporal convergence study of EPI2 and EXPRB32 for a shock solution to the viscous Burgers equation:
we take the RK4 solution (with LF flux, $\dt=5\times 10^{-7}$,  $k=10$, and
$N_e=160$) as a reference solution, and measure the $L^2$ error at $t=1$.
Spatial discretization is carried out using our DG approach with $k=4$ and $N_e=160$. 
We observe second- and third-order convergence rates for EPI2
and EXPRB32, respectively.}
\tablab{expo-burgers-tconv-shock}
\begin{center}
\begin{tabular}{c|cc||cc|cc}
\hline
\multirow{2}{*}{flux} 
& \multirow{2}{*}{$\dt$} 
& \multirow{2}{*}{$Cr$} 
& \multicolumn{2}{c|}{ EPI2 } 
& \multicolumn{2}{c}{ EXPRB32 } \tabularnewline
&  &  &  error & order & error & order \tabularnewline
\hline\hline
\multirow{5}{*}{$LF$} 
  & 0.25 & 559.1 &   2.276E-02 &      $-$ &       8.424E-03 &      $-$ \tabularnewline
  & 0.10 & 213.3 &   3.706E-03 &    1.981 &       5.711E-04 &    2.937 \tabularnewline
  & 0.05 & 105.5 &   8.784E-04 &    2.077 &       6.440E-05 &    3.149 \tabularnewline
  & 0.02 & 42.0  &   1.327E-04 &    2.063 &       3.755E-06 &    3.102 \tabularnewline
  \\
  \multirow{5}{*}{$EF$}   
  & 0.25 & 559.1 &   2.276E-02 &      $-$ &       8.424E-03 &      $-$ \tabularnewline
  & 0.10 & 213.3 &   3.706E-03 &    1.981 &       5.711E-04 &    2.937 \tabularnewline
  & 0.05 & 105.5 &   8.784E-04 &    2.077 &       6.440E-05 &    3.149 \tabularnewline
  & 0.02 & 42.0  &   1.327E-04 &    2.063 &       3.745E-06 &    3.105 \tabularnewline
\hline\hline
\end{tabular}
\end{center}
\end{table}

Considering the simplicity and upwinding nature of Lax-Friedrich flux, 
we use Lax-Friedrich flux for Euler equations (which are hyperbolic)  
in the following examples.

\subsection{Euler equations: Isentropic vortex translation}

We consider the isentropic vortex example in \cite{zhou2003high}, 
where a small vortex perturbation is added to the uniform mean flow 
and translated without changing its shape. 
The superposed flow is given as 
\begin{align*}
u &= u_\infty - \frac{\lambda}{2\pi} \tilde{y} e^{\alpha \LRp{1-r^2}/2}, \quad 
v  = v_\infty + \frac{\lambda}{2\pi}\tilde{x} e^{\alpha\LRp{1-r^2}/2}, \\
T & = T_\infty - \LRp{\frac{1}{2\alpha c_p} } \LRp{\frac{\lambda}{2\pi}}^2 e^{\alpha\LRp{1-r^2}},
\end{align*}
where $r = \norm{\xb - \xb_c - \ub_\infty t}$,
$\tilde{x}=x - x_c - u_\infty t$,
$\tilde{y}=y - y_c - v_\infty t$,
$\ub_\infty=(u_\infty,v_\infty)$, and 
$\xb_c=(5,0)$. Here, $c_p$ is the specific heat ratio at constant pressure 
and $\lambda$ the vortex strength. 
The mean flow is set to be $(u_\infty,v_\infty,\rho_\infty,T_\infty,\pres_\infty) = (0.2,0,1,1,1)$.
We take $\alpha=2$, $\gamma=1.4$, $\lambda=0.05$ and 
$c_p=\frac{\gamma}{\gamma-1}$.
The exact solution is generated from the isentropic relation\footnote{
$$\frac{\rho}{\rho_\infty} = \LRp{\frac{T}{T_\infty}}^{\frac{1}{\gamma-1}}
=\LRp{\frac{\pres}{\pres_\infty}}^{\frac{1}{\gamma}}
$$
}. 
The domain is $\Omega = \LRp{0,10} \times \LRp{-5,5}$ and 
periodic boundary conditions are applied to all directions. 
For a three-dimensional simulation, we take the zero vertical velocity, $w = 0$, and extrude the 2D domain vertically from $0$ to $1$ to obtain $\Omega = \LRp{0,10} \times \LRp{-5,5} \times \LRp{0,1}$. 

\subsubsection{Stability of exponential integrators}
We perform the simulation for $t\in [0,15]$ 
with $k=12$ and $N_e=256$ in Figure \figref{expo-pde-euler2d-vortex}.
Compared to the stable EPI2 solution with $\dt=0.5$ (i.e. $Cr=42.75$) in Figure \figref{pde-euler2d-epi2-Cr42-t15}, 
the EPI2 solution with $\dt=1$ ($Cr=85.56$) in Figure \figref{pde-euler2d-epi2-Cr85-t15} 
is oscillatory (part of the domain away from the vortex). 
One can reduce the oscillation while keeping large Courant number by employing a more accurate, e.g. higher-order, time integrator
\footnote{
Since high-order methods requires more nonlinear evaluations,
  the evolution of the solution can be captured more accurately than the low-order methods. 
 For example, EPI2 needs one nonlinear evaluation at $t^n$, 
 whereas EXPRB32 uses two nonlinear evaluations at $t^n$ and $t^{n+1}$.
}.
To demonstrate this point we show the third-order EXPRB32 solution with $\dt=1$ $(Cr=85.56)$
in Figure \figref{pde-euler2d-exprb32-Cr85-t15} for which oscillations are not visible in the same scale. 
A closer look, see Figure \figref{expo-pde-euler2d-vortex-y0}, 
in which we plot a slice along $y=0$ for all sub-figures in Figure \figref{expo-pde-euler2d-vortex}, shows that the EXPRB32 solution with $(Cr=85.56)$ does reduce oscillations. This is not surprising: {\em though exponential integrators are inherently implicit, large timestep size must be chosen with care in order to avoid adverse affect on the accuracy.}
\begin{figure}[h!t!b!]
\centering
 \subfigure[EPI2 ($Cr=42.75$)]{
  \includegraphics[trim=0.0cm 1.0cm 2.0cm 2cm,clip=true,width=0.3\textwidth]{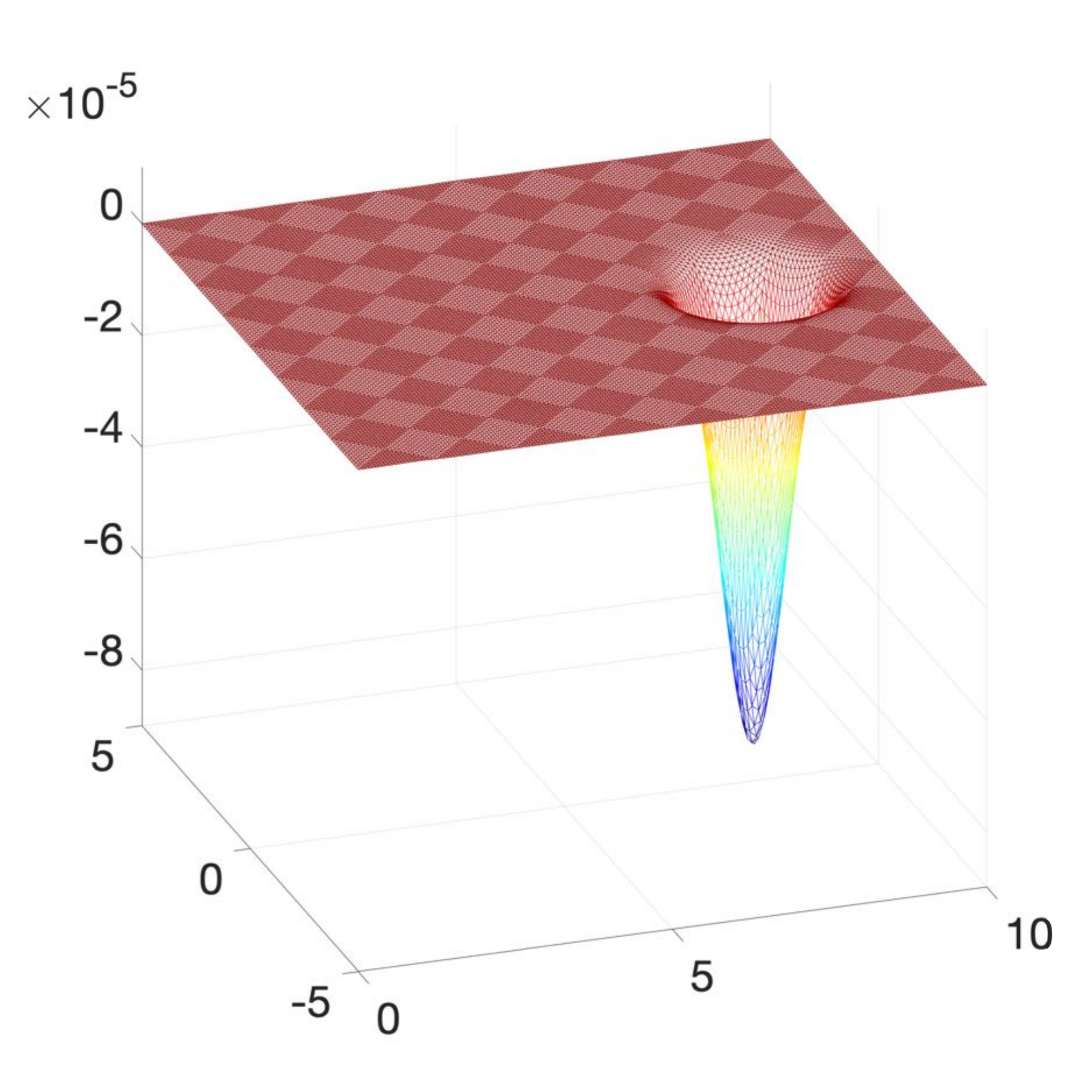}
   \figlab{pde-euler2d-epi2-Cr42-t15}
 }
 \subfigure[EPI2 ($Cr=85.56$)]{
  \includegraphics[trim=0.0cm 1.0cm 2.0cm 2cm,clip=true,width=0.3\textwidth]{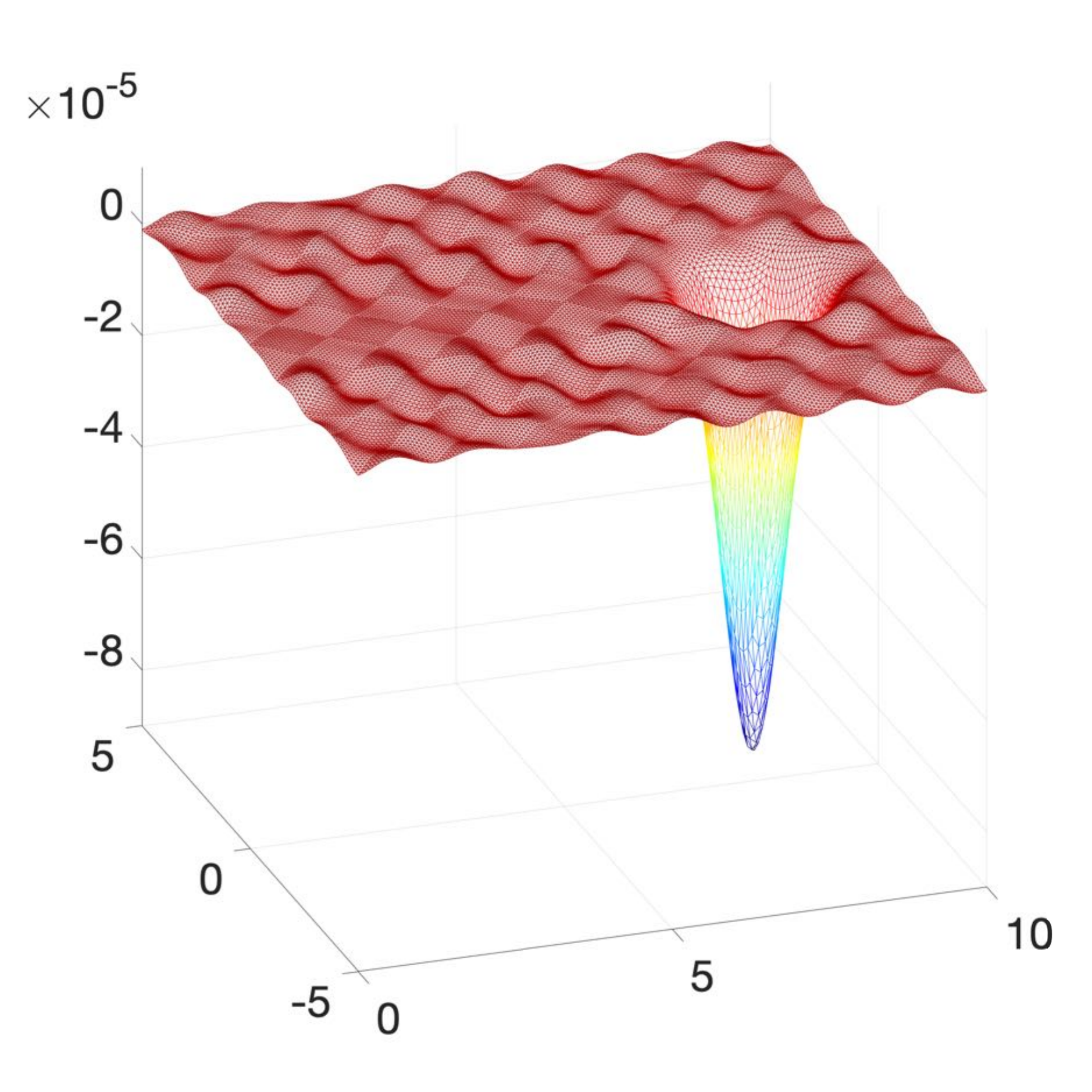} 
   \figlab{pde-euler2d-epi2-Cr85-t15}
 }    
 \subfigure[EXPRB32 ($Cr=85.56$)]{
   \includegraphics[trim=0.0cm 1.0cm 2.0cm 2cm,clip=true,width=0.3\textwidth]{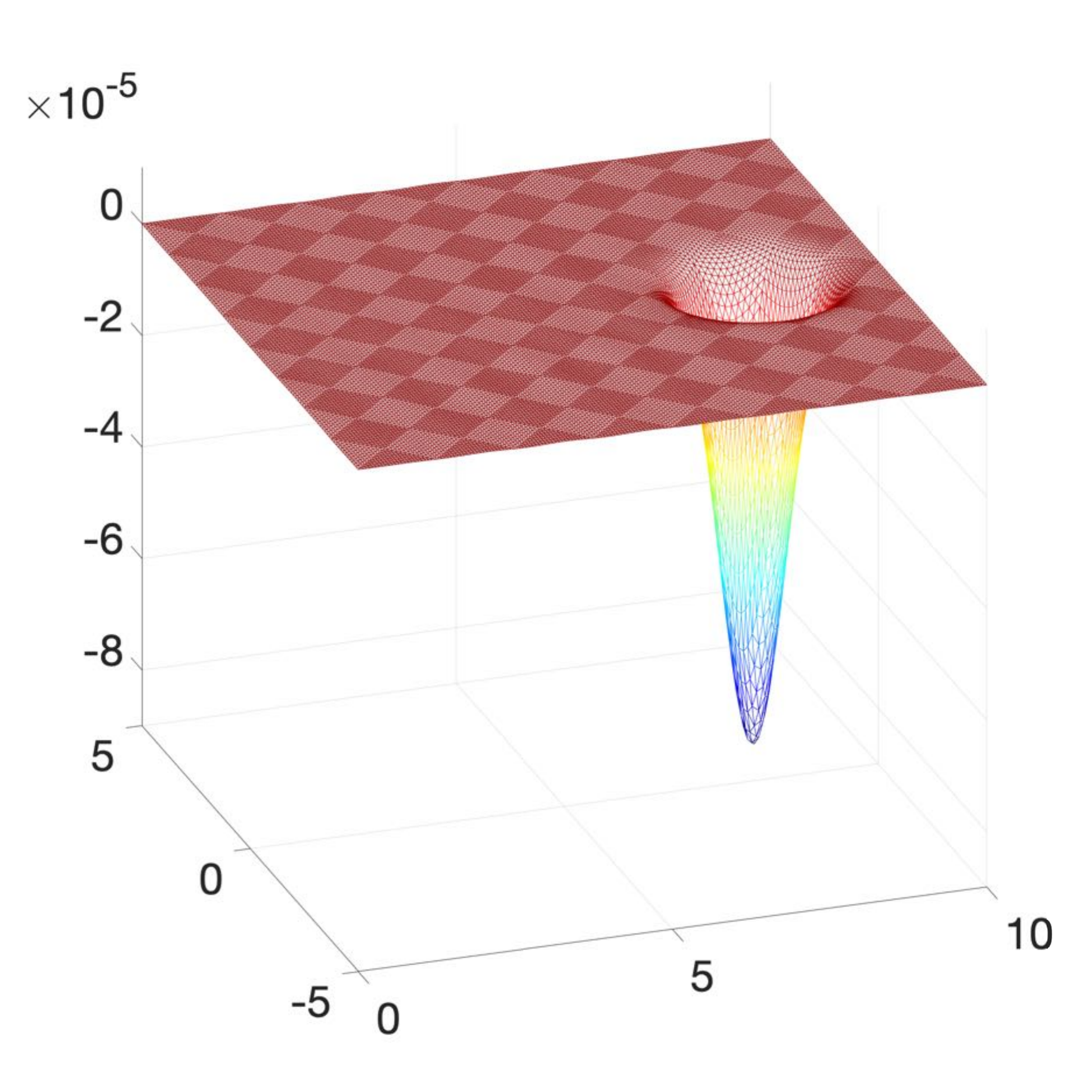}
    \figlab{pde-euler2d-exprb32-Cr85-t15}
 }
\caption{Translating isentropic vortex example for the two-dimensional Euler equations at $t=15$: 
numerical solution using EPI2 with $Cr=42.75$ is in Figure \figref{pde-euler2d-epi2-Cr42-t15}, using EPI2 with $Cr=85.56$ in Figure 
\figref{pde-euler2d-epi2-Cr85-t15}, and using EXPRB32 with $Cr=85.56$ in Figure \figref{pde-euler2d-exprb32-Cr85-t15}.
}
\figlab{expo-pde-euler2d-vortex}
\end{figure}

\begin{figure}[h!t!b!]
\centering

  \includegraphics[trim=0.0cm 0.0cm 0.0cm 0cm,clip=true,width=0.8\textwidth]{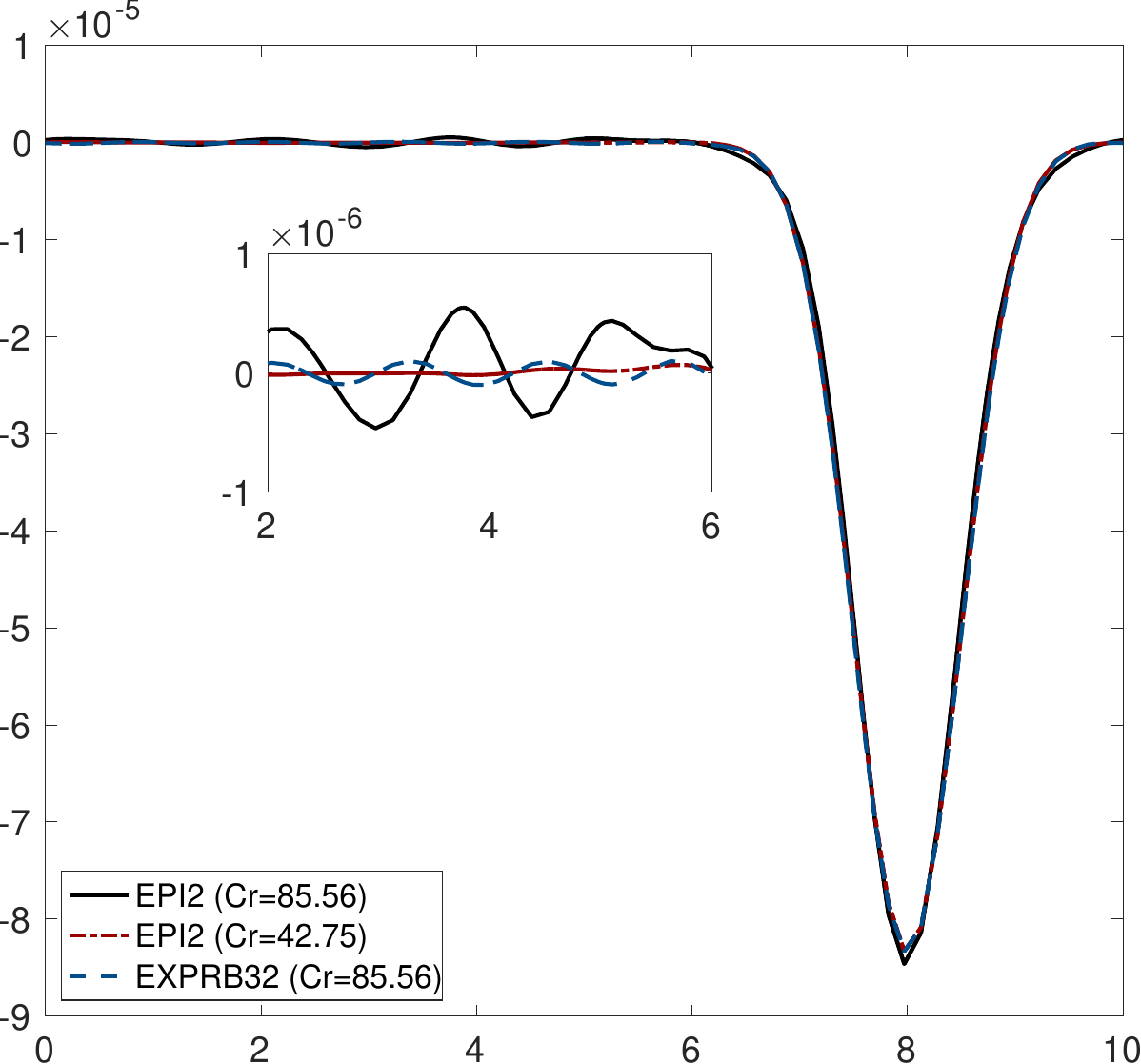}
 
\caption{Translating isentropic vortex example for the two-dimension (black is for EPI2 with $Cr = 85.56$, red-dashed for EPI2 with $Cr = 42.75$, and blue-dashed for EXPRB32 with $Cr=85.56$) Euler equations at $t=15$: a slice along $y=0$ for all sub-figures  in Figure \figref{expo-pde-euler2d-vortex}.}
\figlab{expo-pde-euler2d-vortex-y0}
\end{figure}

\subsubsection{Accuracy and efficient comparison among exponential integrators}

In this section, we take the second, the third, and the fourth-order exponential integrators:
EPI2, EXPRB32, and EXPRB42,\footnote{
EXPRB42  \cite{luan2017fourth} is the fourth-order two-stage method,
\begin{align*}
  q^{(2)} &= q^n + \frac{3}{4} \dt \varphi_1 \LRp{ \frac{3}{4} L \dt } R^n,\\
  q^{n+1} &= q^n + \dt \varphi_1(L \dt) R^n
                 + \frac{32}{9} \dt \varphi_3 \LRp{ L \dt } D_{n,2},
\end{align*}
where $R(q)=Lq + \mc{N}(q)$, $R^n:=R(q^n)$ and $D_{n,2}:= \mc{N}(q^{(2)}) - \mc{N}(q^n)$. 
} respectively, and compare their relative accuracy and efficiency.
Time convergence studies are conducted 
on a uniform mesh and a non-uniform mesh in Figure \figref{expo-pde-euler2d-mesh}.
\begin{figure}[h!t!b!]
\centering
\subfigure[Uniform mesh]{
  \includegraphics[trim=0.0cm 0.0cm 0.0cm 0cm,clip=true,width=0.42\textwidth]{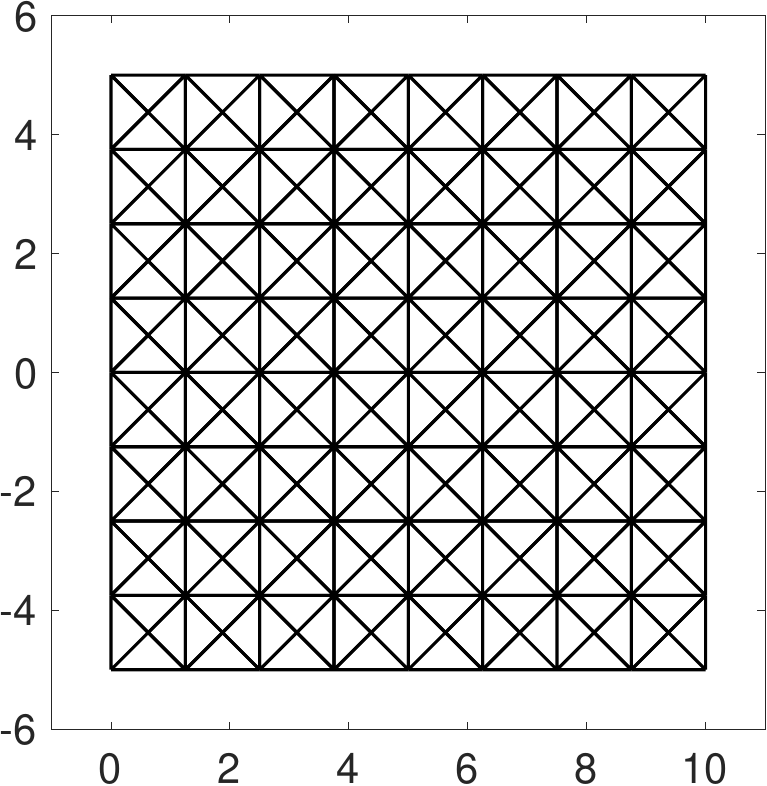}
  \figlab{pde-euler2d-mesh-uniform}
}
\subfigure[Non-uniform mesh]{
  \includegraphics[trim=0.0cm 0.0cm 0.0cm 0cm,clip=true,width=0.42\textwidth]{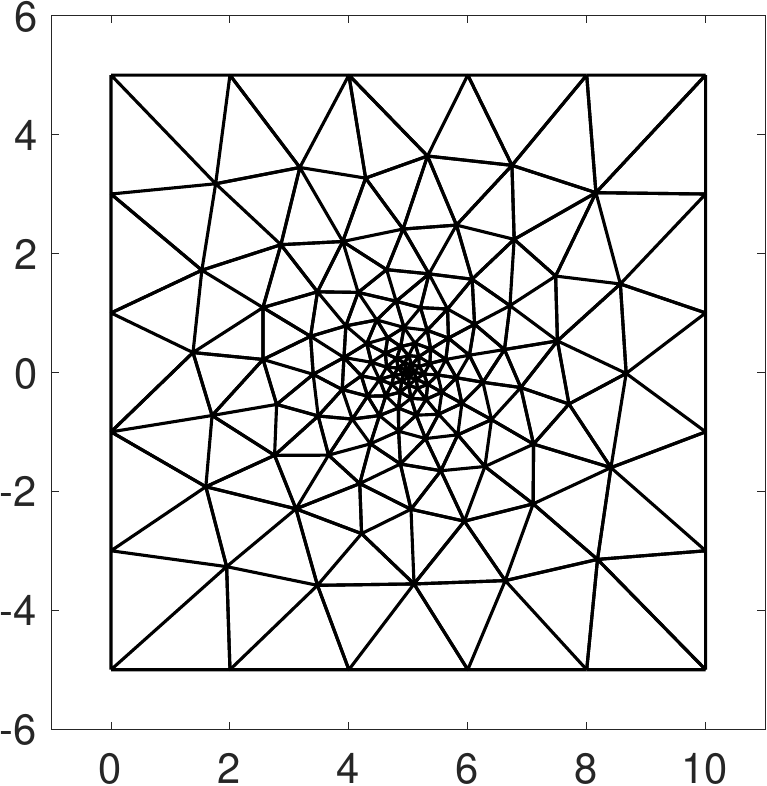} 
  \figlab{pde-euler2d-mesh-nonuniform}
}
\caption{(a) A uniform mesh with $N_e=256$ elements, and (b) a non-uniform meshes with $N_e=250$ elements.}
\figlab{expo-pde-euler2d-mesh}
\end{figure}

We start with a very high-order accurate discretization in space with $k = 16$ for the uniform mesh so that the spatial error (around $10^{-13}$) does not pollute the temporal one. Figure \figref{expo-euler-isentropic-tconv-uni-p16-accuracy} presents the $L^2$-error for the density $\rho$ over a wide range of timestep sizes. 
Note that the numbers on the top of the figure is the corresponding Courant number for the timestep size displayed on the $x$-axis.
As can be seen, EPI2, EXPRB32, and EXPRB42 achieve expected
convergence rate of 2, 3, and 4, respectively. Beyond $10^{-13}$ the
error is dominated by spatial discretization error, which explains why
the error for the last two points (the two smallest timestep size
cases) of the EXPRB42 error curve plateaus.
\begin{figure}[h!t!b!]
    \centering
    \subfigure[Accuracy]{
            \includegraphics[trim=0.0cm 0.cm 0.0cm 0.0cm,clip=true,width=0.43\textwidth]{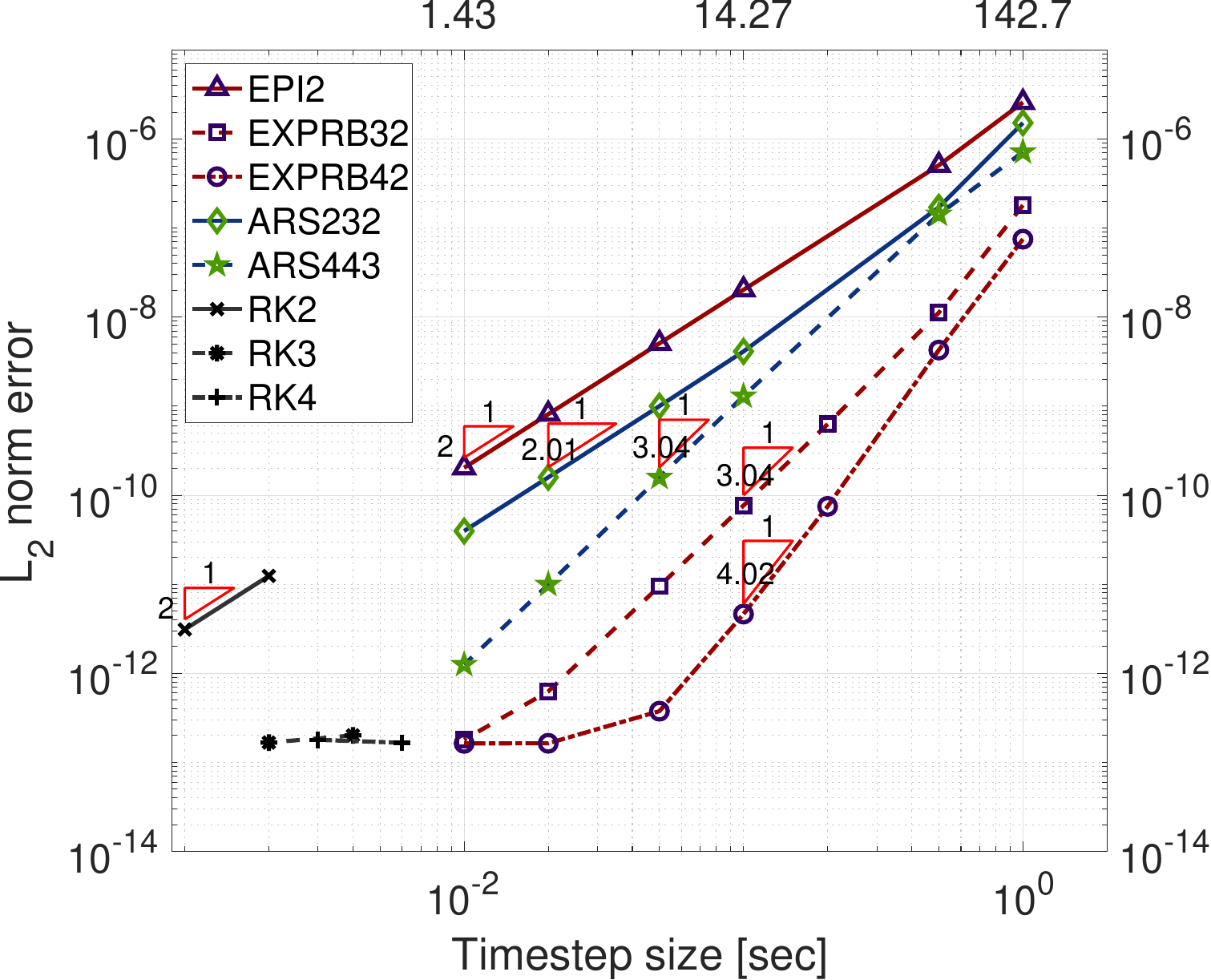}
      \figlab{expo-euler-isentropic-tconv-uni-p16-accuracy}
    }
    \subfigure[Efficiency]{
     \includegraphics[trim=0.0cm 0.cm 0.0cm 0.0cm,clip=true,width=0.413\textwidth]{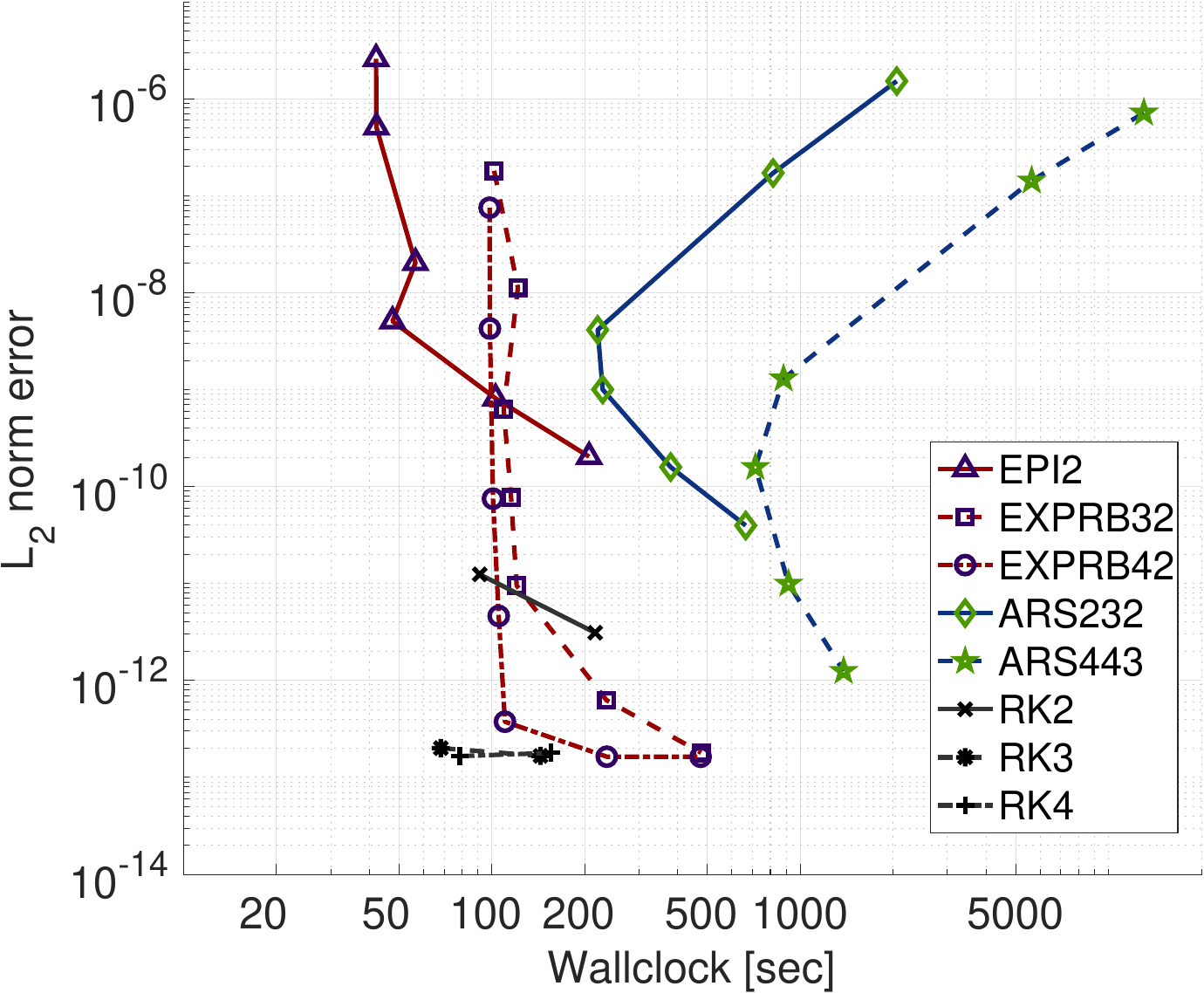}
      \figlab{expo-euler-isentropic-tconv-uni-p16-efficiency}
    }
  \caption{Accuray and effiicency of exponential, IMEX, and RK
    integrators for the isentropic vortex translation in two
    dimensions
    on a uniform mesh  with $N_e = 256$ and $k=16$. }
  \figlab{expo-euler-isentropic-tconv-uni-p16}
\end{figure}

\subsubsection{Accuracy and efficient comparison between exponential and IMEX integrators}

We now compare exponential methods with IMEX (implicit-explicit) time integrators. For IMEX integrators, we integrate the linearized operator $L$ implicitly and the nonlinear operator $\mc{N}$ explicitly.  We consider the second-order ARS232 \cite{ascher1997implicit} and the third-order ARS443 \cite{ascher1997implicit} IMEX schemes.
The $L^2$-error for $\rho$ corresponding to these IMEX methods on the uniform mesh are summarized in 
Figure \figref{expo-euler-isentropic-tconv-uni-p16}. 
As can be seen in Figure \figref{expo-euler-isentropic-tconv-uni-p16-accuracy}, for a given timestep size EXPRB32 is (an order of magnitude) more accurate than ARS443 while EPI2 is (about half order of magnitude) less accurate than ARS232. Efficiency comparison in Figure \figref{expo-euler-isentropic-tconv-uni-p16-efficiency} shows that,  for a given level of accuracy, exponential integrators  EXPRB32 and EPI2 are much  (from two to ten times) more efficient than the IMEX counterparts  ARS232 and ARS443. Though both EXPRB32 and EXPRB42 require two matrix exponential evaluations,
 and hence having similar wallclock, 
  EXPRB42, due to its high-order accuracy, is more efficient than EXPRB32.

\subsubsection{Accuracy and efficient comparison between exponential and RK integrators}

We next compare exponential methods with explicit RK (Runge-Kutta) time integrators. We consider second-order RK2, third-order RK3, and fourth-order RK4 methods. 
Figure \figref{expo-euler-isentropic-tconv-uni-p16-accuracy} shows that 
 RK2 solution converges to the true solution with the second-order accuracy, 
 while RK3 and RK4 solutions immediately saturate at the error level of 
 $\mc{O}(10^{-13})$ as the temporal error is smaller than the spatial
 one. Note that the right most point for each of RK method corresponds
 to (approximately) the largest stable timestep size. Exponential
 methods, again due to their implicit nature, does not have time
 stepsize restriction. As expected Figure
 \figref{expo-euler-isentropic-tconv-uni-p16-efficiency} shows that
 with a same accuracy, exponential integrators are less efficient than
 their same order RK counterparts since the formers require matrix exponential evaluations. 
 Times taken by high-order exponential integrators become comparable to low-order RK counterparts.  {\em This should not understood as a disadvantage. On the contrary, the main advantage of EI is on stiff problems or problem requires large time stepsizes (with CFL number greater than 1) for which explicit RK methods fail. The example shows  that the cost of exponential methods are similar to stable explicit RK methods while stably providing solutions with time stepsizes orders of magnitude larger than the maximal stable time stepsizes for explicit RK methods.}

On the non-uniform mesh, 
we set the center of the vortex to be ${\bf x}_c = (0,0)$ at $t=0$ so that the initial vortex is defined on a coarse region.
This means that more spatial discretization error is introduced than that on the uniform mesh.
In Figure \figref{expo-euler-isentropic-tconv-nonuni-p16}(a),
the saturated error level of $\mc{O}(10^{-11})$ is higher than the counterpart on the uniform mesh. 
All RK solutions immediately reach to the saturated error level of $\mc{O}(10^{-11})$. 
In Figure \figref{expo-euler-isentropic-tconv-nonuni-p16}(b) and 
Table \tabref{expo-euler-isentropic-tconv-nonuni-p16}, 
RK4 is faster than RK2 and RK3. 
Compared to RK4, Exponential DG methods shows slightly better performance at the error level of $\mc{O}(10^{-11})$. 
For example, EPI2 with $Cr=7.49$ is 1.5 times faster than RK4. 


When we lower the solution order from $k=16$ to $k=8$,  
we see the computational gain of Exponential DG methods 
in Table \tabref{expo-euler-isentropic-tconv-nonuni-p8}.
All the numerical solutions saturate at $\mc{O}(10^{-6})$ error level.
The wallclock times of RK2, RK3, and RK4 are $155.8$, $150.0$, and $101.9$. 
EPI2 is three times faster than RK2 and RK3, and two times faster than RK4.
EXPRB32 and EXPRB42 slightly better perform RK4. 


\begin{figure}[h!t!b!]
    \centering
    \subfigure[Accuracy]{
            \includegraphics[trim=0.0cm 0.cm 0.0cm 0.0cm,clip=true,width=0.43\textwidth]{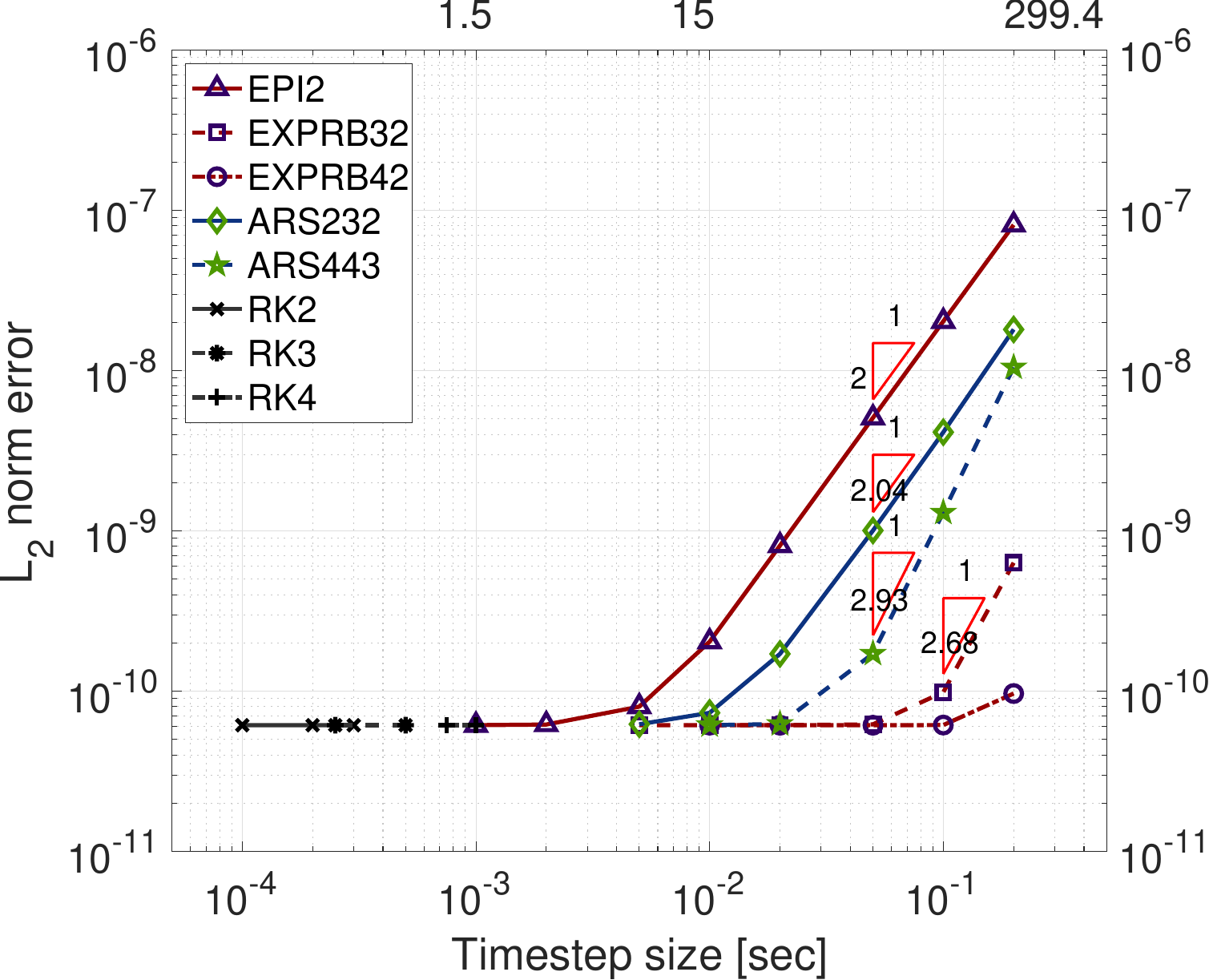}
      \figlab{expo-euler-isentropic-tconv-p16-accuracy}
    }
    \subfigure[Efficiency]{
     \includegraphics[trim=0.0cm 0.cm 0.0cm 0.0cm,clip=true,width=0.405\textwidth]{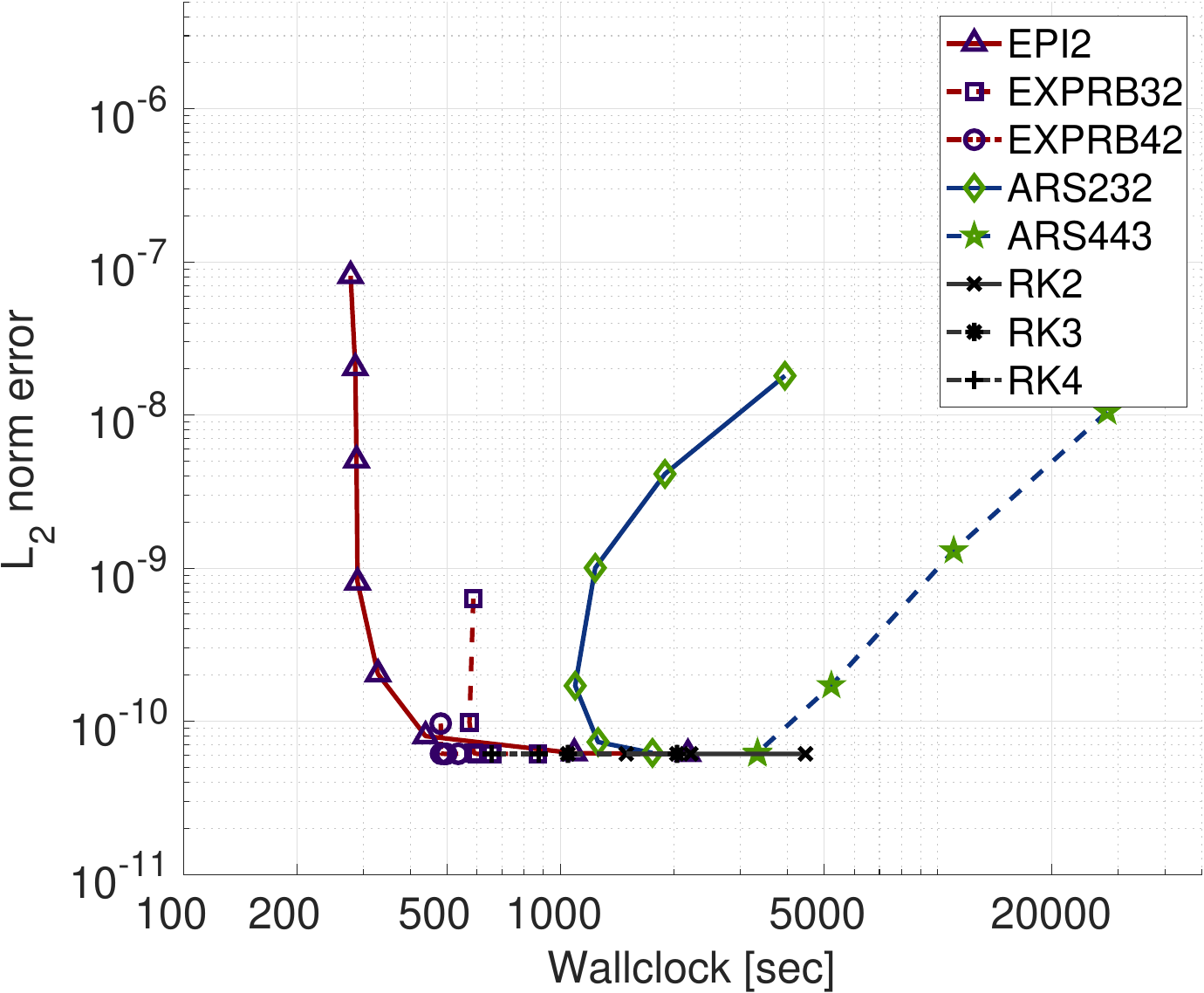}
      \figlab{expo-euler-isentropic-tconv-p16-efficiency}
    }
  \caption{Isentropic vortex translation in two dimensions:
    time convergence study for (a) accuracy and (b) efficiency on non-uniform mesh.
        The computational domain is discretized with $N_e = 250$ and $k=16$. }
  \figlab{expo-euler-isentropic-tconv-nonuni-p16}
\end{figure}


\begin{table}[t]
\caption{Isentropic vortex translation in two dimensions:
    time convergence study on a nonuniform mesh.
    The computational domain is discretized with $N_e = 256$ and $k=16$.
    The error is measured at $T=1$. 
}
\tablab{expo-euler-isentropic-tconv-nonuni-p16}
\begin{center}
\begin{tabular}{*{1}{c}|*{1}{c}|*{2}{c}|*{2}{c}|*{2}{c}|*{1}{c}}
\hline
\multirow{2}{*}{$k=16$} 
  & \multirow{2}{*}{$Cr$}  
  & \multicolumn{2}{c|}{$\rho$}   
  & \multicolumn{2}{c|}{$\rho {\bf u}$}   
  & \multicolumn{2}{c|}{$\rho E$}   
  & \multirow{2}{*}{wc $[s]$}    \tabularnewline
  &     &  error & order &  error & order &  error & order &     \tabularnewline
\hline\hline
\multirow { 3}{*}{   RK2 }&  0.45 &  6.126E-11 &        - & 1.378E-10 &        - & 3.539E-11 &        - &   1491.5 \tabularnewline
                         &  0.30 &  6.127E-11 &   -0.000 & 1.368E-10 &    0.018 & 3.533E-11 &    0.004 &   2213.5 \tabularnewline
                         &  0.15 &  6.127E-11 &    0.000 & 1.366E-10 &    0.002 & 3.531E-11 &    0.001 &   4457.2 \tabularnewline
\multicolumn{9}{c}{} \tabularnewline
\multirow { 2}{*}{   RK3 }&  0.75 & 6.127E-11 &        - & 1.366E-10 &        - & 3.530E-11 &        - &   1044.9 \tabularnewline
                         &  0.37 &  6.127E-11 &    0.000 & 1.366E-10 &    0.000 & 3.530E-11 &    0.000 &   2033.4 \tabularnewline
\multicolumn{9}{c}{} \tabularnewline
\multirow { 2}{*}{   RK4 }&  1.50 & 6.127E-11 &        - & 1.366E-10 &        - & 3.530E-11 &        - &    655.7 \tabularnewline
                         &  1.12 &  6.127E-11 &    0.000 & 1.366E-10 &    0.000 & 3.530E-11 &    0.000 &    873.2 \tabularnewline
\multicolumn{9}{c}{} \tabularnewline
\multirow { 8}{*}{  EPI2 }& 299.42 & 8.115E-08 &        - & 9.834E-08 &        - & 2.858E-07 &        - &    277.3 \tabularnewline
                         & 149.71 &  2.031E-08 &    1.998 & 2.437E-08 &    2.013 & 7.153E-08 &    1.998 &    285.3 \tabularnewline
                         & 74.85 &   5.079E-09 &    2.000 & 6.080E-09 &    2.003 & 1.788E-08 &    2.000 &    287.4 \tabularnewline
                         & 29.94 &   8.093E-10 &    2.004 & 1.068E-09 &    1.898 & 2.840E-09 &    2.008 &    289.3 \tabularnewline
                         & 14.97 &   2.035E-10 &    1.992 & 7.288E-10 &    0.551 & 6.848E-10 &    2.052 &    327.7 \tabularnewline
                         &  7.49 &   8.002E-11 &    1.347 & 1.495E-10 &    2.285 & 1.823E-10 &    1.909 &    438.2 \tabularnewline
                         &  2.99 &   6.188E-11 &    0.281 & 1.369E-10 &    0.096 & 4.527E-11 &    1.520 &   1086.7 \tabularnewline
                         &  1.50 &   6.132E-11 &    0.013 & 1.366E-10 &    0.003 & 3.596E-11 &    0.332 &   2174.3 \tabularnewline
\multicolumn{9}{c}{} \tabularnewline
\multirow { 6}{*}{ EXPRB32 }& 299.42 &  6.324E-10 &        - & 9.584E-10 &        - & 2.218E-09 &        - &    586.5 \tabularnewline
                         & 149.71 &  9.847E-11 &    2.683 & 1.865E-10 &    2.361 & 2.733E-10 &    3.021 &    572.8 \tabularnewline
                         & 74.85 &   6.206E-11 &    0.666 & 2.220E-10 &   -0.251 & 5.044E-11 &    2.438 &    588.4 \tabularnewline
                         & 29.94 &   6.129E-11 &    0.014 & 4.297E-10 &   -0.721 & 4.666E-11 &    0.085 &    599.2 \tabularnewline
                         & 14.97 &   6.133E-11 &   -0.001 & 8.731E-10 &   -1.023 & 7.249E-11 &   -0.636 &    656.4 \tabularnewline
                         &  7.49 &   6.127E-11 &    0.001 & 1.366E-10 &    2.676 & 3.530E-11 &    1.038 &    870.0 \tabularnewline
\multicolumn{9}{c}{} \tabularnewline
\multirow { 5}{*}{ EXPRB42 }& 299.42 & 9.660E-11 &        - & 1.735E-10 &        - & 2.707E-10 &        - &    480.6 \tabularnewline
                         & 149.71 &  6.142E-11 &    0.653 & 1.505E-10 &    0.205 & 3.940E-11 &    2.780 &    481.6 \tabularnewline
                         & 74.85 &   6.127E-11 &    0.004 & 2.662E-10 &   -0.823 & 3.946E-11 &   -0.002 &    496.9 \tabularnewline
                         & 29.94 &   6.128E-11 &   -0.000 & 4.534E-10 &   -0.581 & 4.786E-11 &   -0.211 &    491.3 \tabularnewline
                         & 14.97 &   6.133E-11 &   -0.001 & 9.018E-10 &   -0.992 & 7.430E-11 &   -0.635 &    534.8 \tabularnewline
\multicolumn{9}{c}{} \tabularnewline
\multirow { 6}{*}{ ARS232 }& 299.42 &  1.806E-08 &        - & 1.920E-06 &        - & 1.503E-07 &        - &   3931.8 \tabularnewline
                         & 149.71 &  4.127E-09 &    2.130 & 4.803E-07 &    1.999 & 3.718E-08 &    2.015 &   1889.8 \tabularnewline
                         & 74.85 &   1.005E-09 &    2.038 & 1.201E-07 &    2.000 & 9.266E-09 &    2.005 &   1234.8 \tabularnewline
                         & 29.94 &   1.707E-10 &    1.935 & 1.922E-08 &    2.000 & 1.481E-09 &    2.001 &   1092.3 \tabularnewline
                         & 14.97 &   7.321E-11 &    1.221 & 4.806E-09 &    2.000 & 3.707E-10 &    1.998 &   1256.3 \tabularnewline
                         &  7.49 &   6.214E-11 &    0.237 & 1.209E-09 &    1.991 & 9.790E-11 &    1.921 &   1751.2 \tabularnewline
\multicolumn{9}{c}{} \tabularnewline
\multirow { 5}{*}{ ARS443 }& 299.42 & 1.049E-08 &        - & 1.183E-07 &        - & 3.750E-08 &        - &  28140.8 \tabularnewline
                         & 149.71 & 1.302E-09 &    3.010 & 1.484E-08 &    2.995 & 4.702E-09 &    2.996 &  11017.5 \tabularnewline
                         & 74.85 &  1.714E-10 &    2.925 & 1.863E-09 &    2.994 & 5.787E-10 &    3.022 &   5219.5 \tabularnewline
                         & 29.94 &  6.235E-11 &    1.104 & 1.812E-10 &    2.543 & 5.088E-11 &    2.653 &   3327.1 \tabularnewline
                         & 14.97 &  6.131E-11 &    0.024 & 1.374E-10 &    0.399 & 3.562E-11 &    0.514 &   3319.0 \tabularnewline
\multicolumn{9}{c}{} \tabularnewline
\hline\hline
\end{tabular}
\end{center}
\end{table}

\begin{table}[t]
\caption{Isentropic vortex translation in two dimensions:
    time convergence study on a nonuniform mesh.
    The computational domain is discretized with $N_e = 256$ and $k=8$.
    The error is measured at $t=1$. 
}
\tablab{expo-euler-isentropic-tconv-nonuni-p8}

\begin{center}
\begin{tabular}{*{1}{c}|*{1}{c}|*{2}{c}|*{2}{c}|*{2}{c}|*{1}{c}}
\hline
\multirow{2}{*}{$k=8$} 
  & \multirow{2}{*}{$Cr$}  
  & \multicolumn{2}{c|}{$\rho$}   
  & \multicolumn{2}{c|}{$\rho {\bf u}$}   
  & \multicolumn{2}{c|}{$\rho E$}   
  & \multirow{2}{*}{wc $[s]$}    \tabularnewline
  &     &  error & order &  error & order &  error & order &     \tabularnewline
\hline\hline
\multirow { 2}{*}{  RK2 }&  0.63 & 1.533E-06 &        - & 2.884E-06 &        - & 7.851E-07 &        - &    155.8 \tabularnewline
                         &  0.33 & 1.533E-06 &    0.000 & 2.884E-06 &    0.000 & 7.851E-07 &    0.000 &    306.7 \tabularnewline
\multicolumn{9}{c}{} \tabularnewline
\multirow { 2}{*}{   RK3 }&  0.81 & 1.533E-06 &        - & 2.884E-06 &        - & 7.851E-07 &        - &    150.0 \tabularnewline
                         &  0.41 &  1.533E-06 &    0.000 & 2.884E-06 &    0.000 & 7.851E-07 &    0.000 &    294.3 \tabularnewline
\multicolumn{9}{c}{} \tabularnewline
\multirow { 3}{*}{   RK4 }&  1.36 & 1.533E-06 &        - & 2.884E-06 &        - & 7.851E-07 &        - &    101.9 \tabularnewline
                         &  0.90 &  1.533E-06 &    0.000 & 2.884E-06 &    0.000 & 7.851E-07 &    0.000 &    152.5 \tabularnewline
                         &  0.45 &  1.533E-06 &    0.000 & 2.884E-06 &    0.000 & 7.851E-07 &    0.000 &    307.7 \tabularnewline
\multicolumn{9}{c}{} \tabularnewline
\multirow { 5}{*}{  EPI2 }& 90.42 & 1.536E-06 &        - & 2.887E-06 &        - & 8.365E-07 &        - &     45.2 \tabularnewline
                         & 45.21 &  1.533E-06 &    0.003 & 2.885E-06 &    0.001 & 7.907E-07 &    0.081 &     49.2 \tabularnewline
                         & 22.60 &  1.534E-06 &   -0.001 & 2.893E-06 &   -0.004 & 7.938E-07 &   -0.006 &     52.6 \tabularnewline
                         &  9.04 &  1.533E-06 &    0.001 & 2.884E-06 &    0.003 & 7.849E-07 &    0.012 &     59.9 \tabularnewline
                         &  4.52 &  1.533E-06 &    0.000 & 2.884E-06 &    0.000 & 7.851E-07 &   -0.000 &    113.3 \tabularnewline
\multicolumn{9}{c}{} \tabularnewline
\multirow { 5}{*}{ EXPRB32 }& 90.42 & 1.533E-06 &        - & 2.884E-06 &        - & 7.862E-07 &        - &     91.4 \tabularnewline
                         & 45.21 &  1.533E-06 &    0.000 & 2.884E-06 &    0.000 & 7.862E-07 &    0.000 &     94.9 \tabularnewline
                         & 22.60 &  1.533E-06 &    0.000 & 2.884E-06 &    0.000 & 7.866E-07 &   -0.001 &     99.5 \tabularnewline
                         &  9.04 &  1.533E-06 &    0.000 & 2.884E-06 &    0.000 & 7.849E-07 &    0.002 &    111.0 \tabularnewline
                         &  4.52 &  1.533E-06 &    0.000 & 2.884E-06 &    0.000 & 7.851E-07 &   -0.000 &    221.3 \tabularnewline
\multicolumn{9}{c}{} \tabularnewline
\multirow { 5}{*}{ EXPRB42 }& 90.42 & 1.533E-06 &        - & 2.884E-06 &        - & 7.861E-07 &        - &     77.9 \tabularnewline
                         & 45.21 &    1.533E-06 &    0.000 & 2.884E-06 &    0.000 & 7.859E-07 &    0.000 &     75.7 \tabularnewline
                         & 22.60 &    1.533E-06 &    0.000 & 2.884E-06 &    0.000 & 7.870E-07 &   -0.002 &     81.4 \tabularnewline
                         &  9.04 &    1.533E-06 &    0.000 & 2.884E-06 &    0.000 & 7.849E-07 &    0.003 &    110.9 \tabularnewline
                         &  4.52 &    1.533E-06 &    0.000 & 2.884E-06 &    0.000 & 7.851E-07 &   -0.000 &    221.1 \tabularnewline
\multicolumn{9}{c}{} \tabularnewline
\multirow { 5}{*}{ ARS232 }& 90.42 & 1.515E-06 &        - & 3.469E-06 &        - & 7.177E-07 &        - &    343.7 \tabularnewline
                         & 45.21 &    1.530E-06 &   -0.014 & 2.925E-06 &    0.246 & 7.612E-07 &   -0.085 &    248.6 \tabularnewline
                         & 22.60 &    1.533E-06 &   -0.003 & 2.887E-06 &    0.019 & 7.803E-07 &   -0.036 &    205.9 \tabularnewline
                         &  9.04 &    1.533E-06 &    0.000 & 2.884E-06 &    0.001 & 7.852E-07 &   -0.007 &    226.1 \tabularnewline
                         &  4.52 &    1.533E-06 &    0.000 & 2.884E-06 &    0.000 & 7.851E-07 &    0.000 &    317.7 \tabularnewline
\multicolumn{9}{c}{} \tabularnewline
\multirow { 5}{*}{ ARS443 }& 90.42 & 1.512E-06 &    - & 2.882E-06 &        - & 6.233E-07 &        - &   1616.1 \tabularnewline
                         & 45.21 &  1.526E-06 &   -0.013 & 2.881E-06 &    0.001 & 6.762E-07 &   -0.118 &   1050.9 \tabularnewline
                         & 22.60 &  1.531E-06 &   -0.005 & 2.883E-06 &   -0.001 & 7.490E-07 &   -0.148 &    765.3 \tabularnewline
                         &  9.04 &  1.533E-06 &   -0.001 & 2.884E-06 &   -0.000 & 7.816E-07 &   -0.046 &    684.7 \tabularnewline
                         &  4.52 &  1.533E-06 &    0.000 & 2.884E-06 &    0.000 & 7.846E-07 &   -0.006 &    797.5 \tabularnewline
\multicolumn{9}{c}{} \tabularnewline
\hline\hline
\end{tabular}
\end{center}
\end{table}

To demonstrate the high-order convergence in space,  
we perform the spatial convergence test.
We use a sequence of nested meshes with $N_e = \LRp{256,1024,4096,16384}$ 
for $k=\LRp{1,2,3,4}$ and measure the errors at $t=1$. 
As can be seen in Table \tabref{expo-euler-isentropic-sconv-uni},
the convergence rate of $(k + \half)$ is observed as refining the meshes.

\begin{table}[t]
\caption{Spatial convergence for the isentropic vortex translation in two dimensions.
    We use a sequence of nested meshes with $N_e = \LRc{256, 1024, 4096, 16384}$ for $k=\LRc{1,2,3,4}$. 
    We use EXPRB42 scheme with $\dt=0.01$ and measure the errors at $t=1$.  
}
\tablab{expo-euler-isentropic-sconv-uni}
\begin{center}
\begin{tabular}{*{1}{c}|*{1}{c}|*{2}{c}|*{2}{c}|*{2}{c}}
\hline
\multirow{2}{*}{ } 
  & \multirow{2}{*}{$h$}  
  & \multicolumn{2}{c|}{$\rho$}   
  & \multicolumn{2}{c|}{$\rho {\bf u}$}   
  & \multicolumn{2}{c|}{$\rho E$}   
  \tabularnewline
  &    &  error & order &  error & order &  error & order       \tabularnewline
\hline\hline
\multirow { 4}{*}{    $k=1$ }&  1.00 & 9.942E-04 &        - & 7.389E-03 &        - & 1.595E-03 &        - \tabularnewline
                         &  0.50 & 7.662E-04 &    0.376 & 1.999E-03 &    1.886 & 5.181E-04 &    1.622 \tabularnewline
                         &  0.25 & 4.344E-04 &    0.819 & 5.312E-04 &    1.912 & 2.049E-04 &    1.338 \tabularnewline
                         &  0.12 & 2.312E-04 &    0.910 & 1.383E-04 &    1.941 & 8.571E-05 &    1.257 \tabularnewline
\multicolumn{8}{c}{} \tabularnewline
\multirow { 4}{*}{    $k=2$ }&  1.00 & 2.909E-04 &        - & 7.382E-04 &        - & 4.423E-04 &        - \tabularnewline
                         &  0.50 & 1.201E-04 &    1.276 & 1.418E-04 &    2.380 & 4.965E-05 &    3.155 \tabularnewline
                         &  0.25 & 3.235E-05 &    1.892 & 1.570E-05 &    3.175 & 7.084E-06 &    2.809 \tabularnewline
                         &  0.12 & 5.787E-06 &    2.483 & 2.084E-06 &    2.913 & 1.257E-06 &    2.495 \tabularnewline
\multicolumn{8}{c}{} \tabularnewline
\multirow { 4}{*}{    $k=3$ }&  1.00 & 1.286E-04 &        - & 2.343E-04 &        - & 6.979E-05 &        - \tabularnewline
                         &  0.50 & 1.972E-05 &    2.705 & 1.541E-05 &    3.926 & 5.605E-06 &    3.638 \tabularnewline
                         &  0.25 & 2.123E-06 &    3.215 & 9.655E-07 &    3.996 & 4.450E-07 &    3.655 \tabularnewline
                         &  0.12 & 2.112E-07 &    3.329 & 6.849E-08 &    3.817 & 4.386E-08 &    3.343 \tabularnewline
\multicolumn{8}{c}{} \tabularnewline
\multirow { 4}{*}{    $k=4$ }&  1.00 & 2.686E-05 &        - & 3.180E-05 &        - & 1.311E-05 &        - \tabularnewline
                         &  0.50 & 2.067E-06 &    3.700 & 1.283E-06 &    4.631 & 4.381E-07 &    4.903 \tabularnewline
                         &  0.25 & 9.543E-08 &    4.437 & 4.294E-08 &    4.901 & 1.333E-08 &    5.039 \tabularnewline
                         &  0.12 & 4.540E-09 &    4.394 & 1.458E-09 &    4.880 & 5.294E-10 &    4.654 \tabularnewline
\multicolumn{8}{c}{} \tabularnewline
\hline\hline
\end{tabular}
\end{center}
\end{table}

\subsubsection{Performance of Exponential DG on parallel computers}

Now we study the parallel performance, namely weak and strong scalings, of Exponential DG methods for
three-dimensional Euler equations.  For this purpose, we choose the EPI2 integrator. Parallel simulations are conducted on
Stampede2 at the Texas Advanced Computing Center (TACC) 
 using Skylake (SKX)
nodes. 
Each node of SKX consists of 48 cores of Intel Xeon Platinum
8160 2.1GHz processors and 192GB DDR4 RAM. The interconnect is a 100GB/s Intel Omni-Path (OPA) network with a fat-tree topology. 

We begin with strong scaling in which the problem size is fixed while the number of cores increases.
Table \tabref{expo-euler-strongscale-isentropic-varyCr} compares the efficiencies of two different timestep sizes: $\dt=0.125$ and $\dt=0.25$
on the mesh with $N_e=51200$ (elements) and $k=6$ (solution order). For either of the timestep sizes, the corresponding run with 32 cores and $1600$ elements per core is served as the based line. 
As the number of cores increases (i.e. the number of elements per core decreases) communication-computation overlapping is less effective and thus decreasing the efficiency. 
The efficiency with $\dt=0.125$ is slightly higher than that with $\dt=0.25$ as the latter requires more Krylov iterations than the former:
the total number of Krylov iterations 
$N_{Krylov}$ is 
$561$ for $Cr=15.49 (\dt=0.125)$ and $1024$ for $Cr=32.98 (\dt=0.25)$
\footnote{
We observe $N_{Krylov}=556$ for $Cr=7.75(\dt=0.0625)$. 
The total number of Krylov iterations is proportional to Courant number above a certain Courant number. 
In Table \tabref{expo-euler-weakscale-isentropic-fixdt}, 
$N_{Krylov}$ is about 40 for $Cr<=2.8$. 
However, for $Cr>2.8$, doubling Courant number tends to double $N_{Krylov}$. 
}.
This implies that the spectrum of the linear operator becomes broad by increasing the timestep size\footnote{
Note that the input argument of $\varphi$-function is $\triangle t L$.
}.


\begin{table}[t]
\caption{Isentropic vortex translation in three dimensions:
  strong scaling results for EPI2 with $Cr=16.49  $ ($\dt=0.125$)
  and $Cr=32.98$ ($\dt=0.25$)
  are performed with $N_e=51200$ and $k=6$ up to $t=1$.
}
\tablab{expo-euler-strongscale-isentropic-varyCr}
\begin{center}
\begin{tabular}{c|c||cc|cc}
\hline
\multirow{2}{*}{$\# cores$} & \multirow{2}{*}{$N_e/core$}
& \multicolumn{2}{c|}{$\dt=0.125$}
& \multicolumn{2}{c}{$\dt=0.25$}   \tabularnewline
     &     & Wallclock [s] & Efficiency
            & Wallclock [s] & Efficiency \tabularnewline
\hline\hline
32   & 1600 & 2346  & 100  & 4151  & 100 \tabularnewline
64   & 800  & 1199  & 97.8 & 2121  & 97.9 \tabularnewline
128  & 400  & 607.7 & 96.5 & 1255  & 82.7 \tabularnewline
256  & 200  & 306.9 & 95.6 & 605.7 & 85.7 \tabularnewline
512  & 100  & 237.3 & 61.8 & 419.1 & 61.9 \tabularnewline
1024 & 50   & 122.7 & 59.7 & 219.1 & 59.2 \tabularnewline
2048 & 25   & 61.33 & 59.5 & 118.0 & 55.0 \tabularnewline
4096 & 12.5 & 33.72 & 54.4 & 63.12 & 51.4 \tabularnewline
\hline
\end{tabular}
\end{center}
\end{table}

Next, we conduct a strong scaling test using the exponential DG with EPI2,
 $Cr=8.1$ ($\dt=0.03125$), $N_e = 3,276,800$ and $k=4$.
We choose the number of processors to be $n_p=\LRc{16,32,64,128,256,512,868}\times 48$
so that the number of elements per core approximately becomes $\LRc{4267, 2133, 1067, 533, 267, 133, 79}$, 
i.e., 
every time we double the number of processors, the number of elements is halved.
The speedup factors\footnote{
Speedup is defined as $\frac{T_s}{T_p}$
with $T_s$ serial wallclock time and $T_p$ parallel wallclock time.
} for all cases in Figure \figref{expo-euler-strongscale-isentropic-withEPI2} show that the exponential DG approach delivers good strong scalability up to 41664 cores\textemdash the maximum number of cores in Skylake system in TACC.


\begin{figure}[h!t!b!]
    \centering
      \includegraphics[trim=0.0cm 0.cm 0.0cm 0.0cm,clip=true,width=0.83\textwidth]{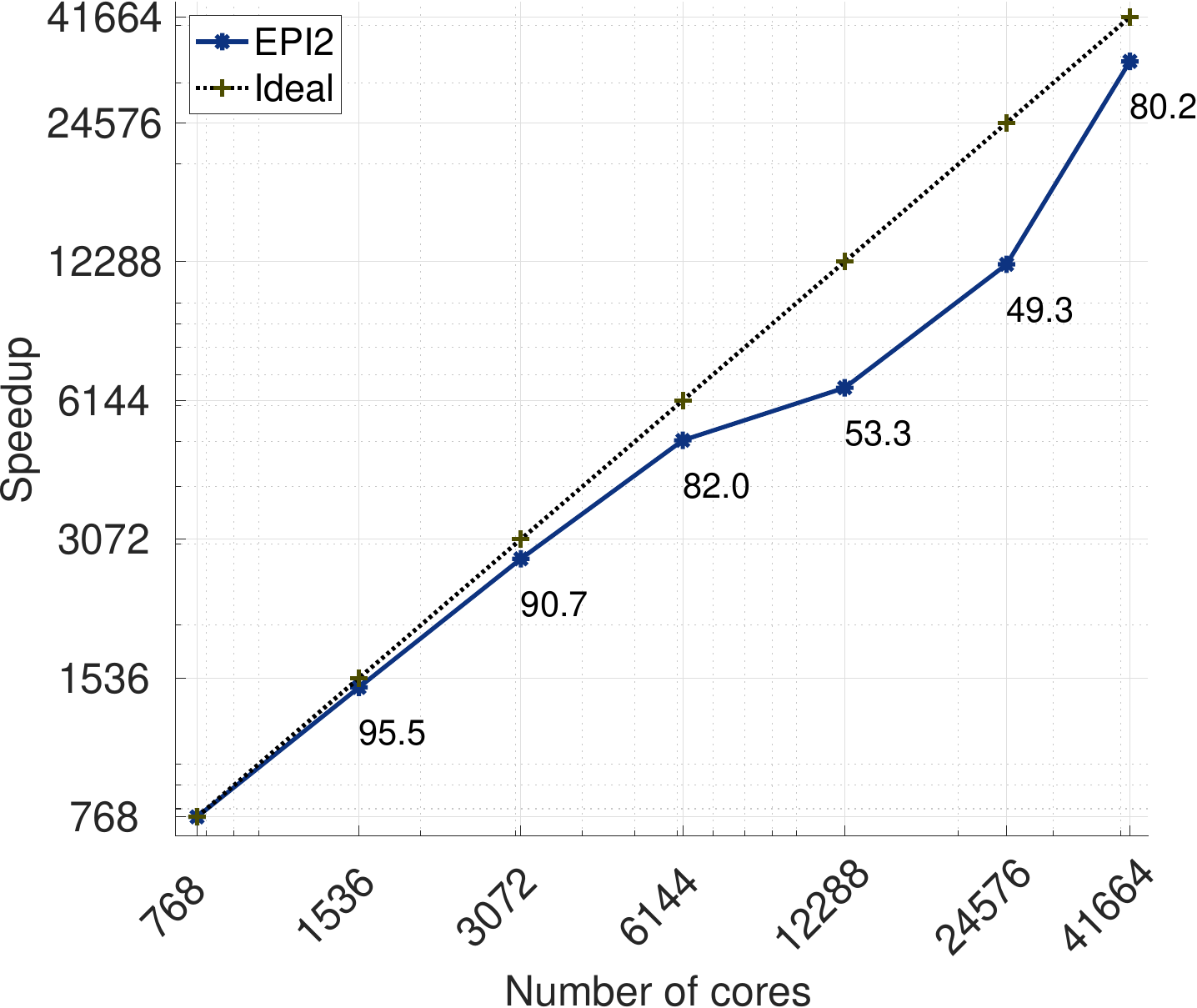}
      \figlab{expo-euler-strongscale-isentropic-epi2}
  \caption{Isentropic vortex translation in three dimensions: 
    strong scaling study for EPI2 with $Cr=8.1$ ($\dt=0.03125$). 
    The computational domain is discretized using $N_e = 3,276,800$ and $k=4$. The number of processors used are $n_p=\LRc{16,32,64,128,256,512,868}\times 48$.}
  \figlab{expo-euler-strongscale-isentropic-withEPI2}
\end{figure}

For the weak scaling test, 
we assign the same amount of work to each processor (by refining the mesh)
while increasing the number of processors. 
Our exponential DG approach uses EPI2, 
$N_e = \LRc{1, 8, 64, 512, 4096}\times 100$, 
and $k=8$.  
The number of processors is chosen in the set $\LRc{4, 32, 256, 2048, 16384}$ so that  the number of elements per core is $25$ (i.e. $18,225$ degrees-of-freedom).
We have tabulated the weak scaling results in Table \tabref{expo-euler-weakscale-isentropic-fixcr}, in which
 each row block shows, for a fixed Courant number, the number of processors, the timestep sizes, the final times, the wallclock times taken, 
 and the number of Krylov iterations.  
 For each fixed Courant number, good weak scalings can be seen through the wallclock times (and $N_{Krylov}$) that do not vary much as the number of processors (and thus the problem size) increases. 
To see this visually, we plot the average time-per-timestep against the number of degrees-of-freedom in Figure \figref{expo-euler-weakscale-isentropic-fixcr}: almost plateau curve for each Courant number indicates favorable weak scaling can be obtained by  the Exponential DG method. As can also be observed, the number of Krylov iterations $N_{Krylov}$, and hence the wallclock time, 
scales linearly with the Courant number. 


\begin{table}[h!t!b!]
\caption{Weak scaling study with fixed Courant numbers for the isentropic vortex translation with three dimensional Euler equations.
    The exponential DG approach consists of EPI2, $N_e = \LRc{1, 8, 64, 512,4096}\times 100$, and $k=8$.
Each curve presents the wallclock time per timestep for various number of degrees-of-freedom with a fixed Courant number. The numbers on each of the curves are the number of cores. 
}
\tablab{expo-euler-weakscale-isentropic-fixcr}
\begin{center}
\begin{tabular}{cc|c|c|c|c}
\hline
$Cr$ &$N_e$ & $\dt$ & Final time & Wallclock [s] & $N_{Krylov}$ \tabularnewline
\hline\hline 
\multirow{4}{*}{$ 5.59$}  
  &  4   & 0.2    & 0.8  & 15.3  & 103 \tabularnewline
  & 32   & 0.1    & 0.4  & 17.8  & 90 \tabularnewline
  & 256  & 0.05   & 0.2  & 24.4  & 81 \tabularnewline
  & 2048 & 0.025  & 0.1  & 27.8  & 81 \tabularnewline
  & 16384& 0.0125 & 0.05 & 29.3  & 82 \tabularnewline
\tabularnewline
\multirow{4}{*}{$ 22.4$}  
  &  4   & 0.2   & 0.8  & 50.4  & 397 \tabularnewline
  & 32   & 0.1   & 0.4  & 89.9  & 381 \tabularnewline
  & 256  & 0.05  & 0.2  & 85.8  & 355 \tabularnewline
  & 2048 & 0.025 & 0.1  & 86.9  & 348 \tabularnewline
  &16384 & 0.0125& 0.05 & 93.1  & 358 \tabularnewline
\tabularnewline
\multirow{4}{*}{$ 55.9$}  
  &  4   & 0.2   & 0.8  & 130.9 & 1024\tabularnewline
  & 32   & 0.1   & 0.4  & 240.4 & 1023 \tabularnewline
  & 256  & 0.05  & 0.2  & 246.2 & 1024 \tabularnewline
  & 2048 & 0.025 & 0.1  & 248.8 & 1024 \tabularnewline
  &16384 & 0.0125& 0.05 & 342.2 & 1408 \tabularnewline
\hline
\end{tabular}
\end{center}
\end{table}

\begin{figure}[h!t!b!]
  \centering
  \includegraphics[trim=0.0cm 0.0cm 0.0cm 0.0cm,clip=true,width=0.8\textwidth]{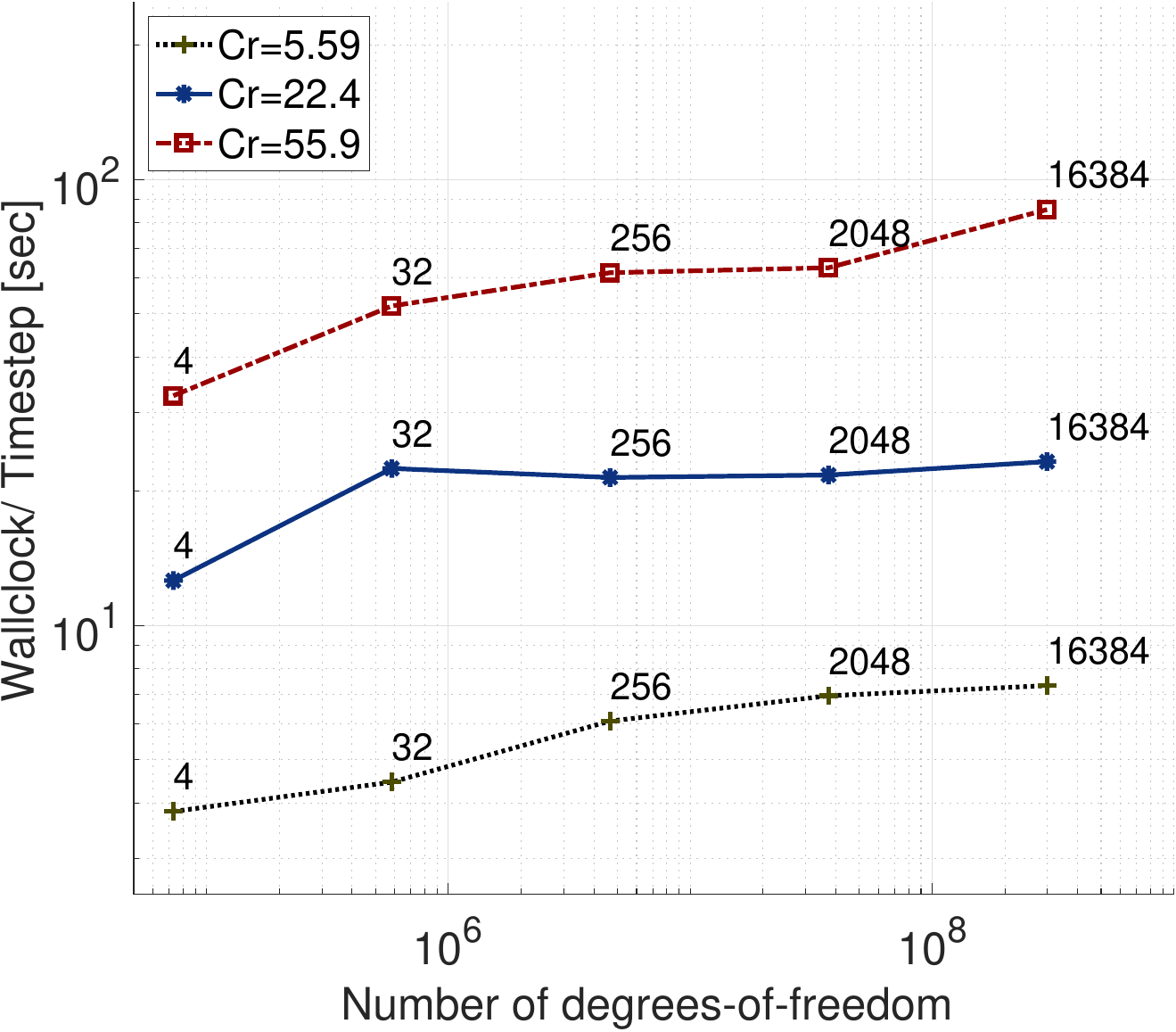} 
  \caption{Weak scaling study with fixed Courant numbers for the isentropic vortex translation with three dimensional Euler equations.
    The exponential DG approach consists of EPI2, $N_e = \LRc{1, 8, 64, 512,4096}\times 100$, and $k=8$. 
    The number of processors for each $N_e$ is $n_p=\LRc{4,32,256,2048,16384}$
    so that each processor has 25 elements ($18,225$ degrees-of-freedom) for all cases. }
  \figlab{expo-euler-weakscale-isentropic-fixcr}
\end{figure}

How does the weak scaling behaves if we fix timestep size $\dt$ instead of Courant number $Cr$? 
In this case, refining the mesh (in order to keep the number of elements per core the same) adds geometrically-induced stiffness to the system,
and thus making the total number of Krylov iterations to increase. 
This is verified in Table \tabref{expo-euler-weakscale-isentropic-fixdt}, which shows linear growth in the total number of Krylov iterations 
as the mesh is refined for $Cr>2.8$. As shown in Figure \figref{expo-euler-weakscale-isentropic-fixdt}, the increase in number of Krylov iterations induces the growth in wallclock time.


\begin{table}[h!t!b!]
\caption{Weak scaling study with fixed $dt$ for the isentropic vortex translation with three dimensional Euler equations.
The exponential DG approach uses EPI2,
$N_e = \LRc{100, 800, 6400, 51200}$, and $k=8$.
The corresponding number of processors are $n_p=\LRc{4,32,256,2048}$, respectively, 
    so that each processor has 25 elements for each case.
}
\tablab{expo-euler-weakscale-isentropic-fixdt}
\begin{center}
\begin{tabular}{c||cc|cc|cc|cc}
\hline
\multirow{2}{*}{$\dt$} 
& \multicolumn{2}{c|}{$N_e=100$}  
& \multicolumn{2}{c|}{$N_e=800$} 
& \multicolumn{2}{c|}{$N_e=6400$} 
& \multicolumn{2}{c}{$N_e=51200$}   \tabularnewline
          &  $Cr$ &  $N_{Krylov}$ 
          &  $Cr$ &  $N_{Krylov}$ 
          &  $Cr$ &  $N_{Krylov}$ 
          &  $Cr$ &  $N_{Krylov}$\tabularnewline
 \hline\hline
 0.025 & 0.7   & 40     & 1.4    & 40    & 2.8       & 45    & 5.6     & 81  \tabularnewline
 0.25  &  7.0  & 128   & 14.0 & 233   & 27.9.    &448   & 55.9   & 1024 \tabularnewline
 2.5    & 69.8 & 1536 &139.7 & 3072 &279.32 &5120 & 558.7 & 10112 \tabularnewline
\hline
\end{tabular}
\end{center}
\end{table}

\begin{figure}[h!t!b!]
  \centering
  \includegraphics[trim=0.0cm 0.0cm 0.0cm 0.0cm,clip=true,width=0.8\textwidth]{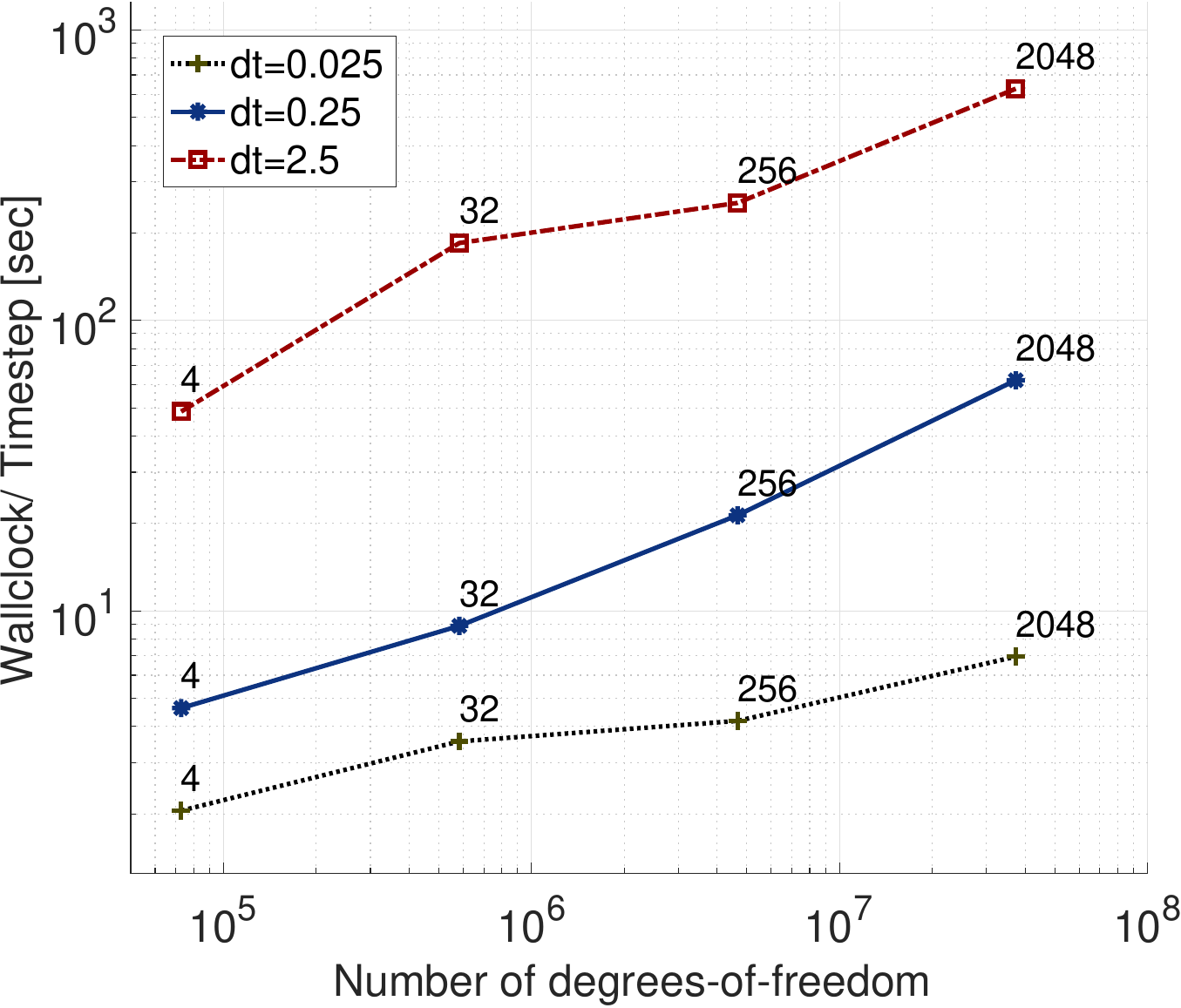} 
  \caption{Weak scaling study with fixed $\dt$ for the isentropic vortex translation with three dimensional Euler equations.
The exponential DG approach uses EPI2,
$N_e = \LRc{100, 800, 6400, 51200}$, and $k=8$.
The corresponding number of processors are $n_p=\LRc{4,32,256,2048}$, respectively, 
    so that each processor has 25 elements for each case.}
  \figlab{expo-euler-weakscale-isentropic-fixdt}
\end{figure}

\subsection{Euler equations: Kelvin-Helmholtz instability}

Kelvin--Helmholtz instability (KHI) is an important mechanism in the development of turbulence.
KHI occurs when two fluids meet across their interface 
with different densities and tangential velocities.
As time goes by, 
small disturbances at the interface grow exponentially,
 and the interface rolls up into KH rotors 
\cite{drazin2004hydrodynamic, springel2010pur, lecoanet2016validated}.
The computational domain is $\Omega =(-5,5) \times (0,5)$. 
We apply periodic boundary condition to the lateral direction, 
whereas no-slip boundary condition to the top and the bottom walls. 

The initial conditions are chosen as 
\begin{align*}
  \rho &= 1 + \half \LRp{ \tanh\LRp{\frac{y - s_1}{a}} - \tanh\LRp{\frac{y - s_2}{a}} },\\
     u &= 0.1 + \LRp{ \tanh\LRp{\frac{y-s_1}{a}} - \tanh\LRp{\frac{y-s_2}{a}} - 1 },\\
     v &= A \sin\LRp{2\pi x} \LRp{ \exp\LRp{-\frac{(y-s_1)^2}{\sigma^2} } 
                                 + \exp\LRp{-\frac{(y-s_2)^2}{\sigma^2} }  }, \\
  \pres &= \gamma^{-1},
\end{align*}
where we take $a=0.05$, $A=0.01$, $\sigma=0.2$,
 $s_1=2$, $s_2=3$,
and  
$c_{\max}=0.1$ and $c_E=1$ for entropy viscosity. 

The numerical simulations are performed with EPI2 and RK4 methods over 
  the uniform mesh with $k=6$ and $N_e=722$ for $t\in\LRs{0,100}$.
  We take $\dt=0.005=:\dt_{RK4}$ for RK4, and 
  $\dt=0.04(=\dt_{RK4}\times 8)$ and $\dt=0.5(=\dt_{RK4}\times 100)$ for EPI2.
The temperature fields are plotted at $t=50$ and $t=100$ 
in Figure \figref{euler-khi-snapshots}. 
The wallock times of RK4, EPI2 (with $\dt=0.04$) and EPI2 (with $\dt=0.5$) are 
$8729s$, $3505s$ and $2656s$, respectively. EPI2 with $\dt=0.04$ (8 times larger time stepsize)
are in good agreement with RK4 but about $2.5$ times faster.
However, EPI2 solution with $\dt=0.5$ is quite deviated from RK4 solution. 
This 
is due to
the way to approximate the entropy residual. 
\footnote{
  Using high-order exponential integrators does not improve the solution quality 
  unlike the isentropic vortex example in Figure 
  \figref{expo-pde-euler2d-vortex}. 
}
Note that in this study, the temporal tendency of the entropy residual 
is approximated by the first-order Euler method, 
i.e., $\dd{\Scal}{t} \approx \frac{\Scal^n - \Scal^{n-1}}{\dt}$.
Thus, using different timestep sizes yields different artificial viscosity. 
The accuracy of the residual computation can be enhanced
by incorporating the second-order approximation such as backward differentiation formula, 
but this is out of the scope of the paper. 

\begin{figure}[h!t!b!]
  \centering
  \subfigure[EPI2 with $\dt_{RK4} \times 100$]{
    \includegraphics[trim=0.0cm 0.cm 0.0cm 0.0cm,clip=true,width=0.43\textwidth]{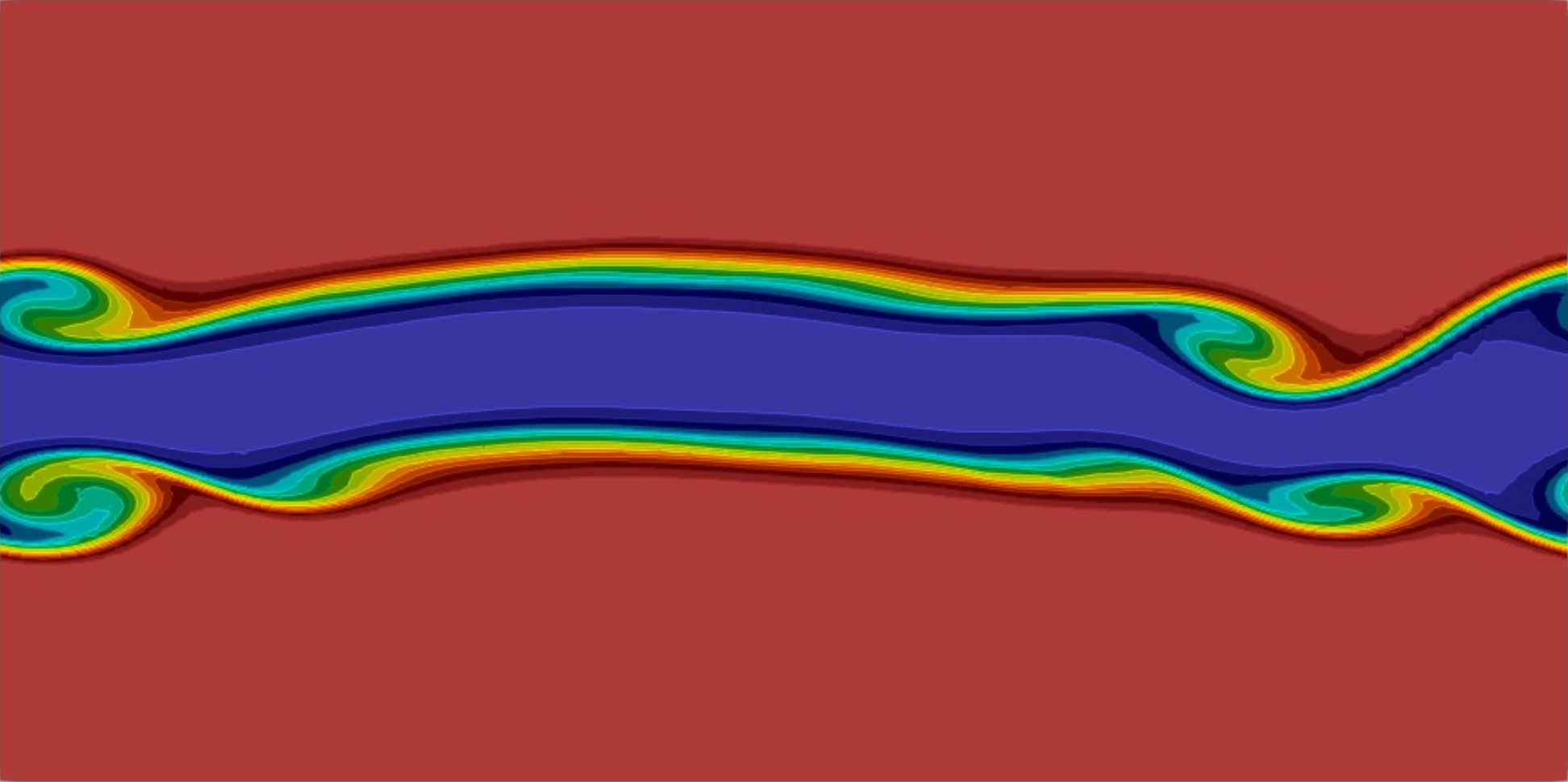}
    \includegraphics[trim=0.0cm 0.cm 0.0cm 0.0cm,clip=true,width=0.43\textwidth]{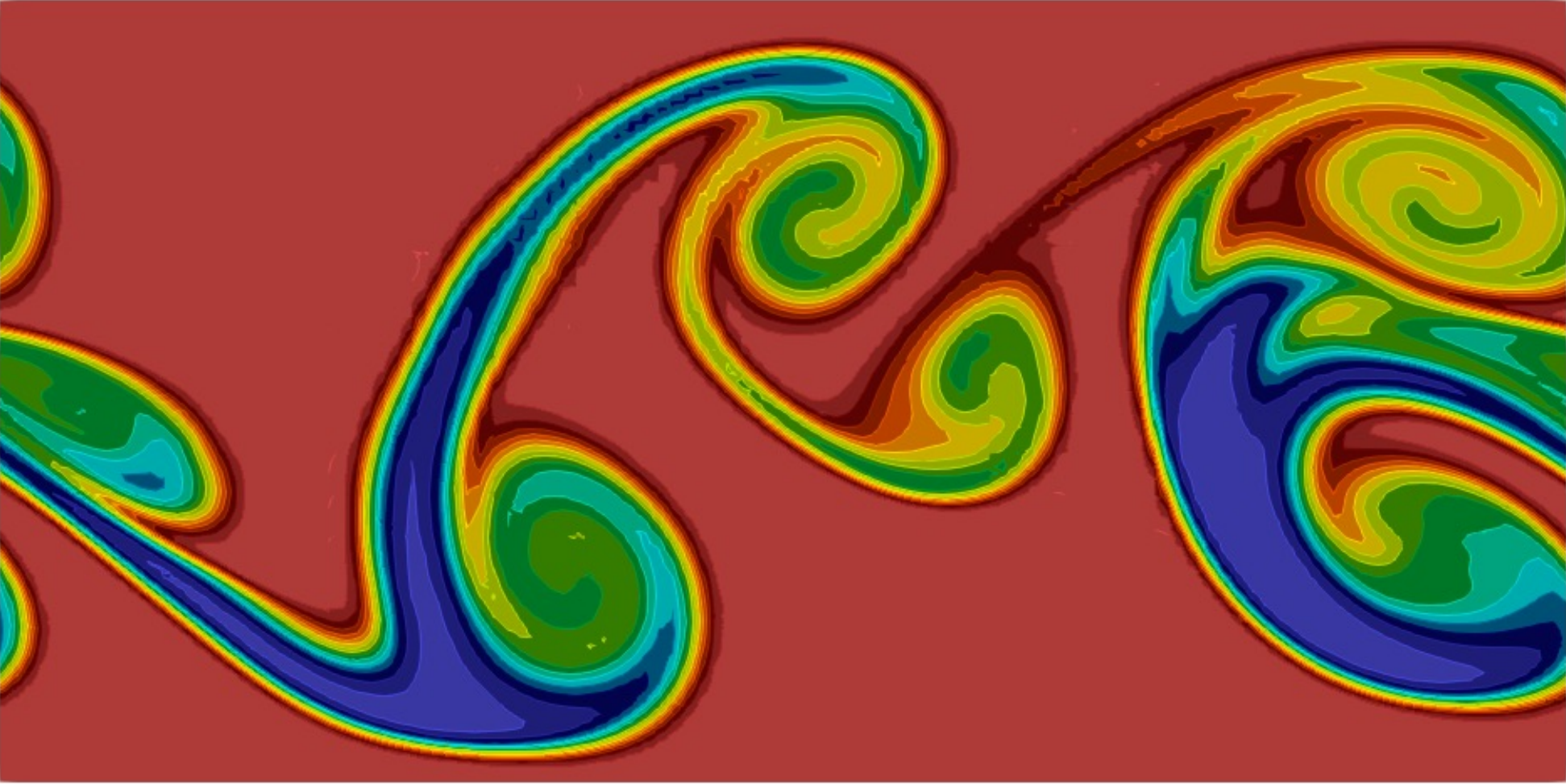}
    \figlab{expo-euler-khi-ss-epi2-x100}
  }  
  \subfigure[EPI2 with $\dt_{RK4} \times 8$]{
    \includegraphics[trim=0.0cm 0.cm 0.0cm 0.0cm,clip=true,width=0.43\textwidth]{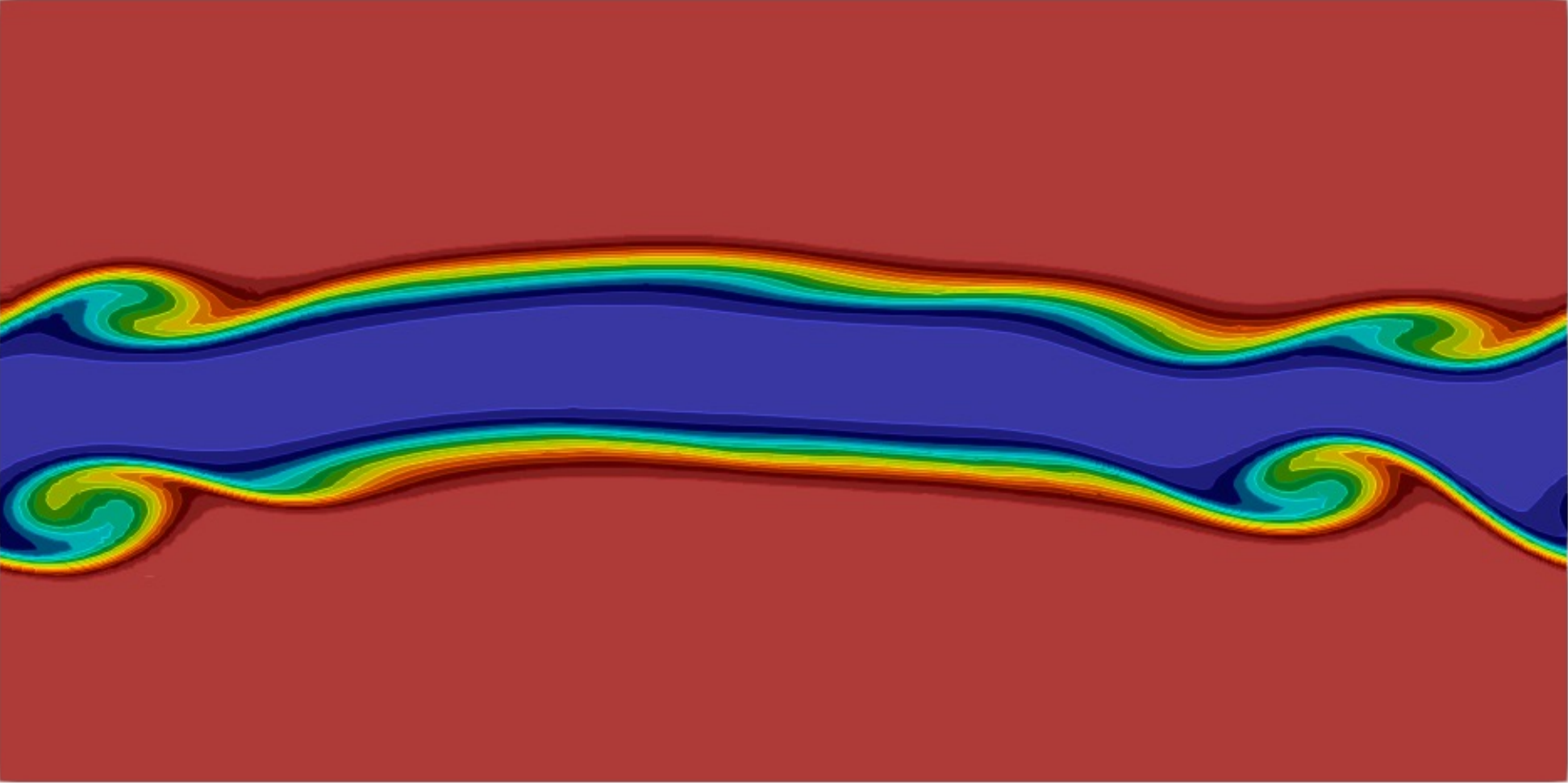}
    \includegraphics[trim=0.0cm 0.cm 0.0cm 0.0cm,clip=true,width=0.43\textwidth]{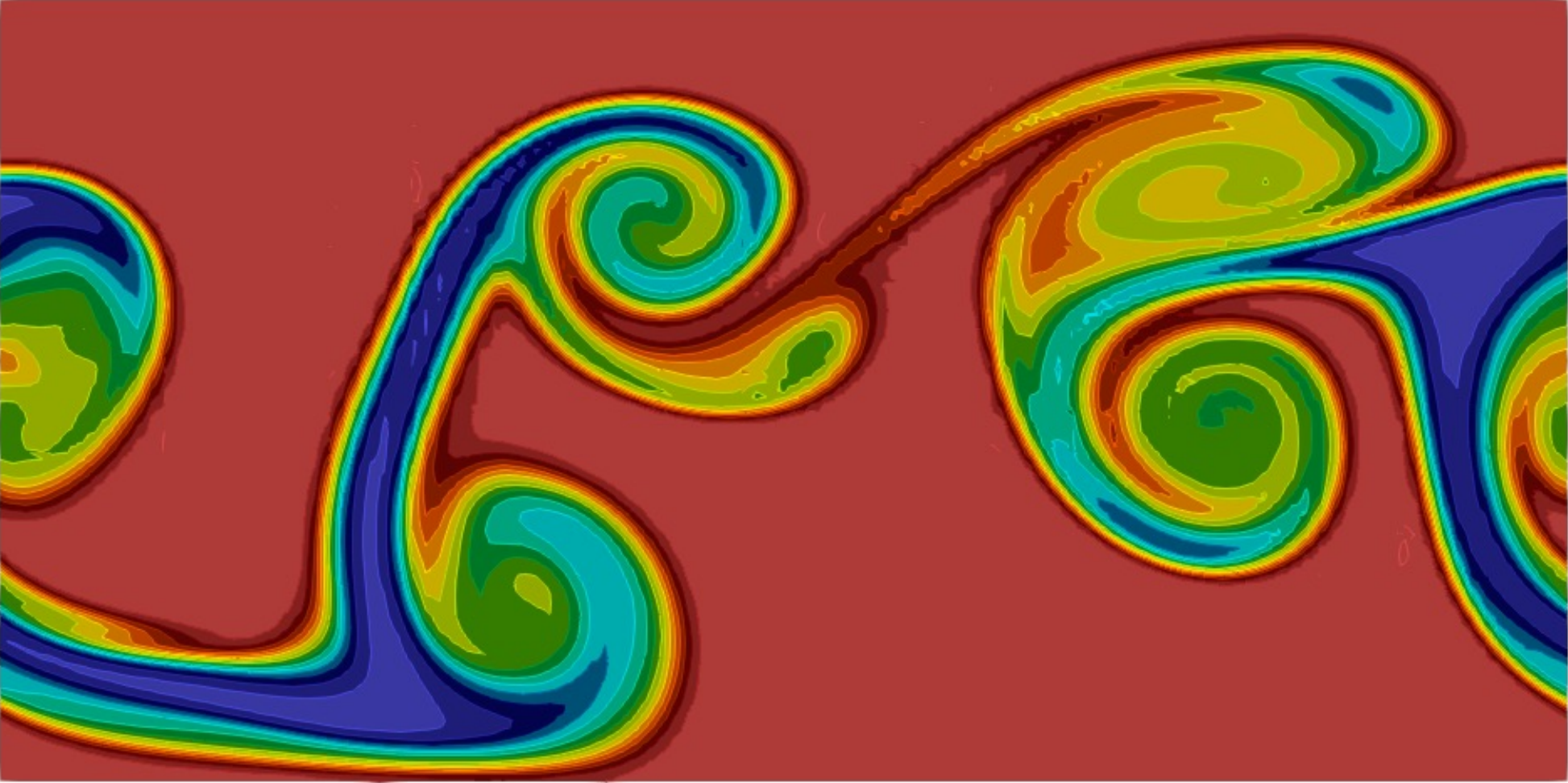}
    \figlab{expo-euler-khi-ss-epi2-x8}
  }  
  \subfigure[RK4]{
    \includegraphics[trim=0.0cm 0.cm 0.0cm 0.0cm,clip=true,width=0.43\textwidth]{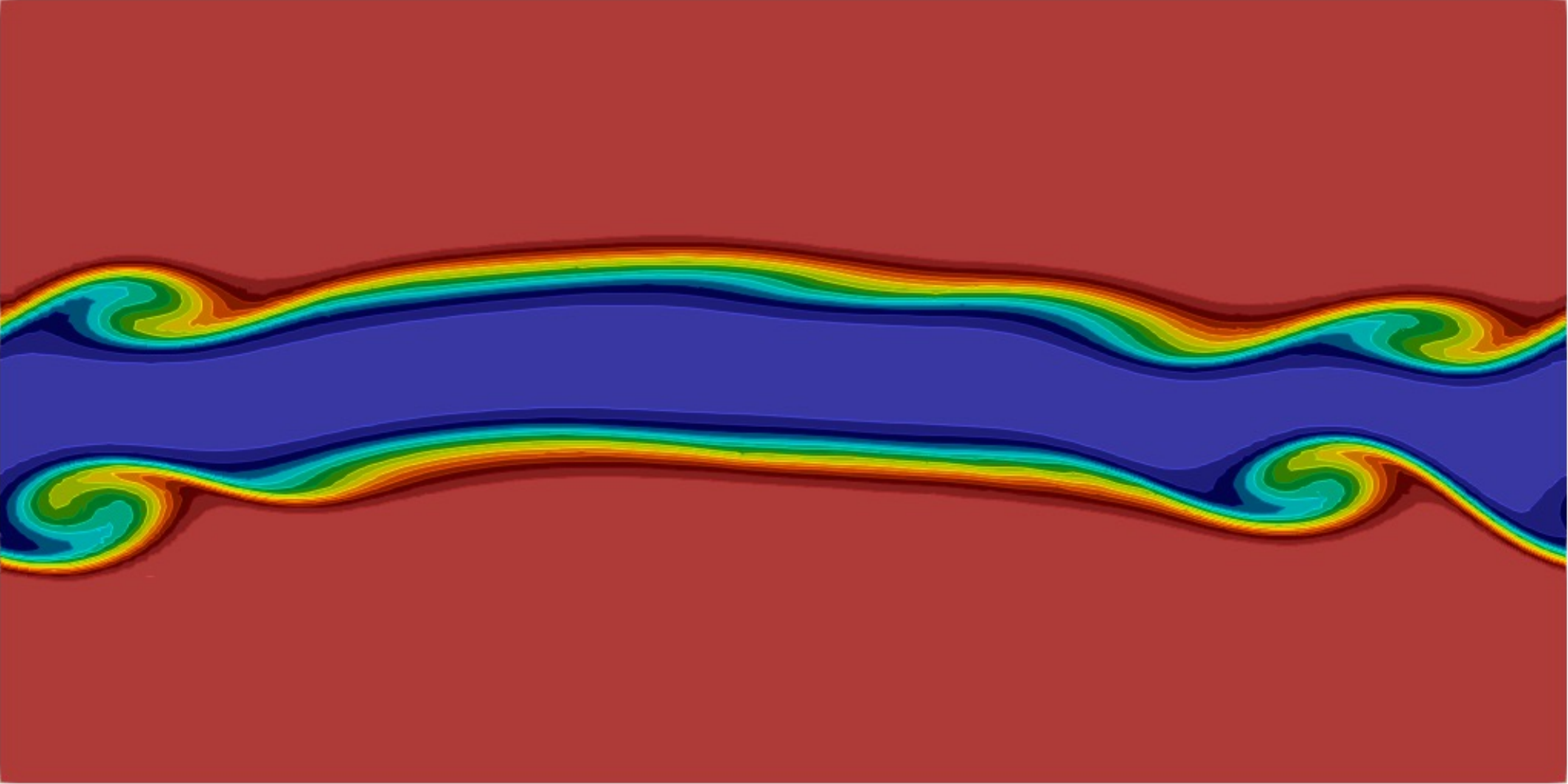}
    \includegraphics[trim=0.0cm 0.cm 0.0cm 0.0cm,clip=true,width=0.43\textwidth]{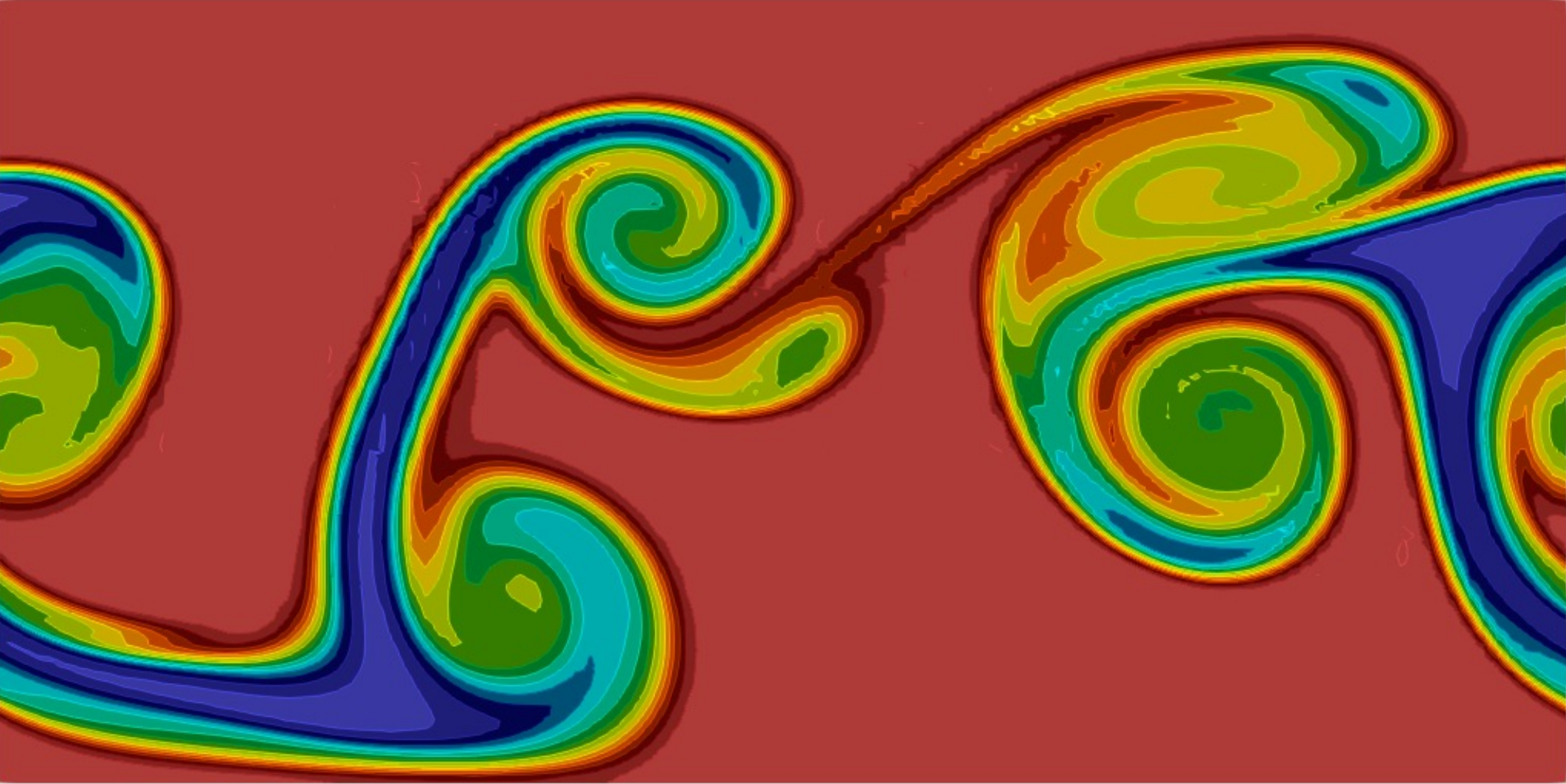}     
    \figlab{expo-euler-khi-ss-rk4}
  }
\caption{Kelvin--Helmholtz instability: evolution of temperature for (a) EPI2 with $dt=0.5(=\dt_{RK4}\times100)$, 
(b) EPI2 with $\dt=0.04(=\dt_{RK4}\times8)$ and (c) RK4 with $\dt=0.005(=:\dt_{RK4})$
at $t=\LRp{50,100}$ 
on a uniform mesh with $N_e = 722$ and $k=6$. 
The temperature ranges from $0.49$ to $1.01$. }
\figlab{euler-khi-snapshots}
\end{figure}


\subsection{Euler equations: shock problems}

Exponential-DG methods are now tested for shock problems
by considering benchmark examples in \cite{liska2003comparison}. 


\subsubsection{Riemann problem: case $\# 4$ }

This example develops four shocks. 
The initial condition with $\gamma=1.4$ is defined to be 
\begin{align*}
  \rho &= 1.1,    u = 0,      v = 0,      \pres = 1.1   & x > 0.5, y > 0.5, \\
  \rho &= 0.5065, u = 0.8939, v = 0,      \pres = 0.35  & x \le 0.5, y > 0.5, \\
  \rho &= 1.1,    u = 0.8939, v = 0.8939, \pres = 1.1   & x \le 0.5, y \le 0.5, \\
  \rho &= 0.5065, u = 0,      v = 0.8939, \pres = 0.35  & x > 0.5, y \le 0.5
\end{align*}
on $\Omega = (0,1)^2$. 
The control parameters of the entropy viscosity are $c_{\max}= 0.1$ and $c_E = 1$.

We conduct the numerical simulation 
with EPI2 and RK4 methods 
for $t\in\LRs{0,0.24}$ with $k=3$ and $N_e=5000$ 
in Figure \figref{euler-rm4-snapshots}. 
We take $\dt_{RK4}=0.0004$ for RK4 and $\dt=0.004(=\dt_{RK4}\times 10)$ for EPI2.
In general, EPI2 solution is comparable with RK4 counterpart. 
With 10 times larger time stepsize, EPI2 (taking $324s$) is 4 times faster than RK4 (taking $1366s$).

\begin{figure}[h!t!b!]
  \centering
  \subfigure[EPI2 with $\dt_{RK4} \times 10$]{
    \includegraphics[trim=0.0cm 0.cm 0.0cm 0.0cm,clip=true,width=0.43\textwidth]{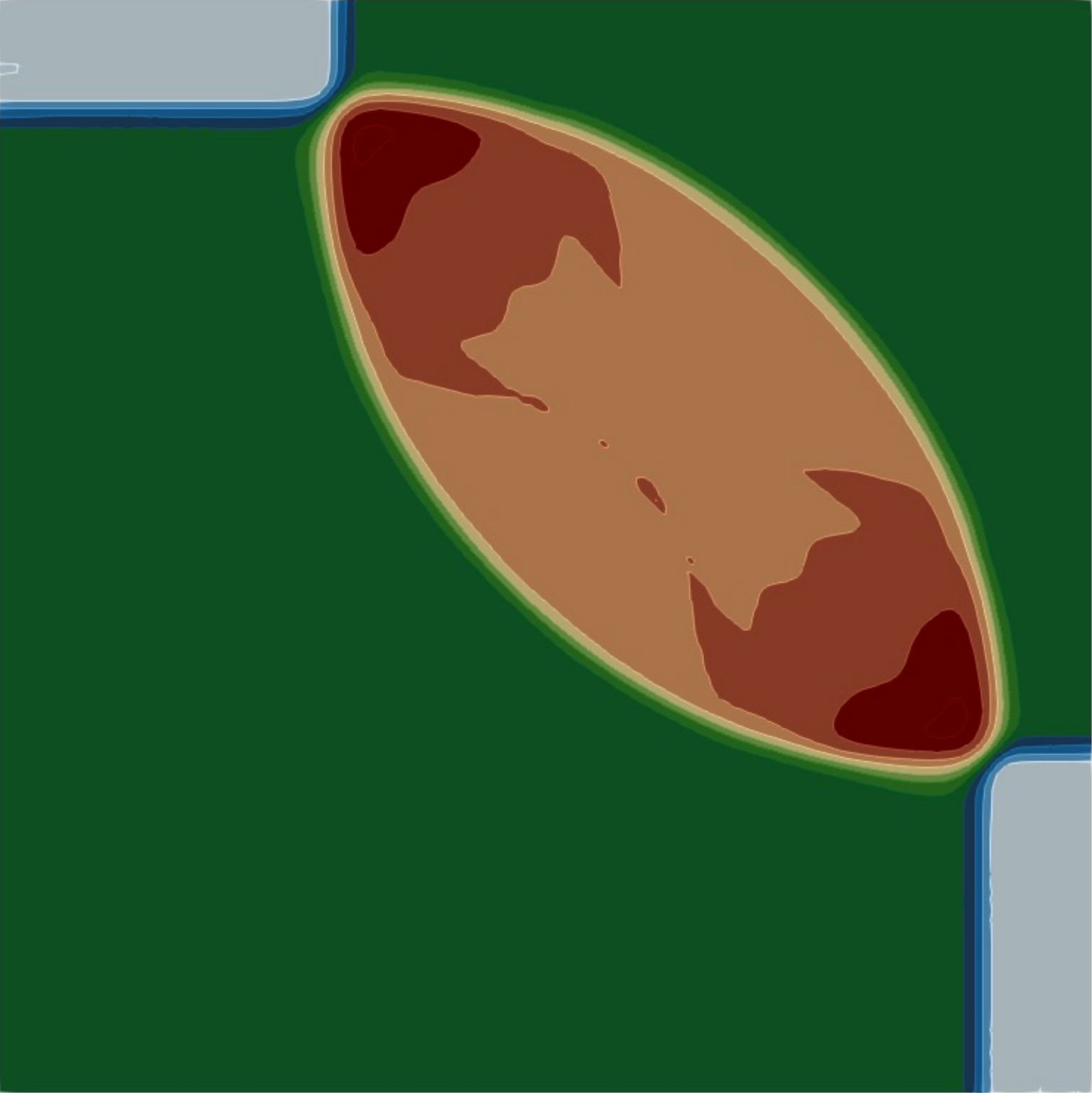}
    \figlab{expo-euler-rm4-ss-epi2}
  }  
  \subfigure[RK4 with $\dt_{RK4}=0.0002$]{
    \includegraphics[trim=0.0cm 0.cm 0.0cm 0.0cm,clip=true,width=0.43\textwidth]{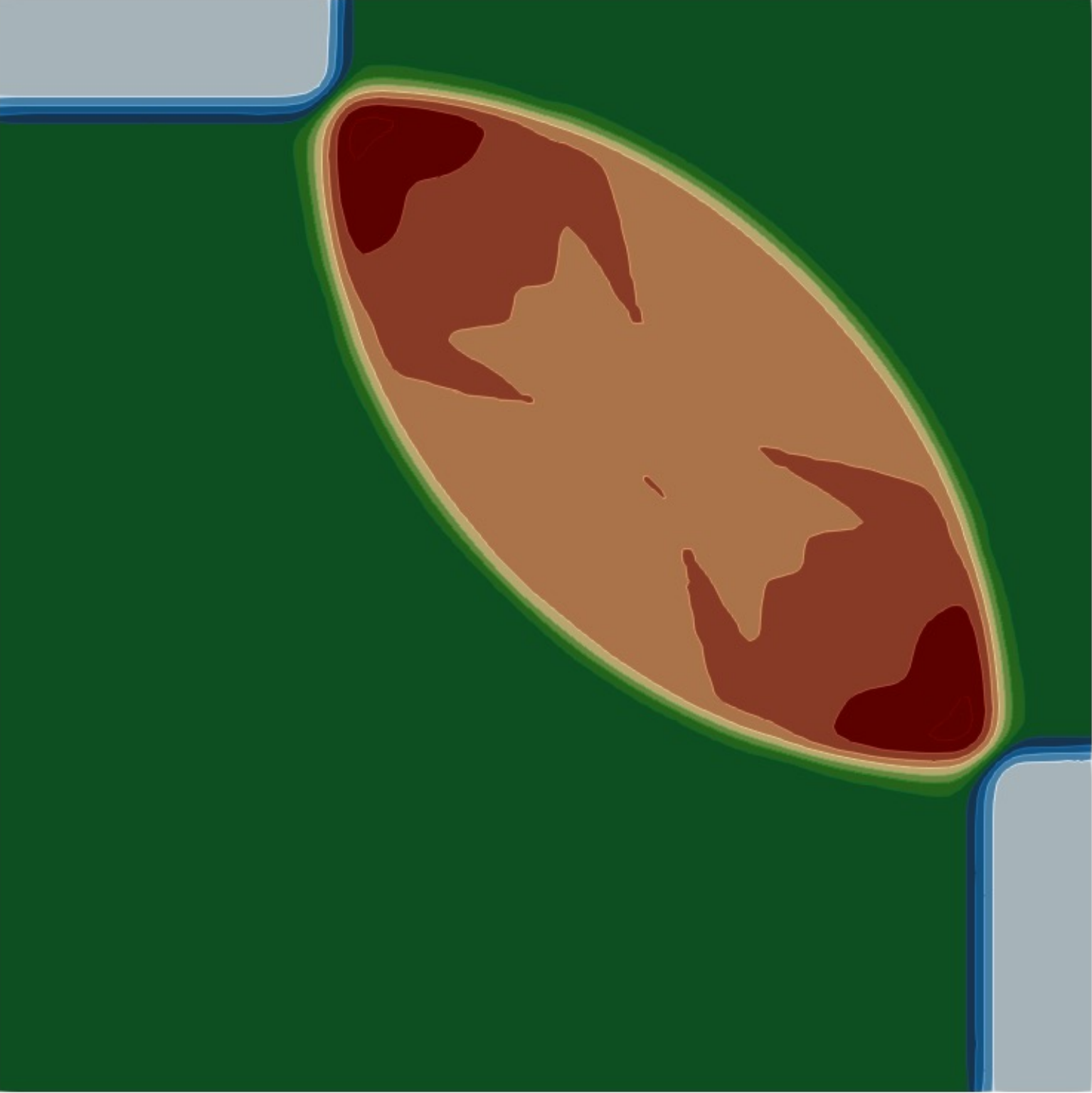}
    \figlab{expo-euler-rm4-ss-ssprk3}
  }  
\caption{Riemann problem: case $\# 4$: density field of 
(a) EPI2 with $\dt=0.002(=\dt_{RK4}\times10)$ and  
(b) RK4 with $\dt_{RK4}=0.0002$
at $t=0.25$ 
on a uniform mesh with $N_e = 5000$ and $k=3$. 
The density ranges from $0.5$ to $1.9$.}
\figlab{euler-rm4-snapshots}
\end{figure}

\subsubsection{Riemann problem: case $\# 12$ }

This example develops two contact waves and two shocks. 
The initial condition with $\gamma=1.4$ is given as 
\begin{align*}
  \rho &= 0.5313,  u = 0,      v = 0,      \pres = 0.4 & x > 0.5, y > 0.5, \\
  \rho &= 1.0,     u = 0.7276, v = 0,      \pres = 1   & x \le 0.5, y > 0.5, \\
  \rho &= 0.8,     u = 0,      v = 0,      \pres = 1   & x \le 0.5, y \le 0.5, \\
  \rho &= 1.0,     u = 0,      v = 0.7276, \pres = 1   & x > 0.5, y \le 0.5
\end{align*}
on $\Omega = (0,1)^2$. 
The control parameters of the entropy viscosity are $c_{\max}= 0.05$ and $c_E = 0.5$.

We conduct the numerical simulation for $t\in\LRs{0,0.24}$ with $k=3$ and $N_e=5000$ 
in Figure \figref{euler-rm12-snapshots}. 
We take $\dt_{RK4}=0.0004$ for RK4 and $\dt=0.004(=\dt_{RK4}\times 10$ for EPI2.
The wallclock times of RK4 and EPI2 are $1051.0s$ and $327.5s$, respectively. 
EPI2, with 10 times larger time stepsize, is about 3 times faster than RK4, 
and produces the comparable solution to RK4 counterpart.

Figure \figref{euler-rm12-snapshots-mu} shows 
the viscosity fields of EPI2 and RK4.
As mentioned in \cite{guermond2011entropy}, 
the viscosity becomes strong in the shocks, 
whereas weak in the rest including contact discontinuities. 
As expected, we also see that the magnitude of the viscosities for EPI2 and RK4 
are different due to the approximation of the entropy residual. 
Also, the entropy viscosity method 
does not completely remove the Gibbs phenomenon associated with high-order spatial discretizaton for shock problems. 
A further study is needed to handle the issue by incorporating several limiters.

\begin{figure}[h!t!b!]
  \centering
  \subfigure[EPI2 with $\dt_{RK4} \times 10$]{
    \includegraphics[trim=0.0cm 0.cm 0.0cm 0.0cm,clip=true,width=0.43\textwidth]{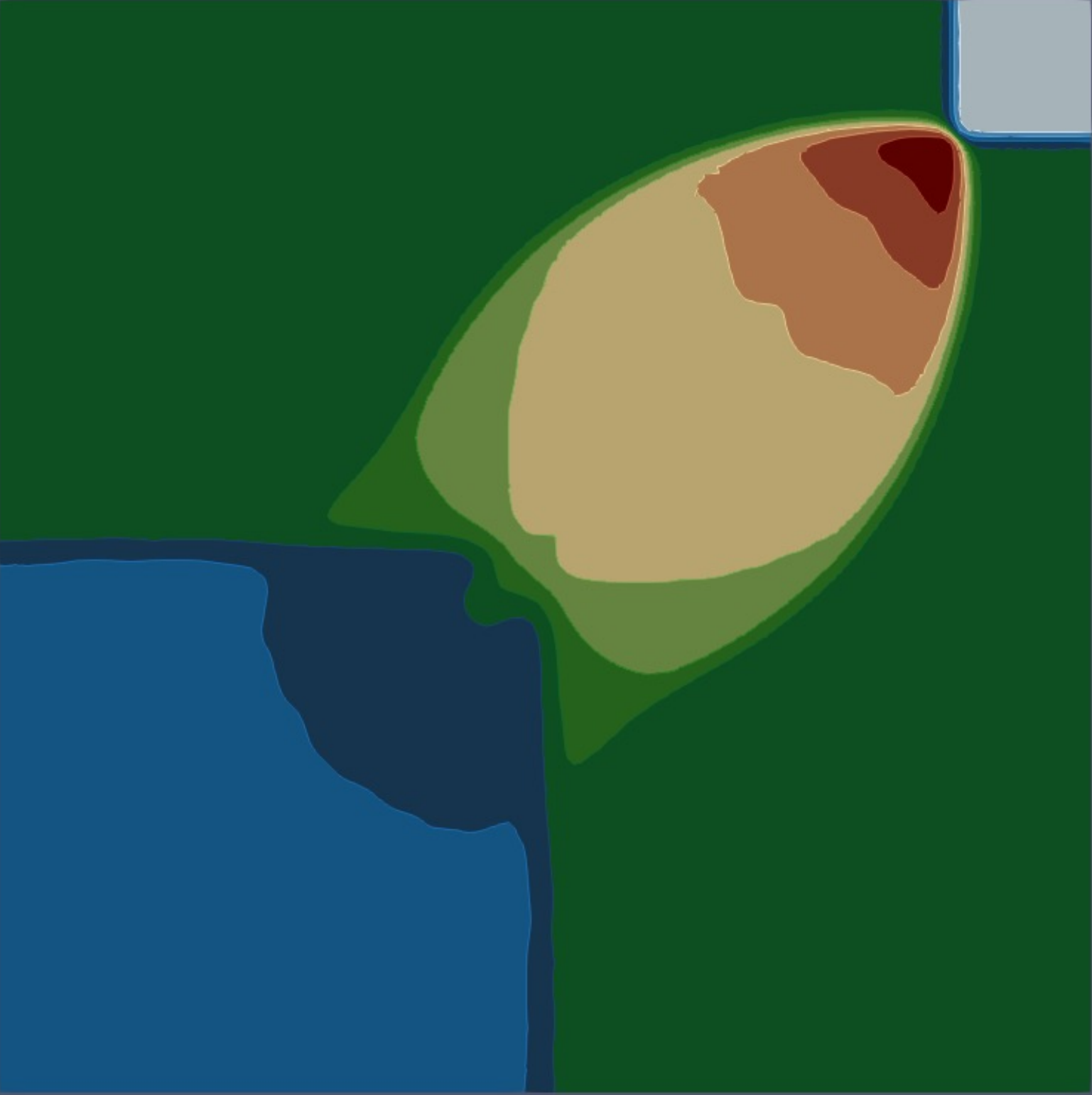}
    \figlab{expo-euler-rm12-ss-epi2}
  }  
  \subfigure[RK4 with $\dt_{RK4}=0.0002$]{
    \includegraphics[trim=0.0cm 0.cm 0.0cm 0.0cm,clip=true,width=0.43\textwidth]{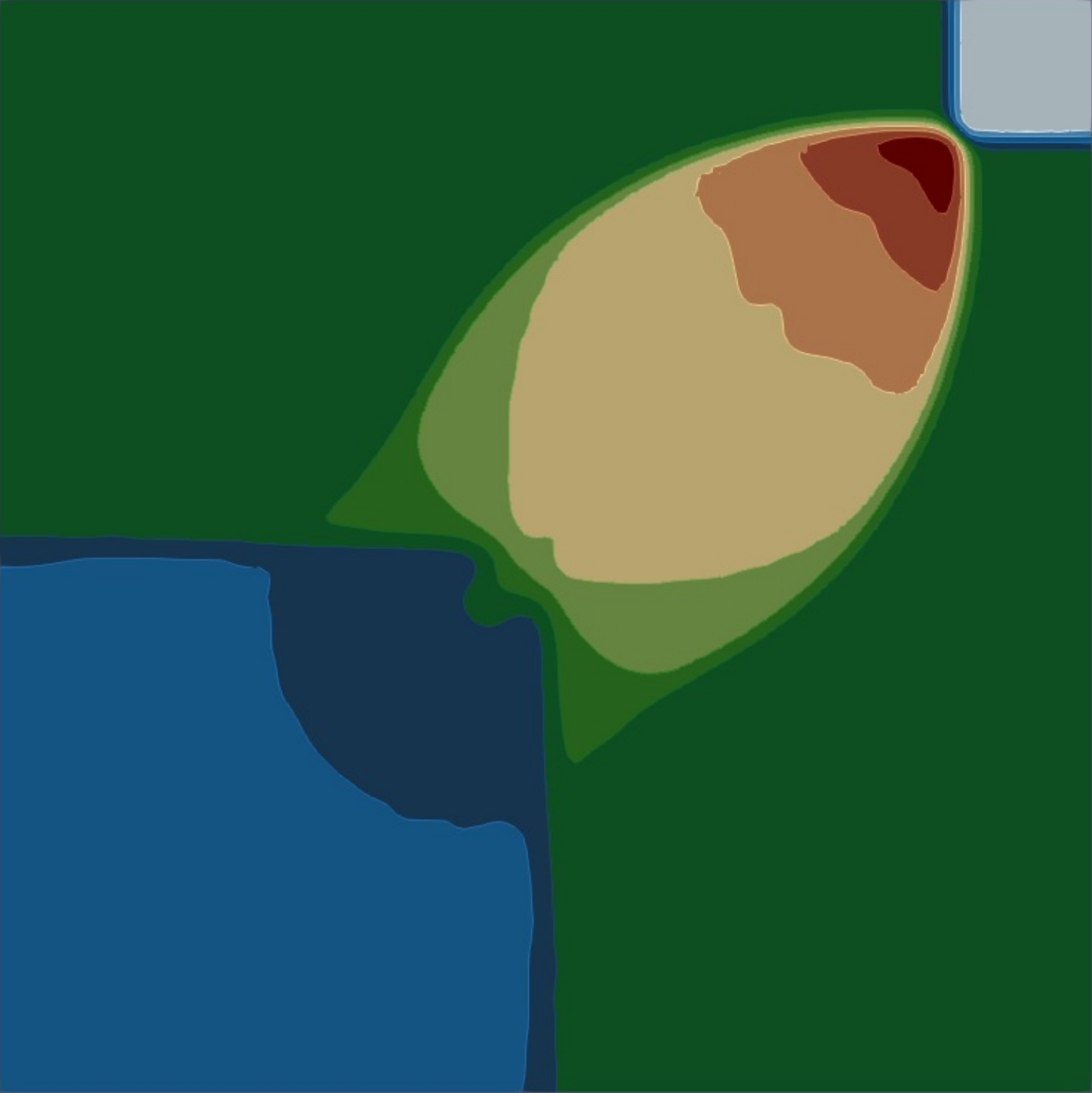}
    \figlab{expo-euler-rm12-ss-ssprk3}
  }  
\caption{Riemann problem: case $\# 12$: density field of 
(a) EPI2 with $\dt=0.004(=\dt_{RK4}\times10)$ and  
(b) RK4 with $\dt_{RK4}=0.0004$
at $t=0.24$ 
on a uniform mesh with $N_e = 5000$ and $k=3$.
The density ranges from $0.5$ to $1.7$. }
\figlab{euler-rm12-snapshots}
\end{figure}

\begin{figure}[h!t!b!]
  \centering
  \subfigure[EPI2 with $\dt_{RK4} \times 10$]{
    \includegraphics[trim=0.0cm 0.cm 0.0cm 0.0cm,clip=true,width=0.43\textwidth]{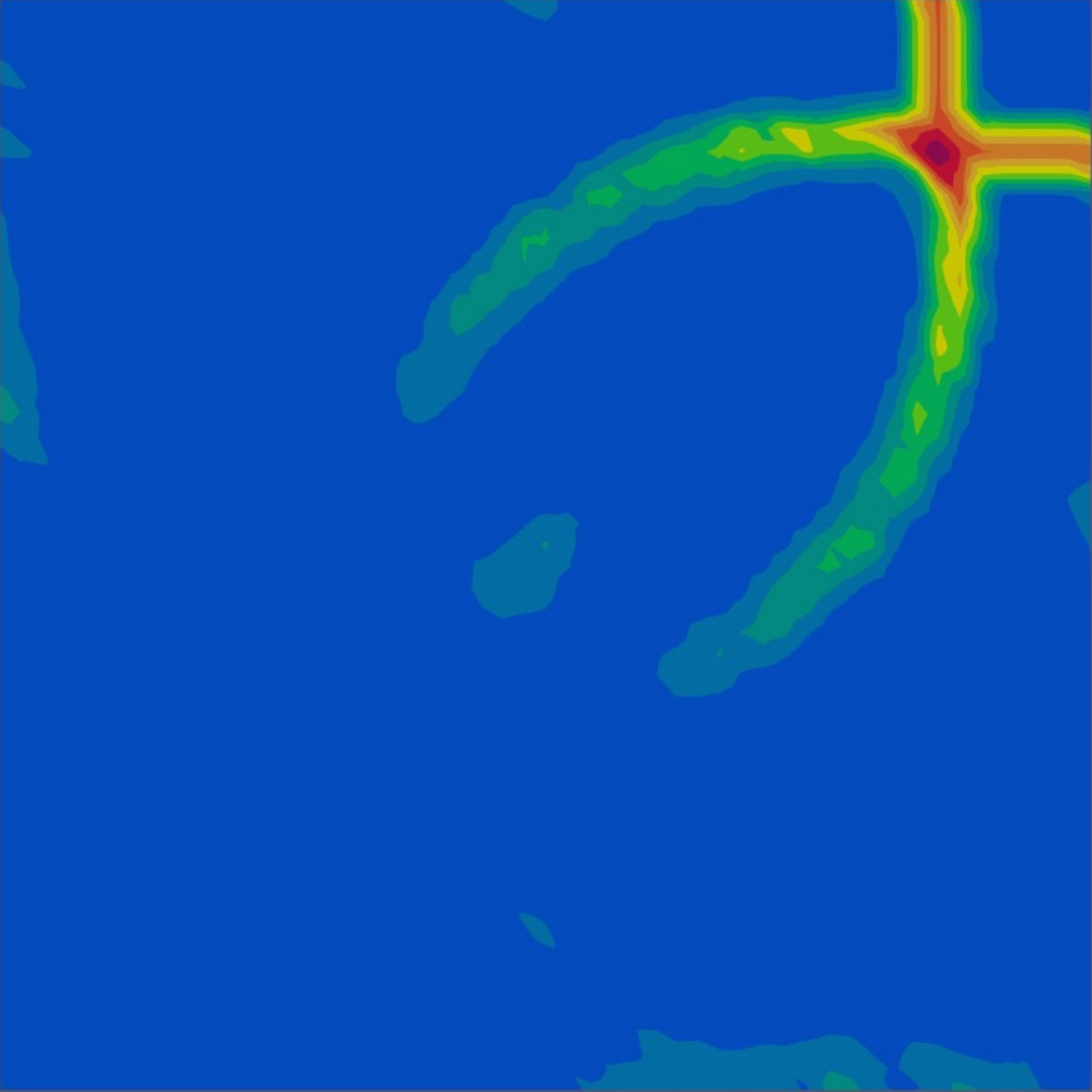}
    \figlab{expo-euler-rm12-ss-epi2-mu}
  }  
  \subfigure[RK4 with $\dt_{RK4}=0.0002$]{
    \includegraphics[trim=0.0cm 0.cm 0.0cm 0.0cm,clip=true,width=0.43\textwidth]{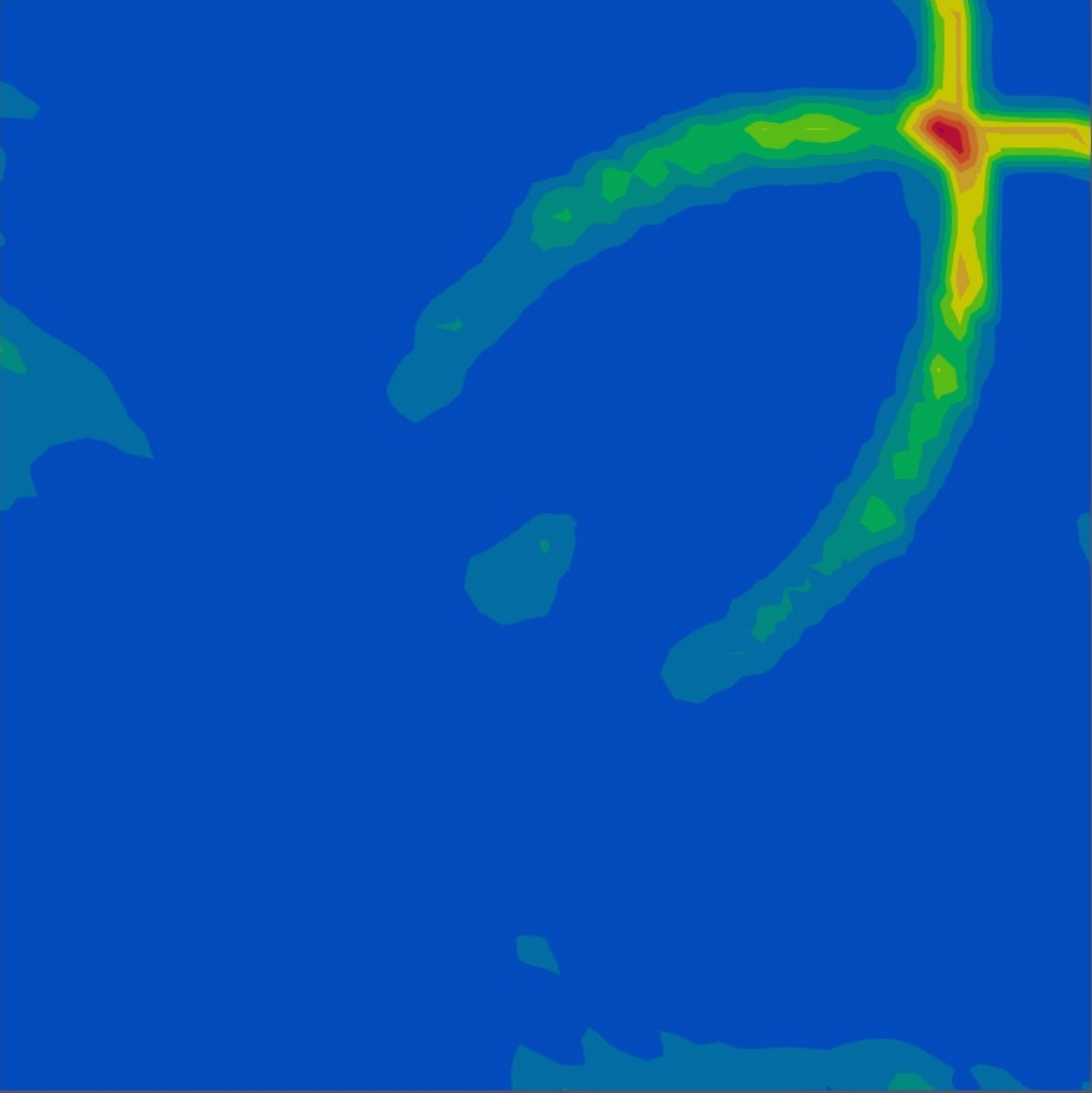}
    \figlab{expo-euler-rm12-ss-ssprk3-mu}
  }  
\caption{Riemann problem: case $\# 12$: viscosity field of 
(a) EPI2 with $\dt=0.004(=\dt_{RK4}\times10)$ and  
(b) RK4 with $\dt_{RK4}=0.0004$
at $t=0.24$ 
on a uniform mesh with $N_e = 5000$ and $k=3$.
The viscosity ranges from $0$ to $0.0035$. }
\figlab{euler-rm12-snapshots-mu}
\end{figure}

\section{Conclusions}
\seclab{Conclusion}
  In this paper, we have developed a Exponential DG framework.
  This is done by splitting the governing differential operator into linear and nonlinear parts 
  to which we apply DG spatial discretization. 
  In particular, we construct the linear part by linearization
  aiming to absorb the stiffness in the system.  
    Since the linear-nonlinear decomposition is done on continuous level, we can avoid taking derivatives of  nonsmooth functions possibly resulting from both spatial and time discretizations.
The resulting semi-discrete system is then intetegrated with exponential integrators.
  Our proposed approach aims to 
  i) circumvent the stringent timestep size arising from explicit integrators; 
  ii) support high-order accuracy in both space and time;
  iii) outperform over IMEX DG methods 
   with no preconditioner; 
  iv) be comparable to explicit RKDG methods for stiff problems; 
  v) be scalable in a modern massively parallel computing architecture.
  We present a detailed stability and convergence analyses for the Burgers equation using the exponential Euler DG scheme.
  
  Numerical results (for Burgers equation and Euler equations) have shown that while explicit RKDG methods suffer from restricted timestep sizes due to numerical stability, 
  Exponential DG framework supports a wide range of Courant numbers.
  We numerically observe that the proposed methods achieve 
  the high-order temporal and spatial convergence rates.
  We also see that Exponential DG is more economical than IMEX DG 
  in the isentropic vortex example on both uniform and non-uniform meshes. 
  For the Euler systems on the non-uniform mesh,
  Exponential DG is comparable to explicit RKDG.
  Moreover, for the shock problems, 
  EPI2 solutions become 3 times faster than
  RK4 solutions when artificial viscosity is employed.
  This is because the diffusion term becomes a dominant source to restrict the timestep size of the explicit methods in the shock problems. For all cases,  if relaxing the accuracy is allowed, while time stepsize beyond the maximum stable time stepsize for explicit RKDG is needed, 
  Exponential DG can be faster than the explicit RKDG.
   
  As have been demonstrated, our proposed framework can exploit current and future parallel computing systems  to solve large scale problems. 
  The key explored in the proposed methods do not require a linear solve
  matrix-free Krylov-based matrix exponential computations
  and the DG compact communication stencil. 
  Indeed, we have numerically shown that 
  Exponential DG methods have favorable strong scaling up to $40K$ cores and weak scaling for $16K$ cores for the Euler isentropic vortex example.
  Ongoing work is to extend the approach to various partial differential equations 
and to scale it beyond hundreds of thousands cores.
\appendix

\section{Local discontinuous Galerkin methods for viscous Burgers equation}
\seclab{ldg-burgers}

We have seen that suboptimal convergence rates for odd orders 
in Table \tabref{expo-burgers-sconv-mms}, 
Table \tabref{expo-burgers-sconv-smooth} and  
Table \tabref{expo-burgers-sconv-shock}.
This is related to the use of central fluxes in the diffusion term. 
To improve a spatial convergence rate, we employ 
local discontinuous Galerkin methods (LDG) \cite{cockburn1998local}.
That is, we define the numerical flux $u^{**}$ in \eqnref{gov-burgers1d-gradu} 
and $q^{**}$ in \eqnref{BurgerL} by 
\begin{align*}
  u^{**} = \average{\u} - \jump{\u}\beta, \quad 
  q^{**} = \average{\q} + \jump{\q}\beta, 
\end{align*}
where we take $\beta=0.5$.
\footnote{
  Lax-Friedrich fluxes are used for $(\tilde{u}\u)^*$ and $\LRp{\frac{u^2}{2}}^*$ in \eqnref{BurgerN}.
}
Indeed, we numerically observe that 
the spatial convergence rates increase to $k+1$ for odd orders 
in the case with the time-independent smooth solution 
as shown in Table \tabref{expo-burgers-sconv-mms-ldg}.
We also see that the spatial convergence rates for odd orders are improved 
up to $0.45$ compared with central flux for time-dependent problems in 
Table \tabref{expo-burgers-sconv-smooth-ldg} and  
Table \tabref{expo-burgers-sconv-shock-ldg}. 
The spatial convergence results are encouraging, thus, 
we will consider to incorporate LDG methods for developing Naiver--Stokes models in the future. 

\begin{table}[t]
  \caption{A time-independent manufactured solution for the viscous Burgers equation:
  a spatial convergence study using $N_e=\LRc{20,40,80,160}$ elements is conducted 
  with central (C) flux  and local discontinous Galkerin (LDG) flux for diffusion term.
  Lax-Friedrich (LF) flux is used for advective term. 
  EXPRB32 scheme with $\dt=5\times10^{-5}$ is used as the time integrator.
  }
  \tablab{expo-burgers-sconv-mms-ldg}
  \begin{center}
    \begin{tabular}{*{1}{c}|*{1}{c}|*{2}{c}|*{2}{c}}
      \hline
      \multirow{2}{*}{ } 
        & \multirow{2}{*}{$h$}  
        & \multicolumn{2}{c|}{LF (C)}   
        & \multicolumn{2}{c|}{LF (LDG)}   
        \tabularnewline
        &    &  error & order &  error & order       \tabularnewline
      \hline\hline
      \multirow { 4}{*}{ $k=1$ }& 1/20  & 4.093E-04 &      $-$ & 4.415E-04 &      $-$ \tabularnewline
                               &  1/40  & 1.223E-04 &    1.743 & 1.112E-04 &    1.990 \tabularnewline
                               &  1/80  & 4.494E-05 &    1.445 & 2.782E-05 &    1.998 \tabularnewline
                               &  1/160 & 1.937E-05 &    1.214 & 6.959E-06 &    1.999 \tabularnewline
      \multicolumn{6}{c}{} \tabularnewline
      \multirow { 4}{*}{ $k=2$ }& 1/20  & 2.630E-06 &      $-$ & 3.586E-06 &      $-$ \tabularnewline
                               &  1/40  & 3.210E-07 &    3.034 & 4.635E-07 &    2.952 \tabularnewline
                               &  1/80  & 3.966E-08 &    3.017 & 5.847E-08 &    2.987 \tabularnewline
                               &  1/160 & 4.916E-09 &    3.012 & 7.364E-09 &    2.989 \tabularnewline
      \multicolumn{6}{c}{} \tabularnewline
      \multirow { 4}{*}{ $k=3$ }& 1/20  & 1.431E-07 &      $-$ & 6.185E-08 &      $-$ \tabularnewline
                               &  1/40  & 1.709E-08 &    3.066 & 3.713E-09 &    4.058 \tabularnewline
                               &  1/80  & 2.003E-09 &    3.093 & 2.270E-10 &    4.032 \tabularnewline
                               &  1/160 & 2.251E-10 &    3.154 & 1.404E-11 &    4.015 \tabularnewline
      \multicolumn{6}{c}{} \tabularnewline
      \multirow { 4}{*}{ $k=4$ }& 1/20  & 5.946E-10 &      $-$ & 8.475E-10 &      $-$ \tabularnewline
                               &  1/40  & 1.827E-11 &    5.024 & 2.593E-11 &    5.030 \tabularnewline
                               &  1/80  & 5.626E-13 &    5.021 & 8.027E-13 &    5.014 \tabularnewline
                               &  1/160 & 1.730E-14 &    5.023 & 2.498E-14 &    5.006 \tabularnewline
      \hline\hline
      \end{tabular}
  \end{center}
  \end{table}

\begin{table}[t]
  \caption{A smooth solution for the viscous Burgers equation:
  a spatial convergence study using  $N_e=\LRc{20,40,80,160}$ elements is conducted 
  with central (C) flux  and local discontinous Galkerin (LDG) flux for diffusion term.
  Lax-Friedrich (LF) flux is used for advective term. 
  EXPRB32 scheme with $\dt=5\times 10^{-5}$ as the time integrator.}
  \tablab{expo-burgers-sconv-smooth-ldg}
  \begin{center}
    \begin{tabular}{*{1}{c}|*{1}{c}|*{2}{c}|*{2}{c}}
      \hline
      \multirow{2}{*}{ } 
        & \multirow{2}{*}{$h$}  
        & \multicolumn{2}{c|}{LF (C)}   
        & \multicolumn{2}{c|}{LF (LDG)}   
        \tabularnewline
        &    &  error & order &  error & order       \tabularnewline
      \hline\hline
      \multirow { 4}{*}{ $k=1$ }& 1/20 & 7.061E-03 &      $-$ & 7.234E-03 &      $-$ \tabularnewline
                               &  1/40 & 2.080E-03 &    1.764 & 1.901E-03 &    1.928 \tabularnewline
                               &  1/80 & 6.731E-04 &    1.627 & 4.928E-04 &    1.948 \tabularnewline
                               &  1/160 & 2.432E-04 &    1.469 & 1.297E-04 &    1.926 \tabularnewline
      \multicolumn{6}{c}{} \tabularnewline
      \multirow { 4}{*}{ $k=2$ }&  1/20 & 2.511E-04 &      $-$ & 3.203E-04 &      $-$ \tabularnewline
                               &  1/40 & 2.741E-05 &    3.196 & 4.035E-05 &    2.989 \tabularnewline
                               &  1/80 & 3.317E-06 &    3.047 & 5.065E-06 &    2.994 \tabularnewline
                               &  1/160 & 4.112E-07 &    3.012 & 6.346E-07 &    2.997 \tabularnewline
      \multicolumn{6}{c}{} \tabularnewline
      \multirow { 4}{*}{ $k=3$ }&  1/20 & 2.574E-05 &      $-$ & 2.508E-05 &      $-$ \tabularnewline
                               &  1/40 & 2.883E-06 &    3.158 & 2.321E-06 &    3.434 \tabularnewline
                               &  1/80 & 3.414E-07 &    3.078 & 2.172E-07 &    3.417 \tabularnewline
                               &  1/160 & 4.108E-08 &    3.055 & 2.002E-08 &    3.439 \tabularnewline
      \multicolumn{6}{c}{} \tabularnewline
      \multirow { 4}{*}{ $k=4$ }&  1/20 & 9.225E-07 &      $-$ & 1.494E-06 &      $-$ \tabularnewline
                               &  1/40 & 2.529E-08 &    5.189 & 8.022E-08 &    4.219 \tabularnewline
                               &  1/80 & 7.135E-10 &    5.147 & 4.164E-09 &    4.268 \tabularnewline
                               &  1/160 & 2.188E-11 &    5.027 & 1.990E-10 &    4.387 \tabularnewline
      \hline\hline
      \end{tabular}
  \end{center}
  \end{table}

  \begin{table}[h!t!b!]
    \caption{A shock solution to the viscous Burgers equation with 
    a spatial convergence study using  $N_e=\LRc{40,80,160}$ elements is conducted 
    with central (C) flux  and local discontinous Galkerin (LDG) flux for diffusion term.
    Lax-Friedrich (LF) flux is used for advective term. 
    EXPRB32 scheme with $\dt=5\times 10^{-5}$ as the time integrator.
    We take the RK4 solution (with $\dt=5\times 10^{-7}$,  $k=10$, and
    $N_e=160$) as a reference solution to measure the $L^2$ error at
    $t=1$. 
    }
    \tablab{expo-burgers-sconv-shock-ldg}
    \begin{center}
      \begin{tabular}{*{1}{c}|*{1}{c}|*{2}{c}|*{2}{c}}
        \hline
        \multirow{2}{*}{ } 
          & \multirow{2}{*}{$h$}  
          & \multicolumn{2}{c|}{LF (C)}   
          & \multicolumn{2}{c|}{LF (LDG)}   
          \tabularnewline
          &    &  error & order &  error & order       \tabularnewline
        \hline\hline
        \multirow { 3}{*}{ $k=1$ }& 1/40  & 1.295E-02 &      $-$ & 1.517E-02 &      $-$ \tabularnewline
                                 &  1/80  & 5.558E-03 &    1.221 & 5.960E-03 &    1.347 \tabularnewline
                                 &  1/160 & 2.001E-03 &    1.474 & 1.838E-03 &    1.698 \tabularnewline
        \multicolumn{6}{c}{} \tabularnewline
        \multirow { 3}{*}{ $k=2$ }&  1/40 & 4.223E-03 &      $-$ & 4.273E-03 &      $-$ \tabularnewline
                                 &   1/80 & 6.722E-04 &    2.651 & 7.225E-04 &    2.564 \tabularnewline
                                 &   1/160& 9.852E-05 &    2.770 & 1.029E-04 &    2.812 \tabularnewline
        \multicolumn{6}{c}{} \tabularnewline
        \multirow { 3}{*}{ $k=3$ }&  1/40 & 9.976E-04 &      $-$ & 1.159E-03 &      $-$ \tabularnewline
                                 &   1/80 & 1.435E-04 &    2.797 & 1.672E-04 &    2.793 \tabularnewline
                                 &   1/160& 1.224E-05 &    3.551 & 1.065E-05 &    3.972 \tabularnewline
        \multicolumn{6}{c}{} \tabularnewline
        \multirow { 3}{*}{ $k=4$ }&  1/40 & 3.606E-04 &      $-$ & 4.442E-04 &      $-$ \tabularnewline
                                 &   1/80 & 2.728E-05 &    3.724 & 2.256E-05 &    4.299 \tabularnewline
                                 &   1/160& 8.084E-07 &    5.077 & 7.304E-07 &    4.949 \tabularnewline
        \hline\hline
        \end{tabular}
    \end{center}
    \end{table}

\section*{Acknowledgements}
SK was partially supported by 
  the U.S. Department of Energy, Office of Science, Advanced Scientific Computing Research Program under contract DE-AC02-06CH11357. The work of TBT was funded in part by DOE grants DE-SC0010518 and DE-SC0011118, NSF Grant NSF-DMS1620352. We are grateful for the supports.


\bibliographystyle{elsarticle-num}
\bibliography{main}

\begin{thebibliography}{10}
\expandafter\ifx\csname url\endcsname\relax
  \def\url#1{\texttt{#1}}\fi
\expandafter\ifx\csname urlprefix\endcsname\relax\def\urlprefix{URL }\fi
\expandafter\ifx\csname href\endcsname\relax
  \def\href#1#2{#2} \def\path#1{#1}\fi

\bibitem{ReedHill73}
W.~H. Reed, T.~R. Hill, Triangular mesh methods for the neutron transport
  equation, Tech. Rep. LA-UR-73-479, Los Alamos Scientific Laboratory (1973).

\bibitem{LeSaintRaviart74}
P.~LeSaint, P.~A. Raviart, On a finite element method for solving the neutron
  transport equation, in: C.~de~Boor (Ed.), Mathematical Aspects of Finite
  Element Methods in Partial Differential Equations, Academic Press, 1974, pp.
  89--145.

\bibitem{johnson1986analysis}
C.~Johnson, J.~Pitk{\"a}ranta, An analysis of the discontinuous {G}alerkin
  method for a scalar hyperbolic equation, Mathematics of computation 46~(173)
  (1986) 1--26.

\bibitem{wheeler1978elliptic}
M.~F. Wheeler, An elliptic collocation-finite element method with interior
  penalties, SIAM Journal on Numerical Analysis 15~(1) (1978) 152--161.

\bibitem{arnold1982interior}
D.~N. Arnold, An interior penalty finite element method with discontinuous
  elements, SIAM journal on numerical analysis 19~(4) (1982) 742--760.

\bibitem{cockburn2000development}
B.~Cockburn, G.~E. Karniadakis, C.-W. Shu, The development of discontinuous
  {G}alerkin methods, in: Discontinuous Galerkin Methods, Springer, 2000, pp.
  3--50.

\bibitem{arnold2002unified}
D.~N. Arnold, F.~Brezzi, B.~Cockburn, L.~D. Marini, Unified analysis of
  discontinuous {G}alerkin methods for elliptic problems, SIAM journal on
  numerical analysis 39~(5) (2002) 1749--1779.

\bibitem{liu2004discontinuous}
R.~Liu, Discontinuous {G}alerkin finite element solution for poromechanics,
  Ph.D. thesis, The University of Texas at Austin (2004).

\bibitem{nair2005discontinuous}
R.~D. Nair, S.~J. Thomas, R.~D. Loft, A discontinuous {G}alerkin global shallow
  water model, Monthly Weather Review 133~(4) (2005) 876--888.

\bibitem{GiraldoWarburton08}
F.~X. Giraldo, T.~Warburton, A high-order triangular discontinous {G}alerkin
  oceanic shallow water model, International Journal For Numerical Methods In
  Fluids 56 (2008) 899--925.

\bibitem{GiraldoRestelli10}
F.~X. Giraldo, M.~Restelli, High-order semi-implicit time-integrators for a
  triangular discontinous {G}alerkin oceanic shallow water model, International
  Journal For Numerical Methods In Fluids 63 (2010) 1077--1102.

\bibitem{wintermeyer2017entropy}
N.~Wintermeyer, A.~R. Winters, G.~J. Gassner, D.~A. Kopriva, An entropy stable
  nodal discontinuous {G}alerkin method for the two dimensional shallow water
  equations on unstructured curvilinear meshes with discontinuous bathymetry,
  Journal of Computational Physics 340 (2017) 200--242.

\bibitem{bassi2005discontinuous}
F.~Bassi, A.~Crivellini, S.~Rebay, M.~Savini, Discontinuous {G}alerkin solution
  of the {R}eynolds-averaged {N}avier--{S}tokes and k--$\omega$ turbulence
  model equations, Computers and Fluids 34~(4-5) (2005) 507--540.

\bibitem{gassner2016split}
G.~J. Gassner, A.~R. Winters, D.~A. Kopriva, Split form nodal discontinuous
  {G}alerkin schemes with summation-by-parts property for the compressible
  {E}uler equations, Journal of Computational Physics 327 (2016) 39--66.

\bibitem{fezoui2005convergence}
L.~Fezoui, S.~Lanteri, S.~Lohrengel, S.~Piperno, Convergence and stability of a
  discontinuous {G}alerkin time-domain method for the 3d heterogeneous
  {M}axwell equations on unstructured meshes, ESAIM: Mathematical Modelling and
  Numerical Analysis 39~(6) (2005) 1149--1176.

\bibitem{bui2012analysis}
T.~Bui-Thanh, O.~Ghattas, Analysis of an hp-nonconforming discontinuous
  {G}alerkin spectral element method for wave propagation, SIAM Journal on
  Numerical Analysis 50~(3) (2012) 1801--1826.

\bibitem{noels2008explicit}
L.~Noels, R.~Radovitzky, An explicit discontinuous {G}alerkin method for
  non-linear solid dynamics: Formulation, parallel implementation and
  scalability properties, International Journal for Numerical Methods in
  Engineering 74~(9) (2008) 1393--1420.

\bibitem{tirupathi2015modeling}
S.~Tirupathi, J.~S. Hesthaven, Y.~Liang, Modeling 3d magma dynamics using a
  discontinuous {G}alerkin method, Communications in Computational Physics
  18~(1) (2015) 230--246.

\bibitem{demkowicz2006computing}
L.~Demkowicz, Computing with hp-adaptive finite elements: volume 1 one and two
  dimensional elliptic and {M}axwell problems, Chapman and Hall/CRC, 2006.

\bibitem{kanevsky2007application}
A.~Kanevsky, M.~H. Carpenter, D.~Gottlieb, J.~S. Hesthaven, Application of
  implicit--explicit high--order {R}unge--{K}utta methods to discontinuous
  {G}alerkin schemes, Journal of Computational Physics 225~(2) (2007)
  1753--1781.

\bibitem{ascher1997implicit}
U.~M. Ascher, S.~J. Ruuth, R.~J. Spiteri, Implicit-explicit {R}unge-{K}utta
  methods for time-dependent partial differential equations, Applied Numerical
  Mathematics 25~(2) (1997) 151--167.

\bibitem{Kennedy2003additive}
C.~A. Kennedy, M.~H. Carpenter, Additive {R}unge-{K}utta schemes for
  convection-diffusion-reaction equations, Applied Numerical Mathematics
  44~(1-2) (2003) 139--181.

\bibitem{pareschi2005implicit}
L.~Pareschi, G.~Russo, Implicit-explicit {R}unge-{K}utta schemes and
  applications to hyperbolic systems with relaxation, Journal of Scientific
  computing 25~(1-2) (2005) 129--155.

\bibitem{FeistauerDolejsiKucera07}
M.~Feistauer, V.~Dolejsi, V.~Kucera, On the discontinuous {G}alerkin method for
  the simulation of compressible flow with wide range of {M}ach numbers,
  Computing and Visualization in Science 10 (2007) 17--27.

\bibitem{RestelliGiraldo09}
M.~Restelli, F.~X. Giraldo, A conservative discontinuous {G}alerkin
  semi-implicit formulation for the {N}avier-{S}tokes equations in
  non-hydrostatic mesoscale modeling, SIAM Journal on Scientific Computing 31
  (2009) 2231--2257.

\bibitem{kang2019imex}
S.~Kang, F.~X. Giraldo, T.~Bui-Thanh, Imex hdg-dg: A coupled implicit
  hybridized discontinuous {G}alerkin and explicit discontinuous {G}alerkin
  approach for shallow water systems, Journal of Computational Physics (2019)
  109010.

\bibitem{caliari2004interpolating}
M.~Caliari, M.~Vianello, L.~Bergamaschi, Interpolating discrete
  advection--diffusion propagators at {L}eja sequences, Journal of
  Computational and Applied Mathematics 172~(1) (2004) 79--99.

\bibitem{celledoni2008symmetric}
E.~Celledoni, D.~Cohen, B.~Owren, Symmetric exponential integrators with an
  application to the cubic {S}chr{\"o}dinger equation, Foundations of
  Computational Mathematics 8~(3) (2008) 303--317.

\bibitem{botchev2006gautschi}
M.~A. Botchev, D.~Harutyunyan, J.~J. van~der Vegt, The {G}autschi time stepping
  scheme for edge finite element discretizations of the {M}axwell equations,
  Journal of computational physics 216~(2) (2006) 654--686.

\bibitem{tokman2002three}
M.~Tokman, P.~Bellan, Three-dimensional model of the structure and evolution of
  coronal mass ejections, The Astrophysical Journal 567~(2) (2002) 1202.

\bibitem{li2018exponential}
S.-J. Li, L.-S. Luo, Z.~J. Wang, L.~Ju, An exponential time-integrator scheme
  for steady and unsteady inviscid flows, Journal of Computational Physics 365
  (2018) 206--225.

\bibitem{kooij2018exponential}
G.~L. Kooij, M.~A. Botchev, B.~J. Geurts, An exponential time integrator for
  the incompressible {N}avier--{S}tokes equation, SIAM journal on scientific
  computing 40~(3) (2018) B684--B705.

\bibitem{schulze2009exponential}
J.~C. Schulze, P.~J. Schmid, J.~L. Sesterhenn, Exponential time integration
  using {K}rylov subspaces, International journal for numerical methods in
  fluids 60~(6) (2009) 591--609.

\bibitem{li2019exponential}
S.-J. Li, L.~Ju, Exponential time-marching method for the unsteady
  {N}avier-{S}tokes equations, in: AIAA Scitech 2019 Forum, 2019, p. 0907.

\bibitem{clancy2013use}
C.~Clancy, J.~A. Pudykiewicz, On the use of exponential time integration
  methods in atmospheric models, Tellus A: Dynamic Meteorology and Oceanography
  65~(1) (2013) 20898.

\bibitem{moret2007rd}
I.~Moret, On {RD}-rational {K}rylov approximations to the core-functions of
  exponential integrators, Numerical Linear Algebra with Applications 14~(5)
  (2007) 445--457.

\bibitem{ragni2014rational}
S.~Ragni, Rational {K}rylov methods in exponential integrators for {E}uropean
  option pricing, Numerical Linear Algebra with Applications 21~(4) (2014)
  494--512.

\bibitem{tal2007restart}
H.~Tal-Ezer, On restart and error estimation for {K}rylov approximation of
  $w=f(a)v$, SIAM Journal on Scientific Computing 29~(6) (2007) 2426--2441.

\bibitem{afanasjew2008implementation}
M.~Afanasjew, M.~Eiermann, O.~G. Ernst, S.~G{\"u}ttel, Implementation of a
  restarted {K}rylov subspace method for the evaluation of matrix functions,
  Linear Algebra and its applications 429~(10) (2008) 2293--2314.

\bibitem{botchev2013block}
M.~A. Botchev, A block {K}rylov subspace time-exact solution method for linear
  ordinary differential equation systems, Numerical linear algebra with
  applications 20~(4) (2013) 557--574.

\bibitem{kooij2017block}
G.~L. Kooij, M.~A. Botchev, B.~J. Geurts, A block {K}rylov subspace
  implementation of the time-parallel {P}araexp method and its extension for
  nonlinear partial differential equations, Journal of computational and
  applied mathematics 316 (2017) 229--246.

\bibitem{niesen2012algorithm}
J.~Niesen, W.~M. Wright, Algorithm 919: A {K}rylov subspace algorithm for
  evaluating the $\phi$-functions appearing in exponential integrators, ACM
  Transactions on Mathematical Software (TOMS) 38~(3) (2012) 22.

\bibitem{gaudreault2018kiops}
S.~Gaudreault, G.~Rainwater, M.~Tokman, Kiops: A fast adaptive {K}rylov
  subspace solver for exponential integrators, Journal of Computational Physics
  372 (2018) 236--255.

\bibitem{luan2019further}
V.~T. Luan, J.~A. Pudykiewicz, D.~R. Reynolds, Further development of efficient
  and accurate time integration schemes for meteorological models, Journal of
  Computational Physics 376 (2019) 817--837.

\bibitem{arnoldi1951principle}
W.~E. Arnoldi, The principle of minimized iterations in the solution of the
  matrix eigenvalue problem, Quarterly of applied mathematics 9~(1) (1951)
  17--29.

\bibitem{koskela2015approximating}
A.~Koskela, Approximating the matrix exponential of an advection-diffusion
  operator using the incomplete orthogonalization method, in: Numerical
  Mathematics and Advanced Applications-ENUMATH 2013, Springer, 2015, pp.
  345--353.

\bibitem{vo2017approximating}
H.~D. Vo, R.~B. Sidje, Approximating the large sparse matrix exponential using
  incomplete orthogonalization and {K}rylov subspaces of variable dimension,
  Numerical linear algebra with applications 24~(3) (2017) e2090.

\bibitem{tokman2006efficient}
M.~Tokman, Efficient integration of large stiff systems of {ODE}s with
  exponential propagation iterative ({EPI}) methods, Journal of Computational
  Physics 213~(2) (2006) 748--776.

\bibitem{michels2014exponential}
D.~L. Michels, G.~A. Sobottka, A.~G. Weber, Exponential integrators for stiff
  elastodynamic problems, ACM Transactions on Graphics (TOG) 33~(1) (2014) 7.

\bibitem{michels2017stiffly}
D.~L. Michels, V.~T. Luan, M.~Tokman, A stiffly accurate integrator for
  elastodynamic problems, ACM Transactions on Graphics (TOG) 36~(4) (2017) 116.

\bibitem{giraldo2010high}
F.~Giraldo, M.~Restelli, High-order semi-implicit time-integrators for a
  triangular discontinuous {G}alerkin oceanic shallow water model,
  International journal for numerical methods in fluids 63~(9) (2010)
  1077--1102.

\bibitem{skaflestad2009scaling}
B.~Skaflestad, W.~M. Wright, The scaling and modified squaring method for
  matrix functions related to the exponential, Applied Numerical Mathematics
  59~(3-4) (2009) 783--799.

\bibitem{horn1994topics}
R.~A. Horn, R.~A. Horn, C.~R. Johnson, Topics in matrix analysis, Cambridge
  university press, 1994.

\bibitem{friesner1989method}
R.~A. Friesner, L.~S. Tuckerman, B.~C. Dornblaser, T.~V. Russo, A method for
  exponential propagation of large systems of stiff nonlinear differential
  equations, Journal of Scientific Computing 4~(4) (1989) 327--354.

\bibitem{gallopoulos1992efficient}
E.~Gallopoulos, Y.~Saad, Efficient solution of parabolic equations by {K}rylov
  approximation methods, SIAM Journal on Scientific and Statistical Computing
  13~(5) (1992) 1236--1264.

\bibitem{hochbruck1998exponential}
M.~Hochbruck, C.~Lubich, H.~Selhofer, Exponential integrators for large systems
  of differential equations, SIAM Journal on Scientific Computing 19~(5) (1998)
  1552--1574.

\bibitem{sidje1998expokit}
R.~B. Sidje, Expokit: A software package for computing matrix exponentials, ACM
  Transactions on Mathematical Software (TOMS) 24~(1) (1998) 130--156.

\bibitem{al2009new}
A.~H. Al-Mohy, N.~J. Higham, A new scaling and squaring algorithm for the
  matrix exponential, SIAM Journal on Matrix Analysis and Applications 31~(3)
  (2009) 970--989.

\bibitem{wilcox2010high}
L.~C. Wilcox, G.~Stadler, C.~Burstedde, O.~Ghattas, A high-order discontinuous
  {G}alerkin method for wave propagation through coupled elastic--acoustic
  media, Journal of Computational Physics 229~(24) (2010) 9373--9396.

\bibitem{Roe81}
P.~L. Roe, Approximate {R}iemann solvers, parametric vectors, and difference
  schemes, Journal of Computational Physics 43~(2) (1981) 357--372.

\bibitem{masatsuka2013like}
K.~Masatsuka, I do Like CFD, vol. 1, Vol.~1, Lulu. com, 2013.

\bibitem{guermond2011entropy}
J.-L. Guermond, R.~Pasquetti, B.~Popov, Entropy viscosity method for nonlinear
  conservation laws, Journal of Computational Physics 230~(11) (2011)
  4248--4267.

\bibitem{zingan2013implementation}
V.~Zingan, J.-L. Guermond, J.~Morel, B.~Popov, Implementation of the entropy
  viscosity method with the discontinuous galerkin method, Computer Methods in
  Applied Mechanics and Engineering 253 (2013) 479--490.

\bibitem{gassner2018br1}
G.~J. Gassner, A.~R. Winters, F.~J. Hindenlang, D.~A. Kopriva, The br1 scheme
  is stable for the compressible navier--stokes equations, Journal of
  Scientific Computing 77~(1) (2018) 154--200.

\bibitem{Babuska71}
I.~Babu\v{s}ka, Error bounds for finite element method, Numerische Mathematik
  16 (1971) 322--333.

\bibitem{BabuskaSuri87}
I.~Babu\v{s}ka, M.~Suri, The $hp$-version of the finite element method with
  quasiuniform meshes, Mathematical Modeling and Numerical Analysis 21 (1987)
  199--238.

\bibitem{BabuskaSuri94}
I.~Babu\v{s}ka, M.~Suri, The $p$ and $h$-$p$ version of the finite element
  method, basic principles and properties, SIAM Review 36~(4) (1994) 578--632.

\bibitem{ErnGuermond04}
A.~Ern, J.-L. Guermond, Theory and Practice of Finite Elements, Vol. 159 of
  Applied Mathematical Sciences, Spinger-{V}erlag, 2004.

\bibitem{Ciarlet02}
P.~G. Ciarlet, The finite element method for elliptic problems, Vol.~40 of
  Classics in Applied Mathematics, SIAM (SIAM), Philadelphia, PA, 2002, reprint
  of the 1978 original [North-Holland, Amsterdam; MR0520174 (58 \#25001)].

\bibitem{engel2001one}
K.-J. Engel, R.~Nagel, One-parameter semigroups for linear evolution equations,
  in: Semigroup forum, Vol.~63, Springer, 2001, pp. 278--280.

\bibitem{hochbruck2010exponential}
M.~Hochbruck, A.~Ostermann, Exponential integrators, Acta Numerica 19 (2010)
  209--286.

\bibitem{cockburn1998local}
B.~Cockburn, C.-W. Shu, The local discontinuous galerkin method for
  time-dependent convection-diffusion systems, SIAM Journal on Numerical
  Analysis 35~(6) (1998) 2440--2463.

\bibitem{zhou2003high}
Y.~Zhou, G.~Wei, High resolution conjugate filters for the simulation of flows,
  Journal of Computational Physics 189~(1) (2003) 159--179.

\bibitem{luan2017fourth}
V.~T. Luan, Fourth-order two-stage explicit exponential integrators for
  time-dependent pdes, Applied Numerical Mathematics 112 (2017) 91--103.

\bibitem{drazin2004hydrodynamic}
P.~G. Drazin, W.~H. Reid, Hydrodynamic stability, Cambridge university press,
  2004.

\bibitem{springel2010pur}
V.~Springel, E pur si muove: {G}alilean-invariant cosmological hydrodynamical
  simulations on a moving mesh, Monthly Notices of the Royal Astronomical
  Society 401~(2) (2010) 791--851.

\bibitem{lecoanet2016validated}
D.~Lecoanet, M.~McCourt, E.~Quataert, K.~J. Burns, G.~M. Vasil, J.~S. Oishi,
  B.~P. Brown, J.~M. Stone, R.~M. O'Leary, A validated non-linear
  {K}elvin--{H}elmholtz benchmark for numerical hydrodynamics, Monthly Notices
  of the Royal Astronomical Society 455~(4) (2016) 4274--4288.

\bibitem{liska2003comparison}
R.~Liska, B.~Wendroff, Comparison of several difference schemes on 1d and 2d
  test problems for the euler equations, SIAM Journal on Scientific Computing
  25~(3) (2003) 995--1017.

\end{thebibliography}







\end{document}